\documentclass[11pt]{article}
\usepackage{fullpage}
\usepackage{hyperref}
\usepackage{color}
\hypersetup{
  colorlinks = true,
  citecolor  = blue,
  linkcolor  = red,
  filecolor  = blue,
  urlcolor   = blue,
}

\usepackage{amsmath,amsthm,amssymb,amscd,epsfig,url,authblk}
\usepackage{microtype}
\usepackage{tikz}
\usepackage{float}
\usepackage{subcaption}
\usepackage{mdframed}
 \numberwithin{equation}{section}
 \mdfdefinestyle{box}{skipabove=0, leftmargin=0,%
    rightmargin=0,%
    innertopmargin=-2mm,%
    innerbottommargin=2mm,
    innerleftmargin=2,
    innerrightmargin=2}
  \mdfdefinestyle{box2}{backgroundcolor=gray!10, skipabove=0, leftmargin=0,%
    rightmargin=0,%
    innertopmargin=-2mm,%
    innerbottommargin=2mm,
    innerleftmargin=2,
    innerrightmargin=2}
\newcommand{\pP}{\mathbb{P}}
\newcommand{\E}{\mathbb{E}}
\theoremstyle{plain}
\newtheorem{theorem}{Theorem}[section]
\newtheorem{prp}[theorem]{Proposition}
\newtheorem{lemma}[theorem]{Lemma}
\newtheorem{proposition}[theorem]{Proposition}
\newtheorem{corollary}[theorem]{Corollary}

\newtheorem{example}[theorem]{Example}
\newtheorem{claim}[theorem]{Claim}

\newtheorem*{fact}{Fact}
\theoremstyle{definition}
\newtheorem{definition}[theorem]{Definition}
\newtheorem{remark}[theorem]{Remark}

\usetikzlibrary{shapes,decorations,arrows,positioning,calc}
\tikzstyle{operator}=[circle, radius=0.2cm, text centered, draw=black]
\tikzstyle{arrow}=[thick, ->, >=stealth]
\DeclareMathOperator*{\tr}{Tr}

\DeclareMathOperator*{\ex}{\mathbb{E}}
\newcommand{\norm}[1]{\left\lVert#1\right\rVert}

\newcommand{\OO}[1]{O\left( #1\right)}
\newcommand{\tO}[1]{\widetilde{O}\left( #1 \right)}

\newcommand{\re}{\mathbb{R}}
\newcommand{\na}{\mathbb{N}}
\newcommand{\T}{\mathbf{T}}
\newcommand{\dd}{\text{-}}
\newcommand{\Nc}{\mathcal{N}}

\newcommand{\inp}[2]{\left\langle#1, #2 \right\rangle}

\newcommand{\mm}{M}

\newcommand{\w}[1]{w(#1)}
\newcommand{\numc}[1]{c_{#1}}
\newcommand{\wei}[1]{n_{#1}}
\newcommand{\weii}[1]{x_{#1}}

\newcommand{\plog}{\text{polylog}}

\newcommand{\U}[1]{U_{#1}}
\newcommand{\V}[1]{V_{#1}}
\newcommand{\W}[1]{W_{#1}}
\newcommand{\B}[1]{B(#1)}

\newcommand{\M}{M}
\newcommand{\N}[1]{N_{#1}}
\newcommand{\LL}[1]{L_{#1}}
\newcommand{\val}[1]{\mathrm{val}(#1)}

\newcommand{\uu}[3]{\left( (#1),#2,#3 \right)}
\newcommand{\type}{t}

\newcommand{\ssmin}{s_{\mathrm{min}}}
\newcommand{\smin}{S_{\mathrm{min}}}

\newcommand{\iso}{W_{\mathrm{iso}}}
\newcommand{\isoo}{\mathsf{Iso}}
\newcommand{\tnoniso}{t_{\mathrm{noniso}}}
\newcommand{\ccal}{C_\alpha}
\newcommand{\dal}{D_\alpha}

\newcommand{\erdos}{Erd\H{o}s-R\'{e}nyi}
\newcommand{\nn}{\nu}
\newcommand{\eps}{\epsilon}
\newcommand{\logn}[1]{\log_N(#1)}
\newcommand{\xx}{x}
\newcommand{\cli}{\mathsf{CLIQUE}} 
\newcommand{\tri}{\mathsf{TRIANGLE}} 
\begin{document}
	
	\title{Graph Matrices: Norm Bounds and Applications}

 \author{Kwangjun Ahn} 
 \affil{Massachusetts Institute of Technology\\
Cambridge, MA\\	 \texttt{kjahn@mit.edu}}
\author{Dhruv Medarametla}
\affil{Troy, MI\\
	 \texttt{dhruvm321@gmail.com}} 
	\author{Aaron Potechin}
	\affil{The University of Chicago\\
Chicago, IL\\
	 \texttt{potechin@uchicago.edu}}
	\maketitle
	\begin{abstract}
	  In this paper, we derive nearly tight probabilistic norm bounds for a class of random matrices we call \emph{graph matrices}.
      While the classical case of symmetric matrices with independent random entries
      (Wigner's matrices) is a special case, in general,
      the entries of our matrices will be dependent in a way
      that can be specified in terms of a fixed-size graph we refer to as the \emph{shape}. For Wigner's matrices, this shape is $K_2$, the clique on 2 vertices. To prove our norm bounds, we 
      use the trace power method.

      In a recent series of papers by Potechin and coauthors, graph matrices played a crucial role in proving average-case lower bounds for the Sum-of-Squares (SoS) hierarchy of proof systems, one of the most powerful, but difficult to analyze, techniques
      in combinatorial optimization. In particular, graph matrices played a crucial role in proving that low-degree SoS cannot refute the existence of a large clique in a random graph and proving that low-degree SoS cannot prove a tight lower bound on the ground state energy of the Sherrington-Kirkpatrick Hamiltonian. In this paper, we give several additional applications of graph matrices. We show that for several technical lemmas in the literature, while the original analyses were quite involved, we can give direct proofs using graph matrices and our norm bounds.
		
	\end{abstract}
	
	\section{Introduction}
 
	\subsection{Background on random matrix theory}
	The theory of random matrices is an important area of mathematics with applications in physics, computer science, and many other areas. For random matrices with independent entries, many results are known. Some highlights are as follows:
	\begin{enumerate}
	    \item For symmetric random $n \times n$ matrices whose entries in the upper right triangle are independent random variables with mean $0$ and variance $1$:
	    \begin{enumerate}
	        \item Wigner's Semicircle Law \cite{wigner1958distribution} describes the limit of the mean density of the spectra of these matrices as $n \to \infty$.
	        \item The distribution of the maximum eigenvalue of these matrices is described by the Tracy-Widom distribution \cite{Tracy1994}.
	        \item For the special case of the Gaussian Orthogonal Ensemble where the entries are independent Gaussians, Mehta and Gaudin \cite{MEHTA1960420} found the distribution for how far apart consecutive eigenvalues are from each other.
	    \end{enumerate}
	    \item For random $n \times n$ matrices whose entries are independent random variables with mean $0$ and variance $1$ (where there is no symmetry requirement), Girko's Circular Law \cite{girko} describes the limit of the mean density of the spectra of these matrices as $n \to \infty$.
	\end{enumerate}
	However, for random matrices whose entries are dependent, much less is known.
	\subsection{Properties of graph matrices}
	In this paper, we analyze graph matrices, a class of random matrices which was first explicitly defined in the previous version of this paper \cite{medarametla_et_al:LIPIcs:2016:6663} and which has the following properties:
	\begin{enumerate}
	    \item The entries of the matrix depend on a random input whose size is described by a parameter $n$ and whose entries are labeled by indices in $[n]$ \footnote{In the most general case, the size of the input may be described by several parameters $n_1$, \ldots, $n_{t}$ rather than a single parameter $n$. We discuss this in Section \ref{gendef}.}. For example, this random input could be a random graph on $n$ vertices. \item The dependence of the entries of the matrix on the random input can be described by a small graph which we call a \emph{shape} \footnote{In the most general case, a shape will be a hypergraph with labeled edges rather than a graph. We discuss this in Section \ref{gendef}.}.
	    \item The rows and columns of the matrix are indexed by tuples of indices in $[n]$ and the matrix (as a function of the input) is symmetric under permutations of $[n]$. As a result, if we permute the entries of the input by applying a permutation of $[n]$, this applies a corresponding permutation to the rows and columns of the matrix.
	\end{enumerate}
	Our main results are probabilistic norm bounds on all graph matrices. (See Definitions \ref{def:graphmatrix} and Definition \ref{def:gengraphmatrix} for the formal definitions of graph matrices and see Theorems \ref{thm:mainresult}, \ref{thm:formal:mainresult}, and \ref{boundgeneral} for the statements of the norm bounds).
	The norm bounds we prove are tight up to $\plog(n)$ factors and they are general enough to permit several applications.
	
	\subsection{Applications of graph matrices}

	One important application of graph matrices is analyzing the Sum-of-Squares (SoS) hierarchy of proof systems. The Sum-of-Squares hierarchy, independently investigated by Shor, Nesterov, Parrillo, Lasserre, and Grigoriev \cite{shor1987class,nesterov2000squared,parrilo2000structured,lasserre2001global, grigoriev2001complexity, grigoriev2001linear}, is a powerful tool for solving combinatorial optimization problems as well as various inference problems arising from physics and machine learning.
	SoS is an appealing technique because it has the following nice properties:
	\begin{enumerate}
	    \item Broadly applicable: SoS can be applied to any system of polynomial equations over the real numbers, a problem class which subsumes most combinatorial problems as a special case.
	    \item Surprisingly powerful: SoS captures the Goemans-Williamson algorithm for max-cut \cite{goemans1995improved}, the Goemans-Linial relaxation for sparsest cut (which was shown to give an $O(\sqrt{\log{n}})$ approximation by Arora, Rao, and Vazirani \cite{arora2009expander}), and the subexponential-time algorithm for unique games \cite{arora2010subexponential, barak2011rounding,guruswami2011lasserre}. SoS has also been used to give algorithms for various inference problems arising from physics and machine learning.
	    Examples include planted sparse vector \cite{Barak:2014:RSR:2591796.2591886}, dictionary learning \cite{Barak:2015:DLT:2746539.2746605}, tensor decomposition \cite{ge2015decomposing, hopkins2016fast}, tensor completion \cite{barak2016noisy, potechin2017exact}, quantum separability \cite{Barak:2017:QES:3055399.3055488}, and robust estimation/mixtures of Gaussians \cite{10.1145/3188745.3188748, 10.1145/3188745.3188970,bakshi2020outlier,diakonikolas2020robustly}.
	    \item In some sense, simple: All that SoS uses is polynomial equalities and the fact that squares are non-negative over the real numbers.
	\end{enumerate}
	For more information about the sum of squares hierarchy, see the survey of Barak and Steurer \cite{BarakS14}.
	
	Unfortunately, SoS is often difficult to analyze. One reason for this is that the moment matrix, a matrix which appears in the analysis of SoS, is generally very complicated. For average-case problems, the moment matrices often admit a decomposition into graph matrices, which can be useful for the analysis. Indeed, graph matrices played a crucial role in the SoS lower bounds for planted clique (i.e. proving that low-degree SoS cannot refute the existence of a large clique in a random graph) \cite{Norms,DBLP:journals/corr/HopkinsKP15,DBLP:journals/corr/DeshpandeM15,DBLP:journals/corr/RaghavendraS15, FinalPlantedClique}. Recently, graph matrices were used to prove SoS lower bounds for the Sherrington-Kirkpatrick problem (i.e. proving that low-degree SoS cannot prove a tight lower bound on the ground state energy of the Sherrington-Kirkpatrick Hamiltoniam) \cite{PAPSOS}. 
	
	As we demonstrate in Section \ref{sec:appl}, graph matrices can also capture the proofs of several other technical statements, including the following:
	\begin{enumerate}
	    \item Graph matrices can be used to capture the proof of Theorem 7 of \cite{barak2012hypercontractivity}, which says that if $A$ is an $m\times n$ matrix with each entry drawn independently from $\Nc(0,1)$ then with high probability, $\|A\|_{2 \to 4} \leq 3  m+\tO{ \max\left\{n\sqrt{m},~ n^2\right\}}$ where $\|A\|_{2\to 4} = \max_{\left\| x \right\|_2=1} \norm{Ax}_4$ is the $2 \to 4$ norm\footnote{For background on the $2 \to 4$ norm and the more general $p \to q$ norm $\|A\|_{p \to q} = \max_{\left\| x \right\|_p=1} \norm{Ax}_q$, see the paper \cite{barak2012hypercontractivity}. While we limit our discussion to the $2 \to 4$ norm as it is the most important such norm, we expect that graph matrices can also be used to analyze the $p \to q$ norm for other $p$ and $q$.} of $A$ (see Section~\ref{fournorm} for details). This result is the key idea behind the SoS algorithm for planted sparse vector \cite{Barak:2014:RSR:2591796.2591886}.
	    
	    Moreover, graph matrices can be used to capture the proof of the key lemma for analyzing a faster spectral algorithm for planted sparse vector~\cite{hopkins2016fast};
	    see Section~\ref{fast:sparse} for details.
	    
	    \item Graph matrices can be used to capture the proofs of key lemmas for analyzing the tensor decomposition algorithm of \cite{ge2015decomposing} as well as the faster version in \cite{hopkins2016fast}.
	    See Sections~\ref{overtensor} and \ref{fastovertensor} for details.
	\end{enumerate}
For these applications, our proofs based on graph matrices are \emph{mechanical}, whereas the original proofs for those technical statements are involved and often require clever arguments.

	\subsection{Updates from the previous version of this paper}
	This paper is a major update of the paper ``Bounds on the Norms of Uniform Low Degree Graph Matrices'' which appeared in RANDOM 2016 \cite{medarametla_et_al:LIPIcs:2016:6663}. The graph matrices that our previous paper focused on are those described in Section \ref{simpledef}; these matrices were motivated by the analyses of SoS on the planted clique problem. The work done in our previous paper is covered within the first \ref{pf:main} sections of this paper as well as Appendix \ref{lowerboundsection}.
	
	In this updated paper, we generalize the definition of graph matrices to allow different types of vertices, different input distributions, and hyperedges. We describe these generalizations in Section \ref{gendef} and prove norm bounds on these generalized graph matrices in Section \ref{pf:generalgeneral}. As discussed above, these generalized norm bounds allow us to capture the proofs of several technical statements in the literature; these proofs are contained in Section \ref{sec:appl}. 
	
	\subsection{Notation}
	Given a matrix $M$,  we write $M(i,j)$ to denote the $(i,j)$-th element of $M$; we take $\norm{M}$ to be the spectral norm of $M$, i.e. $\norm{M}=\max_{\norm{v}_2=1} \norm{Mv}_2$. 
	For a graph $G$, we denote its vertex set by $V(G)$ and denote its edge set  by $E(G)$.

    We define $\widetilde{O}(f(n))$ to mean $O(f(n)(\log n)^c)$ for some $c \geq 0$.

	\section{Graph matrices and their norm bounds} \label{simpledef}

    In this section, we will define graph matrices for the setting where the input is an undirected \erdos~random graph.
Before getting into the formal definitions, we begin with some motivating examples.

\subsection{First motivating example: clique indicator} \label{sec:mot1}
Given an undirected graph $G$  on the vertex set $[n]:=\{1,2,\dots, n\}$, we say $H\subset G$ is a clique if  every two vertices in $H$ are adjacent to each other.
Suppose that we are interested in studying the  $n(n-1) \times n(n-1)$ clique indicator matrix $\cli$ defined as follows: 
\begin{align} \label{def:clique}
        \cli\big((i_1,i_2),(j_1,j_2)\big):=\begin{cases} 1, & \text{if}~i_1,i_2,j_1,j_2~\text{are distinct and form a clique in $G$,}\\
    0, &\text{otherwise.}
    \end{cases}
\end{align}
Note that here we require that $i_1 \neq i_2$ and $j_1 \neq j_2$ but $(i_1,i_2)$ and $(j_1,j_2)$ can intersect.

 More specifically, we consider the case where $G$ is a random graph generated as per the \erdos~model: for each pair of vertices,  there is an edge between them with probability $1/2$ and all of these events are independent. In this case, the clique indicator $\cli$ is a random matrix but its entries are not independent of each other. 

One natural way to study the clique indicator is via discrete Fourier analysis.    
We encode the edges of $G$ as the $\pm 1$ variables  $\{\chi_{\{l_1,l_2\}}\}$ where $\chi_{\{l_1,l_2\}}:=+1$ $\text{if}~\{l_1,l_2\} \in E(G)$, and $-1$ otherwise.
We can now decompose each entry of $\cli$ as a sum of products of the variables as follows:
\begin{fact} 
The following identity holds for any four distinct indices $i_1,i_2,j_1,j_2 \in [n]$:
\begin{align}\label{mot:decomp}
    \cli \big((i_1,i_2),(j_1,j_2)\big) = \frac{1}{2^6}\cdot \sum_{\substack{R:~\text{graph on} ~ \{i_1,i_2,j_1,j_2 \}}} \prod_{\{ l_1,l_2  \}\in E(R)}\chi_{ \{ l_1,l_2 \} }\,.
\end{align}
\end{fact}
\begin{proof}
If the four vertices $i_1,i_2,j_1,j_2$ form a clique, each summand $\prod_{\{l_1,l_2\}\in E(R)}\chi_{ \{ l_1,l_2 \} }$ is equal to $1$, which implies that the right hand side of \eqref{mot:decomp} is $1$.
Otherwise, some edge $e$ between $i_1,i_2,j_1,j_2$ is missing. In this case, there are perfect cancellations between the graphs $R$ where $e \in E(R)$ and the graphs $R$ where $e \notin E(R)$ so the right hand side is equal to $0$.
\end{proof}

Now let us inspect the terms in the decomposition \eqref{mot:decomp}.
For concreteness, consider the case $i_1=1,~i_2=2,~j_1=3,~j_2=4$:
\begin{align}\label{mot:decomp:1234}
    \cli \big((1,2),(3,4)\big) = \frac{1}{2^6}\cdot \sum_{\substack{R:~\textnormal{graph on} ~ \{1,2,3,4 \}}} \prod_{\{ l_1,l_2  \}\in E(R)}\chi_{ \{ l_1,l_2 \} }\,.
\end{align}
Notably, each term in the decomposition corresponds to a graph $R$ on four vertices $1,2,3,4$. Let us call each such graph $R$ a ribbon (see Definition~\ref{def:ribbon} for precise details). To represent the fact that $R$ corresponds to row $(1,2)$ and column $(3,4)$ of $\cli$, we associate tuples $A_R = (1,2)$ and $B_R = (3,4)$ to $R$ which we call the left and right vertices of $R$.

For example, one term in \eqref{mot:decomp} is $\chi_{\{1,3\}}\chi_{\{2,3\}}\chi_{\{2,4\}}$. We represent this term with the ribbon $R$ where  $A_R = (1,2)$, $B_R = (3,4)$, and $E(R) = \{\{1,3\},\{2,3\},\{2,4\}\}$. See Figure~\ref{fig:ribbon:ex} for an illustration.
	
	\begin{figure}[H]
\centering
				\scalebox{0.65}{
					\begin{tikzpicture}[shorten >=1pt,auto,node distance=3cm,
					thick,main node/.style={circle,fill=blue!20,draw,font=\sffamily\LARGE\bfseries}]
					\node[main node] (1) at (-3,3)  {$1$};
					\node[main node] (2) at (-3,0)  {$2$};
					\node[main node] (3) at (3,3)  {$3$};
					\node[main node] (4) at (3,0)  {$4$};  
					\node[draw=none,fill=none] (8) at (-3,4.5) {\fontsize{20}{22.4}\selectfont$A_R$};
					 
					\node[draw=none,fill=none] (10) at (3,4.5) {\fontsize{20}{22.4}\selectfont$B_R$};
						\path[every node/.style={font=\sffamily},line width=1pt]
				  (1) edge node [] {} (3)
				    (2) edge node [] {} (3)
				    (2) edge node [] {} (4);
				\draw[draw=black] ($(1.135)+(-0.3,0.3)$) rectangle  ($(2.-45)+(0.3,-0.3)$);
		\draw[draw=black] ($(3.135)+(-0.3,0.3)$) rectangle  ($(4.-45)+(0.3,-0.3)$);		\end{tikzpicture}
				}
				\caption{An illustration of the ribbon $R$ where the vertex set $V(R)$ is $\{1,2,3,4\}$, $A_R=(1,2)$, $B_R=(3,4)$, and $E(R) = \{\{1,3\},\{2,3\},\{2,4\}\}$} \label{fig:ribbon:ex}
				\end{figure}
 
 Now let us represent \eqref{mot:decomp} using ribbons.
 To that end, we define a Fourier character for each ribbon $R$:
 \begin{align*}
     \chi_{R}:= \prod_{\{l_1,l_2\} \in E(R)} \chi_{\{l_1,l_2\}}\,.
 \end{align*}
Expressing the terms in \eqref{mot:decomp} using these ribbons and their Fourier characters, we have that
\begin{align} \label{new:decom}
       \cli \big((1,2),(3,4)\big) = \frac{1}{2^6}\cdot \sum_{\substack{\text{ribbon } R: V(R) = \{1,2,3,4\}\\  A_R = (1,2),~B_R = (3,4) }} \chi_R\,.
\end{align}

In fact, there is nothing special about the choice $i_1=1,~i_2=2,~j_1=3,~j_2=4$, and \eqref{new:decom} is true for any choice of $i_1,i_2,j_1,j_2$ as long as they are distinct.
This motivates us to consider the shape $\alpha$ of a ribbon $R$: to obtain the shape $\alpha$ of a ribbon $R$, we keep the graph structure of $R$ the same while replacing the vertices $1,2,3,4$ by (undetermined) variables $u_1,u_2,v_1,v_2$, respectively. 
Moreover, we replace $A_R$ and $B_R$ with $\U{\alpha}:=(u_1,u_2)$ and $\V{\alpha}:=(v_1,v_2)$, respectively.
See Definition~\ref{def:shape} for precise details.

For example, if we consider the ribbon in Figure~\ref{fig:ribbon:ex}, its shape $\alpha$ has $V(\alpha) = \{u_1,u_2,v_1,v_2\}$, $\U{\alpha} = (u_1,u_2)$, $\V{\alpha} = (v_1,v_2)$, and $E(\alpha) = \{\{u_1,v_1\},\{u_2,v_1\},\{u_2,v_2\}\}$.  See Figure~\ref{fig:shape:ex} for an illustration.

	\begin{figure}[H]
\centering
				\scalebox{0.65}{
					\begin{tikzpicture}[shorten >=1pt,auto,node distance=3cm,
					thick,main node/.style={circle,fill=blue!20,draw,font=\sffamily\LARGE\bfseries}]
					\node[main node] (1) at (-3,3)  {$u_1$};
					\node[main node] (2) at (-3,0)  {$u_2$};
					\node[main node] (3) at (3,3)  {$v_1$};
					\node[main node] (4) at (3,0)  {$v_2$};  
					\node[draw=none,fill=none] (8) at (-3,4.5) {\fontsize{20}{22.4}\selectfont$\U{\alpha}$};
					 
					\node[draw=none,fill=none] (10) at (3,4.5) {\fontsize{20}{22.4}\selectfont$\V{\alpha}$};
						\path[every node/.style={font=\sffamily},line width=1pt]
				  (1) edge node [] {} (3)
				    (2) edge node [] {} (3)
				    (2) edge node [] {} (4);
				\draw[draw=black] ($(1.135)+(-0.3,0.3)$) rectangle  ($(2.-45)+(0.3,-0.3)$);
		\draw[draw=black] ($(3.135)+(-0.3,0.3)$) rectangle  ($(4.-45)+(0.3,-0.3)$);		\end{tikzpicture}
				}
				\caption{An illustration of the shape $\alpha$ where the vertex set $V(\alpha)$ is $\{u_1,u_2,v_1,v_2\}$, $\U{\alpha}=(u_1,v_1)$, $\V{\alpha}=(v_1,v_2)$, and $E(\alpha) = \{\{u_1,v_1\},\{u_2,v_1\},\{u_2,v_2\}\}$.}
				\label{fig:shape:ex}
				\end{figure}

Based on shapes, we can write a matrix analogue of the decomposition  \eqref{new:decom}.
Given a shape $\alpha$, we group all of the ribbons with a given shape $\alpha$ together and then construct a corresponding matrix $M_{\alpha}$:
\begin{align*}
\M_{\alpha}\big((i_1,i_2),(j_1,j_2)\big) :=\begin{cases} \chi_{R_{i_1,i_2,j_1,j_2}}  & \text{if}~i_1,i_2,j_1,j_2~\text{are distinct,}\\
    0, &\text{otherwise,} 
    \end{cases}
\end{align*}
where $R_{i_1,i_2,j_1,j_2}$ is the ribbon with shape $\alpha$ such that  $A_R=(i_1,i_2)$ and $B_R= (j_1,j_2)$.
We call this matrix $\M_{\alpha}$ a graph matrix  (see Definition~\ref{def:graphmatrix} for  precise details).  
Using these graph matrices, we have the following matrix analogue of \eqref{new:decom}:
\begin{align}
\cli = \frac{1}{2^6}\cdot \sum_{\substack{\text{shape }\alpha: V(\alpha) = \{u_1,u_2,v_1,v_2\},\\
\U{\alpha} = (u_1,u_2), ~\V{\alpha} = (v_1,v_2)}}{M_{\alpha}}\,.
\end{align}
As we discuss in Example \ref{revisitingcliqueindicatorexample}, this decomposition of $\cli$ into graph matrices allows us to analyze it by separating out the parts of $\cli$ which have large norm from the parts of $\cli$ which have small norm.
Now let us take a look at another example.
\begin{remark}
For the particular Z-shaped $\alpha$ shown in Figure~\ref{fig:shape:ex}, the corresponding graph matrix $M_{\alpha}$ has been analyzed more precisely. In particular, Cai and Potechin \cite{CP20} determined the limit of the mean density of the spectrum of the singular values of this matrix.
\end{remark}
\subsection{Second motivating example: counting triangles}
\label{sec:mot2}
Under the same random graph model, suppose that now we are interested in counting the number of triangles that include two chosen vertices.  
More specifically, for any two distinct vertices $i,j\in[n]$, define $\triangle_{i,j}$ be the number of triangles containing $i,j$, and consider the  $n\times n $ matrix $\tri$  defined as follows: 
\begin{align} \label{def:tri}
        \tri\big(i_1,j_1\big):=\begin{cases} \triangle_{i_1,j_1}, & \text{if}~i_1,j_1~\text{are distinct},\\
    0, &\text{otherwise.}
    \end{cases}
\end{align}
Again, we can study $\tri$ via discrete Fourier analysis. Similar to \eqref{mot:decomp},
we have the following decomposition:
\begin{align}\label{mot2:decomp}
    \tri \big(i_1, j_1\big) = \frac{1}{2^3}\cdot  \sum_{k_1\in [n]\setminus\{i_1,j_1\}}\sum_{\substack{R:~\text{graph on} ~ \{i_1,j_1,k_1 \}}} \prod_{\{ l_1,l_2  \}\in E(R)}\chi_{ \{ l_1,l_2 \} }\,.
\end{align} 

After some inspection of \eqref{mot2:decomp}, one will realize that the terms in the decomposition cannot be characterized using the ribbons and shapes we described in Section~\ref{sec:mot1}. This is because there is a vertex  $k_1$ in a graph $R$ that does not come from the row/column indices $i_1,j_1$. To capture this, we introduce \emph{middle} vertices in ribbons and shapes.

To illustrate, let us consider the case $i_1=1$, $j_1=2$ and $k_1=3$.
In this case, for each ribbon in \eqref{mot2:decomp}, we can set $A_R = (1)$ and $B_R = (2)$ but there is still a vertex $3$ which is outside of $A_R$ and $B_R$.
We will call this vertex a middle vertex and denote the set of middle vertices by $C_R$ (note that in general there could be more than one middle vertex). Here we have $C_R=\{3\}$. See Definition~\ref{def:ribbon} for precise details.

For example, one term in \eqref{mot2:decomp} is $\chi_{\{1,3\}}\chi_{\{2,3\}} $. We represent this term with the ribbon $R$ where  $A_R = (1)$, $B_R = (2)$, $C_R = \{3\}$ and $E(R) = \{\{1,3\},\{2,3\} \}$. See Figure~\ref{fig:ribbon:ex2} for an illustration.
	
	\begin{figure}[H]
\centering
				\scalebox{0.65}{
					\begin{tikzpicture}[shorten >=1pt,auto,node distance=3cm,
					thick,main node/.style={circle,fill=blue!20,draw,font=\sffamily\LARGE\bfseries}]
					\node[main node] (1) at (-3,3)  {$1$};
					\node[main node] (2) at (3,3)  {$2$};
					 	\node[main node] (3) at (0,3)  {$3$}; 
					\node[draw=none,fill=none] (8) at (-3,4.5) {\fontsize{20}{22.4}\selectfont$A_R$};
					 
					\node[draw=none,fill=none] (9) at (0,4.5) {\fontsize{20}{22.4}\selectfont$C_R$};
				 
					\node[draw=none,fill=none] (10) at (3,4.5) {\fontsize{20}{22.4}\selectfont$B_R$};
						\path[every node/.style={font=\sffamily},line width=1pt]
				  (1) edge node [] {} (3)
				    (2) edge node [] {} (3);
				\draw[draw=black] ($(1.135)+(-0.3,0.3)$) rectangle  ($(1.-45)+(0.3,-0.3)$);
		\draw[draw=black] ($(2.135)+(-0.3,0.3)$) rectangle  ($(2.-45)+(0.3,-0.3)$);		\end{tikzpicture}
				}
				\caption{An illustration of the ribbon $R$ where the vertex set $V(R)$ is $\{1,2,3\}$, $A_R=(1)$, $B_R=(2)$, $C_R=\{3\}$, and $E(R) = \{\{1,3\},\{2,3\} \}$.} \label{fig:ribbon:ex2}
				\end{figure}

Following Section~\ref{sec:mot1}, we can now similarly abstract away the identity of the indices $i_1=1,~j_1=2$, and $k_1=3$ and consider the shape of a ribbon. 
To obtain the shape $\alpha$ of a ribbon $R$, we keep the graph structure of $R$ the same while replacing the vertices $1,2,3$ by (undetermined) variables $u_1,v_1,w_1$, respectively. 
Moreover, we replace $A_R$, $B_R$, and $C_R$ with $\U{\alpha}:=(u_1)$, $\V{\alpha}:=(v_1)$, and $\W{\alpha}:=\{w_1\}$ respectively.
See Definition~\ref{def:shape} for precise details.

For example, if we consider the ribbon in Figure~\ref{fig:ribbon:ex2}, its shape $\alpha$ has $V(\alpha) = \{u_1,v_1,w_1\}$, $\U{\alpha} = (u_1)$, $\V{\alpha} = (v_1)$, $\W{\alpha} = \{w_1\}$, and $E(\alpha) = \{\{u_1,w_1\},\{v_1,w_1\} \}$.  See Figure~\ref{fig:shape:ex2} for an illustration.

	\begin{figure}[H]
\centering
				\scalebox{0.65}{
					\begin{tikzpicture}[shorten >=1pt,auto,node distance=3cm,
					thick,main node/.style={circle,fill=blue!20,draw,font=\sffamily\LARGE\bfseries}]
					\node[main node] (1) at (-3,3)  {$u_1$};
					\node[main node] (2) at (3,3)  {$v_1$};
					 	\node[main node] (3) at (0,3)  {$w_1$}; 
					\node[draw=none,fill=none] (8) at (-3,4.5) {\fontsize{20}{22.4}\selectfont$\U{\alpha}$};
					 
					\node[draw=none,fill=none] (9) at (0,4.5) {\fontsize{20}{22.4}\selectfont$\W{\alpha}$};
				 
					\node[draw=none,fill=none] (10) at (3,4.5) {\fontsize{20}{22.4}\selectfont$\V{\alpha}$};
						\path[every node/.style={font=\sffamily},line width=1pt]
				  (1) edge node [] {} (3)
				    (2) edge node [] {} (3);
				\draw[draw=black] ($(1.135)+(-0.3,0.3)$) rectangle  ($(1.-45)+(0.3,-0.3)$);
		\draw[draw=black] ($(2.135)+(-0.3,0.3)$) rectangle  ($(2.-45)+(0.3,-0.3)$);		\end{tikzpicture}
				}
				\caption{An illustration of a shape $\alpha$ where the vertex set $V(\alpha)$ is $\{u_1,v_1,w_1\}$, $\U{\alpha}=(u_1)$, $\V{\alpha}=(v_1)$, $\W{\alpha}=\{w_1\}$ and $E(\alpha) = \{\{u_1,w_1\},\{v_1,w_1\} \}$.} \label{fig:shape:ex2}
				\end{figure}

As in Section~\ref{sec:mot1}, we can write a matrix analogue of the decomposition \eqref{mot2:decomp} using shapes.
However, since there is a middle vertex, there are many different ribbons which have a given shape $\alpha$, row index $i_1$, and column index $j_1$. To handle this, we need to take the sum of the Fourier characters for these ribbons (see Definition~\ref{def:graphmatrix} for precise details). This gives us the following expression: 
\begin{align*}
\M_{\alpha}\big(i_1, j_1\big) :=\begin{cases} \displaystyle \sum_{k_1\in [n] \setminus \{i_1,j_1\} }{\chi_{R_{i_1,j_1,k_1}}},   & \text{if}~i_1,j_1~\text{are distinct,}\\
    0, &\text{otherwise.} 
    \end{cases}
\end{align*} 
where $R_{i_1,j_1,k_1}$ is the ribbon with shape $\alpha$ such that  $A_R=(i_1)$, $B_R= (j_1)$, and $C_R = \{k_1\}$.  
Using these graph matrices, we can decompose $\tri$ as follows:
\begin{align}
\tri = \frac{1}{2^3}\cdot \sum_{\substack{\text{shape }\alpha: V(\alpha) = \{u_1,v_1,w_1\},\\
\U{\alpha} = (u_1), ~\V{\alpha} = (v_1),~\W{\alpha} =\{w_1\}   }}{M_{\alpha}}\,.
\end{align}

As this example shows, by introducing middle vertices into graph matrices, we can increase their expressive power. 
We now give formal definitions of the concepts presented above.

\subsection{Definition of graph matrices}  \label{subsec:simple}

We consider graph matrices to be matrices whose rows and columns are indexed by monomials.
We begin by considering the ground set where the indices of the variables are taken from.
\begin{definition}[Ground set] \label{simp:ground}
We set the ground set for the variable indices to be $[n]:=\{1,2,\cdots, n\}$.
In other words, we consider monomials consisting of the variables $\{x_i\}_{i\in [n] }$.

\end{definition}

For the simplified setting in this section, we focus on the case where the constraint $x^2_i = 1$ is present for each $i$ so we have that $x_i = \pm{1}$.
In this case, all monomials can be reduced to multilinear monomials.
Since we can represent each multilinear monomial $\prod_{a \in A}{x_a}$ by the set of indices $A$, we make the following definition about matrix indices.
	\begin{definition}[Matrix indices] \label{simp:indices}
We define a matrix index $A$ to be a tuple of distinct indices $A = (a_1,\ldots,a_{m})$ for some $m\in \na$, where $a_1,\dots, a_{m} \in [n]$. 
	We denote by $|A|$ the size of the index tuple, i.e., $|A|=m$.
	Moreover, we let $V(A):=\{a_1,a_2,\dots,a_{m}\}$. 
	\end{definition}
\begin{remark}
Note that there are $n*(n-1)*\ldots*(n-m+1) = \frac{n!}{(n-m)!}$ different matrix indices $A$ of size $m$. Thus, if a graph matrix has rows indexed by monomials of degree $m$ and columns indexed by monomials of degree $m'$ then it will be an $\frac{n!}{(n-m)!} \times \frac{n!}{(n-m')!}$ matrix.
\end{remark}
\begin{remark} \label{rmk:why}
A careful reader might wonder why we chose tuples instead of subsets to represent monomials.
For instance, matrix indices $(1,2,3)$ and $(1,3,2)$ both represent the monomial $x_1x_2x_3$, and hence, it may seem redundant to consider tuples over subsets.
However, it turns out that adopting tuples in lieu of subsets results in a more canonical definition of graph matrices. In particular, if we used sets for the matrix indices rather than tuples then some graph matrices would not be symmetric under permutations of $[n]$ but using tuples ensures that all graph matrices are symmetric under permutations of $[n]$.
\end{remark}
Given matrix indices, we now discuss how each entry is defined. We first formally define the input random graph:
\begin{definition}[\erdos~input distribution]
    We define the \erdos~input distribution to be the random graph $G$ having $[n]$ as its vertices such that for each pair of distinct vertices $i,j\in [n]$, the edge $e=\{i,j\}$ is present with probability $1/2$ and all of these events are independent. 
    We denote the \erdos-input distribution by  $G(n,1/2)$.
\end{definition}
Having defined the input distribution, we now formalize how the entries of graph matrices depend on this distribution.
To that end, we first define Fourier characters over input graphs through the following two definitions:
	\begin{definition}
	For a given graph $G$ on $[n]$ and a pair of distinct indices $e=\{i,j\}\subset [n]$,  we define $\chi_{e}(G)$ to be $1$ if $e \in E(G)$ and $-1$ if $e\notin E(G)$.
	\end{definition}
	\begin{definition}[Fourier characters]
		For a given graph $G$ on $[n]$ and a set $E=\{e_1,e_2,\dots, e_l\}$ consisting of pairs of distinct indices in $[n]$, we define the Fourier character $\chi_E(G) = \prod_{i=1}^l{\chi_{e_i}(G)} = (-1)^{|E \setminus E(G)|}$.
	\end{definition}

With the Fourier characters we have defined over the input $G\sim G(n,1/2)$, the $(A,B)$-th entry of a graph matrix is defined as a Fourier character on an edge set that depends on $A$ and $B$.
To be more specific about how the edge set depends on $A$ and $B$, we now formally define \emph{ribbons}:
\begin{mdframed}[style=box]\begin{definition}[Ribbons] \label{def:ribbon}
	A ribbon $R = (A_R,B_R,C_R,E(R))$  consists of two matrix indices $A_R$ and $B_R$ (possibly $V(A_R)\cap V(B_R)\neq\emptyset$), an additional set of indices $C_R$ such that $C_R \cap (V(A_R) \cup V(B_R)) = \emptyset$, and a set $E(R)$  of pairs of distinct indices from $V(A_R) \cup V(B_R)\cup  C_R$. We make the following definitions about ribbons:
		\begin{enumerate}
			\item (Graphical representation of a ribbon) We identify $R$ with its graphical representation defined as a graph with vertices $V(R) = V(A_R) \cup V(B_R) \cup C_R$ and edges $E(R)$.
			Here we regard $V(A_R)$ and $V(B_R)$ as distinguished sets, and refer to them as the left and right vertices of $R$, respectively. We refer to $C_R$ as the middle vertices of $R$.
				
			\item (Fourier character of a ribbon) We define the Fourier character $\chi_{R}$ to be $\chi_{E(R)}$.
		
			\item (Matrix associated to a ribbon) Given a ribbon $R$ and an input graph $G\sim G(n,1/2)$, we define $M_R$  to be a $\frac{n!}{(n-|A_R|)!}\times \frac{n!}{(n-|B_R|)!}$ matrix such that $M_R(A_R,B_R) = \chi_{R}(G)$ and $M_R(A,B) = 0$ if $A \neq A_R$ or $B \neq B_R$.
		\end{enumerate}
	\end{definition}
	\end{mdframed}
See Appendix~\ref{app:example} for examples of ribbons.

Although  Definition~\ref{def:ribbon} looks complicated at first sight, ribbons are indeed natural objects to consider when we consider the monomial-indexed matrices defined over the input distribution $G(n,1/2)$:
	
	\begin{proposition}[Ribbons as orthonormal basis for  matrix functions] \label{prop:ortho}
	Given $a,b\in \na$, the collection of matrices $M_R$ such that $|A_R|=a$, $|B_R|=b$, and $C_R$ does not contain any isolated vertices forms an orthonormal basis for the space of matrix-valued functions such that:
		\begin{enumerate}
		    \item Each function maps an input graph $G \sim G(n,1/2)$ to an $\frac{n!}{(n-a)!}\times \frac{n!}{(n-b)!}$ matrix whose rows and columns are indexed by matrix indices $A$ and $B$ of size $a$ and $b$, respectively.
		    \item We take the inner product $\langle{M,M'}\rangle = \E_{G \sim G(n,\frac{1}{2})}\left[\sum_{A,B}{M(A,B)M'(A,B)}\right]$.
		\end{enumerate}
	\end{proposition}
	We are now ready to define graph matrices.
	Informally, given the sizes $a$ and $b$ of the row and column matrix index, respectively, a graph matrix of a certain ``shape'' is an $\frac{n!}{(n-a)!}\times \frac{n!}{(n-b)!}$ matrix which is the sum of $M_R$ over all ribbons $R$ whose graphical representation has the designated ``shape.''
	To formalize this, we make the following definitions about shapes:
	\begin{definition}[Index shapes] \label{simp:indexshape}
		We define an index shape $U$ to be a tuple of distinct variables $U = (u_1,\ldots,u_{m})$ for some $m\in \na$.
		We denote by $|U|$ the size of the index shape, i.e., $|U|=m$.
	Moreover, let  $V(U):=\{u_1,u_2,\cdots,u_{m}\}$.
	\end{definition}
	\begin{mdframed}[style=box]
	\begin{definition}[Shapes] \label{def:shape}
		A shape $\alpha = (\U{\alpha},\V{\alpha},\W{\alpha},E(\alpha))$  consists of index shapes $\U{\alpha}$ and $\V{\alpha}$ (possibly having common variables), an additional set   $\W{\alpha}$ of variables distinct from  the ones in $V(\U{\alpha})\cup V(\V{\alpha})$, and a set $E(\alpha)$ consisting of pairs of distinct variables from $V(\U{\alpha}) \cup V(\V{\alpha})\cup \W{\alpha}$.
		We use $u$, $v$, and $w$ to denote variables in $\U{\alpha}$, $\V{\alpha}$, and $\W{\alpha}$, respectively.
		
			Furthermore, we identify $\alpha$ with its graphical representation defined as a graph with vertices $V(\alpha)=V(\U{\alpha}) \cup V(\V{\alpha}) \cup \W{\alpha}$ and edges $E(\alpha)$.
		Here we regard $V(\U{\alpha})$ and $V(\V{\alpha})$ as distinguished sets and refer to them as the left and right vertices of $\alpha$, respectively. We refer to $\W{\alpha}$ as the middle vertices of $\alpha$.
		
		Lastly, we sometimes use the notation $u\in \alpha$ in place of $u\in V(\alpha)$ to avoid confusion between $V_{\alpha}$ and $V(\alpha)$.
	\end{definition}
	\end{mdframed}
	See Appendix~\ref{app:example} for examples of shapes.

	With Definition~\ref{def:shape}, we can now be formal about what it means for a ribbon $R$ to have a shape $\alpha$. 
	To do that, we first define realizations of variables and shapes:
	\begin{definition}[Realizations of variables] 
		Given a set of variables $U$, we say a map $\sigma: U \to [n]$ is a realization of $U$ if $\sigma$ an injective map. Moreover, for a realization $\sigma$, we make the following definitions:
		\begin{enumerate}
		    \item For an index shape $U$, 
		We define $\sigma(U)$ to be the matrix index obtained by applying $\sigma$ to each variable of $U$, i.e., if $U = (u_1,\ldots,u_{|U|})$, then $\sigma(U) = (\sigma(u_1),\ldots,\sigma(u_{|U|}))$.
	    \item Moreover, for  a set of pairs of distinct variables $E$, we define $\sigma(E)$ to be the set of pairs of distinct indices obtained from $E$ by applying $\sigma$ to each pair. 
		\end{enumerate}
	\end{definition}
	\begin{definition}[Realizations of shapes]
		Given a shape $\alpha$ and a realization $\sigma: V(\alpha) \to [n]$, we define $\sigma(\alpha)$ to be the ribbon $\sigma(\alpha) = (\sigma(\U{\alpha}), \sigma(\V{\alpha}), \sigma(\W{\alpha}), \sigma(E(\alpha)))$.
	\end{definition}
	See Appendix~\ref{app:example} for examples of realizations of shapes. 
	With these definitions, we can formally define shapes of ribbons:
	\begin{definition}[Shapes of ribbons]
	Given a shape $\alpha$ and a ribbon $R$, we say $R$ has shape $\alpha$ if there exists a realization $\sigma$ such that $\sigma(\alpha)=R$.
	\end{definition}
	We are now ready to define graph matrices:
\begin{mdframed}[style=box]
	\begin{definition}[Graph matrices] \label{def:graphmatrix}
		Given a shape $\alpha$, we define the graph matrix $M_{\alpha}$ to be the $\frac{n!}{(n-|\U{\alpha}|)!}\times \frac{n!}{(n-|\V{\alpha}|)!}$ matrix
	such that for matrix indices $A$ and $B$ with $|A|=|\U{\alpha}|$ and $|B|= |\V{\alpha}|$,
		\begin{align} \label{def:1}
		M_{\alpha}(A,B) = \sum_{\substack{\sigma \text{ is a realization of } \alpha,\\~\sigma(\U{\alpha}) = A,~ \sigma(\V{\alpha})=B}}{\chi_{\sigma(E(\alpha))}}\,.
		\end{align}
	\end{definition}
\end{mdframed}

	\begin{remark} \label{rmk:altdef1}
    It is sometimes more convenient to work with the following alternative definition of graph matrices:
		\begin{align} \label{def:2}
		M_{\alpha} = \sum_{R:~ R \text{ has shape $\alpha$}}{M_{R}}\,.
		\end{align}
	This alternative definition is convenient when we are decomposing a matrix $M$ where we have the Fourier coefficient $c_R$ for each ribbon $R$ (which is generally the case when we are using the pseudo-calibration technique of \cite{FinalPlantedClique} to analyze the Sum-of-Squares hierarchy). In this case, letting $\alpha$ be the shape of $R$, the coefficient for $M_{\alpha}$ is just $c_R$. This would not be the case for our definition.
	
	Fortunately, the two definitions~\eqref{def:1} and \eqref{def:2} only differ by a constant multiplicative factor.
	More precisely, let $\mathrm{Aut}(\alpha)$ be the set of graph automorphisms on $\alpha$ which keep $U_{\alpha}$ and $V_{\alpha}$ fixed. 
	One can verify that $\eqref{def:1} =  |\mathrm{Aut}(\alpha)|\cdot \eqref{def:2}$.
	Thus, the norm bounds which we prove in this paper also apply to this alternative definition and are in fact slightly stronger.
	\end{remark}
\subsection{Examples of graph matrices}\label{sec:ex}
In this section, we provide several examples of graph matrices. Our first example is the $\pm 1$ Wigner matrix:
	\begin{example}[$\pm 1$ Wigner matrix] \label{ex:3}
	Consider $V(\V{\alpha})=\{u_1,v_1\}$, $\U{\alpha}=(u_1)$, $\V{\alpha} =(v_1)$, and $E(\alpha)=\big\{\{u_1,v_1\}\big\}$. 
	In other words, the graphical representation of $\alpha$ is the graph with two vertices $u_1,v_1$ connected by the edge $\{u_1,v_1\}$; see Figure~\ref{fig:3a}. 
	In this case, $\M_\alpha$ is an $n\times n$ symmetric random matrix such that each off-diagonal entry is $\pm 1$ with probability $1/2$ and the diagonal entries are zeros.
	\begin{figure}[H]
				\centering
				\scalebox{0.65}{
					\begin{tikzpicture}[shorten >=1pt,auto,node distance=3cm,
					thick,main node/.style={circle,fill=blue!20,draw,font=\sffamily\LARGE\bfseries}]
					\node[main node] (1) at (-1.5,0)  {$u_1$};
					\node[main node] (2) at (1.5,0)  {$v_1$};
					\node[draw=none,fill=none] (3) at (-1.5,1.5) {\fontsize{20}{22.4}\selectfont$\U{\alpha}$};
					\node[draw=none,fill=none] (4) at (1.5,1.5) {\fontsize{20}{22.4}\selectfont$\V{\alpha}$};
					\path[every node/.style={font=\sffamily},line width=1pt]
				    (1) edge node [] {} (2);
		\draw[draw=black] ($(1.135)+(-0.3,0.3)$) rectangle  ($(1.-45)+(0.3,-0.3)$);
		\draw[draw=black] ($(2.135)+(-0.3,0.3)$) rectangle  ($(2.-45)+(0.3,-0.3)$);
	\end{tikzpicture}
					
				}
				\caption{ The shape corresponding to the $\pm 1$ Wigner matrix from Example~\ref{ex:3}.}
				\label{fig:3a}
			\end{figure}

	\end{example}
	
Our second example is the matrix $J-I$, where $J$ refers to the all $1$'s matrix. Note that this matrix is deterministic as $E(\alpha)$ is empty.

    \begin{example}[Off-diagonal $1$'s matrix] \label{ex:off-diag}
	Consider $V(\V{\alpha})=\{u_1,v_1\}$, $\U{\alpha}=(u_1)$, $\V{\alpha} =(v_1)$, and $E(\alpha)=\big\{\big\}$. 
	In other words, the graphical representation of $\alpha$ is the empty graph with vertices $u_1$ and $v_1$; see Figure~\ref{fig:4}. 
	In this case, $\M_\alpha$ is an $n\times n$ symmetric matrix with $0$'s on the diagonals and $1$'s off the diagonals.
	\begin{figure}[H]
				\centering
				\scalebox{0.65}{
					\begin{tikzpicture}[shorten >=1pt,auto,node distance=3cm,
					thick,main node/.style={circle,fill=blue!20,draw,font=\sffamily\LARGE\bfseries}]
					\node[main node] (1) at (-1.5,0)  {$u_1$};
					\node[main node] (2) at (1.5,0)  {$v_1$};
					\node[draw=none,fill=none] (3) at (-1.5,1.5) {\fontsize{20}{22.4}\selectfont$\U{\alpha}$};
					\node[draw=none,fill=none] (4) at (1.5,1.5) {\fontsize{20}{22.4}\selectfont$\V{\alpha}$};
					\path[every node/.style={font=\sffamily},line width=1pt];
		\draw[draw=black] ($(1.135)+(-0.3,0.3)$) rectangle  ($(1.-45)+(0.3,-0.3)$);
		\draw[draw=black] ($(2.135)+(-0.3,0.3)$) rectangle  ($(2.-45)+(0.3,-0.3)$);
	\end{tikzpicture}
					
				}
				\caption{ The shape corresponding to $J-I$ from Example~\ref{ex:off-diag}.}
				\label{fig:4}
			\end{figure}

	\end{example}
	
Our next example is the identity matrix. For this example, $\U{\alpha}\cap \V{\alpha}$ is nonempty. 

    \begin{example}[Identity matrix] \label{ex:5}
	Consider $V(\V{\alpha})=\{\nu_1\}$, $\U{\alpha}=(\nu_1)$, $\V{\alpha} =(\nu_1)$, and $E(\alpha)=\big\{\big\}$. 
	Here the graphical representation of $\alpha$ is a single vertex with both $\U{\alpha}$ and $\V{\alpha}$ containing this vertex; see Figure~\ref{fig:5}. 
	$\M_\alpha$ is the $n\times n$ identity matrix.
	\begin{figure}[H]
				\centering
				\scalebox{0.65}{
					\begin{tikzpicture}[shorten >=1pt,auto,node distance=3cm,
					thick,main node/.style={circle,fill=blue!20,draw,font=\sffamily\LARGE\bfseries}]
					\node[main node] (1) at (0,0)  {$\nu_1$};
					\node[draw=none,fill=none] (3) at (0,1.5) {\fontsize{20}{22.4}\selectfont$\U{\alpha}\cap \V{\alpha}$};
					\path[every node/.style={font=\sffamily},line width=1pt];
		\draw[draw=black] ($(1.135)+(-0.3,0.3)$) rectangle  ($(1.-45)+(0.3,-0.3)$);
		\draw[draw=black] ($(1.135)+(-0.4,0.4)$) rectangle  ($(1.-45)+(0.4,-0.4)$); 
	\end{tikzpicture}
					
				}
				\caption{ The shape corresponding to $I$ from Example~\ref{ex:5}.}
				\label{fig:5}
			\end{figure}

	\end{example}
	
We end with an example where $\W{\alpha}$ is nonempty.

\begin{example}[Shifted degree matrix] \label{ex:6}
	Consider $V(\V{\alpha})=\{\nu_1,w_1\}$, $\U{\alpha}=(\nu_1)$, $\V{\alpha} =(\nu_1)$, $\W{\alpha} =(w_1)$, and $E(\alpha)=\big\{\{\nu_1,w_1\}\big\}$. 
	Here the graphical representation of $\alpha$ is the graph with two vertices $\nu_1$ and $w_1$ with both $\U{\alpha}$ and $\V{\alpha}$ containing only $\nu_1$; see Figure~\ref{fig:6}. 
	$\M_\alpha$ is the $n\times n$ random matrix equal to $deg(G)-\frac{n-1}{2}I$, where $G \sim G(n,1/2)$ and $deg(G)$ is the diagonal matrix containing the degree of each vertex of $G$.
	\begin{figure}[H]
				\centering
				\scalebox{0.65}{
					\begin{tikzpicture}[shorten >=1pt,auto,node distance=3cm,
					thick,main node/.style={circle,fill=blue!20,draw,font=\sffamily\LARGE\bfseries}, ,two/.style={diamond,fill=red!5,draw,font=\sffamily\LARGE\bfseries}]
					\node[main node] (1) at (0,0)  {$\nu_1$};
					\node[main node] (2) at (0,-2.5)  {$w_1$};
					\node[draw=none,fill=none] (3) at (0,1.5) {\fontsize{20}{22.4}\selectfont$\U{\alpha}\cap \V{\alpha}$};
					\node[draw=none,fill=none] (3) at (0,-4) {\fontsize{20}{22.4}\selectfont$\W{\alpha}$};
					\path[every node/.style={font=\sffamily},line width=1pt]
				    (1) edge node [] {} (2);
		\draw[draw=black] ($(1.135)+(-0.3,0.3)$) rectangle  ($(1.-45)+(0.3,-0.3)$);
		\draw[draw=black] ($(1.135)+(-0.4,0.4)$) rectangle  ($(1.-45)+(0.4,-0.4)$);
	\end{tikzpicture}
					
				}
				\caption{ The shape corresponding to the shifted degree matrix from Example~\ref{ex:6}.}
				\label{fig:6}
			\end{figure}

	\end{example}
For more complicated examples of graph matrices, see Appendix~\ref{app:example2}.
	\subsection{Norm bounds on graph matrices}
	We are now ready to state our main theorem, which gives a norm bound on every graph matrix $\M_\alpha$ in terms of some combinatorial quantities related to the graphical representation of $\alpha$.
	In particular, it turns out that the minimum size of a vertex separator between $\U{\alpha}$ and $\V{\alpha}$ is the relevant quantity for the bound, where the vertex separator is defined as follows:
	\begin{definition}[Vertex separators]
	If $G$ is a graph and $U,V \subseteq V(G)$, we say $S$ is a vertex separator between $U$ and $V$ if  all paths from $U$ to $V$ intersect $S$.
	\end{definition}
	\begin{remark}
	We allow paths of length $0$, so any separator $S$ between $U$ and $V$ must contain $U \cap V$.
	\end{remark}
	\noindent Here we give an informal statement of our main theorem. We defer the formal statement to Theorem~\ref{thm:formal:mainresult}.
			\begin{mdframed}[style=box2]
	\begin{theorem}[Informal] \label{thm:mainresult}
		Let $\alpha = (\U{\alpha},\V{\alpha},\W{\alpha},E(\alpha))$  be  a shape without  isolated middle vertices, and let $\ssmin$ be the minimum size of a vertex separator between  $\U{\alpha}$ and $\V{\alpha}$. Then with high probability $\norm{\M_\alpha}\leq \widetilde{O}(n^{\frac{1}{2}(|V(\alpha)|-\ssmin)})$.
		
		Moreover, if $\alpha$ has $u$ isolated middle vertices, then $\norm{\M_\alpha}\leq \widetilde{O}(n^{\frac{1}{2}(|V(\alpha)|+u-\ssmin)})$.
		\end{theorem}
		\end{mdframed} 
	
		We illustrate these norm bounds by revisiting examples we considered in Sections \ref{sec:mot1}, \ref{sec:mot2}, and ~\ref{sec:ex}. We refer readers to Appendix~\ref{app:example2} for more examples.
		
		\begin{example}[Revisiting Example~\ref{ex:3}: $\pm 1$ Wigner matrix] 
	For the shape $\alpha$ from Example~\ref{ex:3}, we have  $\ssmin= 1$.
	By Theorem~\ref{thm:mainresult}, with high probability $\norm{\M_\alpha}\leq \widetilde{O}(n^{\frac{1}{2}})$. 
	Note that this is consistent with the well-known tight upper bound of $2(1+o(1))\sqrt{n}$~\cite{bai1988necessary,soshnikov1999universality} up to a  poly-logarithmic term.
	\end{example}
	
		\begin{example}[Revisiting Example~\ref{ex:off-diag}: off-diagonal $1$'s matrix] 
	For the shape $\alpha$ from Example~\ref{ex:off-diag}, we have  $\ssmin= 0$.
	By Theorem~\ref{thm:mainresult}, with high probability $\norm{\M_\alpha}\leq \widetilde{O}(n)$. 
	Calculating the norm of $J-I$ tells us that $||M_\alpha||=n-1$, again showing consistency up to a  poly-logarithmic term.
	\end{example}
	
		\begin{example}[Revisiting Example~\ref{ex:5}: identity matrix]  
	For the shape $\alpha$ from Example~\ref{ex:5}, we have  $\ssmin= 1$.
	By Theorem~\ref{thm:mainresult}, with high probability $\norm{\M_\alpha}\leq \widetilde{O}(1)$. 
	We know that $||M_\alpha||=1$, so this result is also consistent.
	\end{example}
	
		\begin{example}[Revisiting Example~\ref{ex:6}: shifted degree matrix] 
	For the shape $\alpha$ from Example~\ref{ex:5}, we have  $\ssmin= 1$.
	By Theorem~\ref{thm:mainresult}, with high probability $\norm{\M_\alpha}\leq \widetilde{O}(n^{\frac{1}{2}})$. Using Chernoff bounds, with high probability the norm is on the order of $\sqrt{n\log(n)}$ so this is consistent as well. 
	\end{example}
	
    \begin{example}[Revisiting clique indicator]\label{revisitingcliqueindicatorexample}
    Let's consider the clique indicator 
    \begin{align}
    \cli = \frac{1}{2^6}\cdot \sum_{\substack{\alpha:~\text{shape}\\
    \U{\alpha} = (u_1,u_2), ~\V{\alpha} = (v_1,v_2)}}{M_{\alpha}}\,.
    \end{align}
    which we described in Section \ref{sec:mot1}.
        \begin{enumerate}
            \item For shapes $\alpha$ with no edges between $U_{\alpha}$ and $V_{\alpha}$, $s_{min} = 0$ so $\norm{M_\alpha}$ is $\widetilde{O}(n^2)$.
            \item For shapes $\alpha$ such that the maximum matching between $U_{\alpha}$ and $V_{\alpha}$ has size $1$, $s_{min} = 1$ so $\norm{M_{\alpha}}$ is $\widetilde{O}(n^{\frac{3}{2}})$. These matrices have considerably larger norm than a random $n(n-1) \times n(n-1)$ matrix and this was a major hurdle for proving tight degree $4$ SoS lower bounds for planted clique \cite{Norms,DBLP:journals/corr/HopkinsKP15,DBLP:journals/corr/DeshpandeM15,DBLP:journals/corr/RaghavendraS15}.
            \item For shapes $\alpha$ such that there is a matching with $2$ edges between $U_{\alpha}$ and $V_{\alpha}$, $s_{min} = 2$ so $\norm{M_{\alpha}}$ is $\widetilde{O}(n)$. As far as their norms are concerned, these matrices behave like random matrices.
        \end{enumerate}
    \end{example}
	\begin{example}[Revisiting triangle counter]
	Let's consider the matrix
	\begin{align}
    \tri = \frac{1}{2^3}\cdot \sum_{\substack{\alpha:~\text{shape}\\
    \U{\alpha} = (u_1), ~\V{\alpha} = (v_1),~\W{\alpha} =\{w_1\}   }}{M_{\alpha}}\,.
    \end{align}
    which we described in Section \ref{sec:mot2}.
	    \begin{enumerate}
	        \item For the shape $\alpha$ with no edges, $s_{min} = 0$ and $u = 1$ so $\norm{M_{\alpha}}$ is $\widetilde{O}(n^2)$.
	        \item For shapes $\alpha$ with one edge, either $s_{min} = u = 1$ (if the edge is between $u$ and $v$) or $s_{min} = u = 0$ (if the edge is between $w$ and $u$ or $v$). In either case, $u - s_{min} = 0$ so $\norm{M_{\alpha}}$ is $\widetilde{O}(n^{\frac{3}{2}})$. These matrices correspond to the facts that vertices with higher degrees are contained in more triangles and a pair of vertices $\{u,v\}$ can only be contained in triangles if $(u,v) \in E(G)$. 
	        \item For shapes $\alpha$ with $2$ or $3$ edges, $s_{min} = 1$ and $u = 0$ so $\norm{M_{\alpha}}$ is $\widetilde{O}(n)$. 
	    \end{enumerate}
	\end{example}		
 
   \paragraph{Agenda for the proof of Theorem~\ref{thm:mainresult}:} The proof requires some technical preparations, and the techniques will be detailed over the next few sections (Sections~\ref{tech}, \ref{pm1},  \ref{sec:modified}, and \ref{pf:main}):
    \begin{enumerate}
        \item In Section~\ref{tech}, we describe the main technical tool we use for proving probabilistic norm bounds, namely the \emph{trace power method}.
        
        \item In Section~\ref{pm1}, we consider $\pm 1$ Wigner matrices (Example~\ref{ex:3}) and illustrate how we can obtain probabilistic norm bounds using the trace power method. While this analysis is not new and is less precise than the sharp analyses from the prior works~\cite{bai1988necessary,soshnikov1999universality}, it serves as a useful warm-up for analyzing more general graph matrices.
        
        \item In Section~\ref{sec:modified}, we discuss how to modify the techniques for $\pm 1$ Wigner matrices to cover more general graph matrices.
        In particular, we describe a \emph{vertex partitioning argument} which is crucial in simplifying the counting arguments for more general graph matrices.
        \item In Section~\ref{pf:main}, we finally prove Theorem~\ref{thm:mainresult} based on the techniques developed up to that section.
    \end{enumerate}
    Having outlined our agenda for the proof, we begin with the trace power method.

	\section{The trace power method and constraint graphs} \label{tech}

	This section describes our main technical tools for establishing probabilistic spectral norm bounds, the trace power method and constraint graphs.
	\subsection{The trace power method}
	
	The most well-known technique for obtaining probabilistic spectral norm bounds on random matrices is called \emph{the trace power method}, perhaps more well-known as \emph{the moment method} (see e.g. \cite[Chapter 2.3.4]{tao2012topics}). We use the trace power method in the form of the following lemma:
\begin{mdframed}[style=box]
\begin{lemma}[The trace power method]\label{basictracepower}
		Assume that $M$ is a random matrix and that we have bounds $\{\B{2q}: q \in \mathbb{N}\}$ such that for all $q \in \mathbb{N}$, $\E\left[\tr\left((M{M^\top})^q\right)\right] \leq \B{2q}$.
		Then, for all $\epsilon > 0$,
		\[
		\pP\left[\norm{M} > \min_{q \in \mathbb{N}}{\left\{\sqrt[2q]{\frac{\B{2q}}{\epsilon}}\right\}}\right] < \epsilon.
		\]
	\end{lemma}
\end{mdframed}	
	\begin{proof}
		To prove this, we first show that for all $q \in \mathbb{N}$, ${\norm{M}}^{2q} \leq \tr\left((M{M^\top})^q\right)$.
			Indeed, this is an easy consequence from linear algebra.
			First, observe that 
			\[
			\norm{M} = \max_{u:\norm{u}_2=1}{ \norm{Mu}_2 } = \max_{u:\norm{u}_2=1}{ \sqrt{{u^\top}{M^\top}Mu} } = \sqrt{\lambda_{\max}({M^\top}M)} = \sqrt{\lambda_{\max}(M{M^\top})}\,.
			\]
			Now let $\lambda_1,\ldots,\lambda_n\geq 0$ be the eigenvalues of $M{M^\top}$.
			Then, for all $q \in \mathbb{N}$, $\tr\left((M{M^\top})^q\right) = \sum_{i=1}^{n}{{\lambda_i}^q} \geq {\lambda_{\max}(M{M^\top})}^q = {\norm{M}}^{2q}$, as needed.
		
		Having established the above inequality,  Lemma~\ref{basictracepower} is an immediate consequence of  Markov's inequality: for all $q \in \mathbb{N}$ and $\epsilon > 0$,
		\begin{align*}
		\pP\left[\norm{M} > \sqrt[2q]{\frac{B(2q)}{\epsilon}}\right] &= 
		\pP\left[{\norm{M}}^{2q} > \frac{B(2q)}{\epsilon}\right] \\
		&\leq \pP\left[\tr\left((M{M^\top})^q\right) > \frac{\ex\left[\tr\left((M{M^\top})^q\right)\right]}{\epsilon}\right] < \epsilon.
		\end{align*}
	\end{proof}
The main message of the trace power method is as follows:
to obtain upper bounds on $\norm{M_{\alpha}}$, it is sufficient to bound $\E\left[\tr\left((M_{\alpha}M_{\alpha}^\top)^q\right)\right]$ for an appropriate $q \in \mathbb{N}$ (which may depend on $n$).
	\subsection{Constraint graphs} 
	 To compute $\E\left[\tr\left((M_{\alpha}M_{\alpha}^\top)^q\right)\right]$, we expand it out into a large sum. 
	 We will then encode the terms of this sum using a concept which we call \emph{constraint graphs}. 
	 For this, we need a few more definitions.
	\begin{definition} \label{def:copies}
		Given a shape $\alpha$ and a $q \in \mathbb{N}$, we define $H(\alpha,2q)$ to be the multi-graph which is formed as follows:
		\begin{enumerate}
			\item Take $q$ copies $\alpha_1,\ldots,\alpha_{q}$ of $\alpha$ and take $q$ copies $\alpha^{\top}_1,\ldots,\alpha^{\top}_q$ of $\alpha^{\top}$, where $\alpha^\top$ is the shape obtained from $\alpha$ by switching the role of $\U{\alpha}$ and $\V{\alpha}$.
			\item For all $i \in [q]$, we glue them together by setting $V_{\alpha_i} = U_{\alpha^{\top}_i}$ and $V_{\alpha^{\top}_i} = U_{\alpha_{i+1}}$ (where $\alpha_{q+1} = \alpha_1$).
			\item If there are overlapping edges after gluing, we keep them as multiple edges.
		\end{enumerate}
		We define $V(\alpha,2q) = V(H(\alpha,2q))$ and we define $E(\alpha,2q) = E(H(\alpha,2q))$.
	\end{definition}
	
	We now provide a few examples of $H(\alpha, 2q)$ for different $\alpha$.
	For each example, we will draw  $H(\alpha, 2q)$ according to the following rules:
	\begin{itemize}
	    \item We will denote the vertices in $\U{\alpha_i} \setminus \V{\alpha_i}$, $\V{\alpha_i} \setminus \U{\alpha_i}$, and $\W{\alpha_i}$ by $u_{j;i}$, $v_{j;i}$, and $w_{j;i}$, respectively.
	    Moreover, we will denote the vertices in $\W{\alpha_i^\top}$ by $w'_{j;i}$.
	    \item We denote the vertices in $\U{\alpha}\cap \V{\alpha}$ by $\nn_{j;0}$.
	The reason why the second subscript is $0$ is because such vertices appear as a single copy in $H(\alpha,2q)$ as we shall see in Figure~\ref{fig:5b}.
	\item Consequently, $\alpha_i^\top$ will consist of $\U{\alpha_i^\top} = (v_{j;i})$, $\V{\alpha_i^\top}=( u_{j;i+1})$, and $\W{\alpha_i^\top} =\{w'_{j:i}\}$ (where $u_{j;q+1} = u_{j;1}$).
	\item As for the edges, we draw the edges of $\alpha_i$ by red solid lines and the edges of $\alpha_i^\top$ by blue dashed lines. 
We make an exception for the edges such that both endpoints are contained in $\U{\alpha}$ or both endpoints are contained in $\V{\alpha}$.
We will draw these edges by purple double lines. We use double lines as these edges are double edges in the constraint graph where the overlap happened due to the gluing step of Definition~\ref{def:copies}.
	\end{itemize}
	
	Let us first consider the following  two examples where $\W{\alpha}=\emptyset$.
	\begin{figure}[H] 
		\centering
		\begin{subfigure}[H]{0.20\textwidth}
			\centering
			\scalebox{0.5}{
				\begin{tikzpicture}[shorten >=1pt,auto,node distance=3cm,
				thick,main node/.style={circle,fill=blue!20,draw,font=\sffamily\LARGE\bfseries}]
				\node[main node] (1) at (-2,0)  {$u_1$};
				\node[main node] (2) at (2,0)  {$v_1$};
				\path[every node/.style={font=\sffamily},line width=1pt]
				(1) edge node [] {} (2);
					\node[draw=none,fill=none] (5) at (-2,1.5) {\fontsize{20}{22.4}\selectfont$\U{\alpha}$};
					\node[draw=none,fill=none] (6) at (2,1.5) {\fontsize{20}{22.4}\selectfont$\V{\alpha}$};
					\draw[draw=black] ($(1.135)+(-0.3,0.3)$) rectangle  ($(1.-45)+(0.3,-0.3)$);
		\draw[draw=black] ($(2.135)+(-0.3,0.3)$) rectangle  ($(2.-45)+(0.3,-0.3)$);	\end{tikzpicture}
			}
			\caption{Figure \ref{fig:4a}: $\alpha$.}
			\label{fig:4a}
		\end{subfigure}
		\qquad\qquad
		\begin{subfigure}[H]{0.60\textwidth}
			\centering
			\scalebox{0.60}{
				\begin{tikzpicture}[shorten >=1pt,auto,node distance=3cm,
				thick,main node/.style={circle,fill=blue!20,draw,font=\sffamily\LARGE\bfseries}]
				
				\node[main node] (1) at (0,3)  {$u_{1;1}$};
				\node[main node] (3) at (2.25,2.25)  {$v_{1;1}$};
				\node[main node] (5) at (3,0)  {$u_{1;2}$};
				\node[main node] (7) at (2.25,-2.25)  {$v_{1;2}$};
				\node[main node] (9) at (0,-3)  {$u_{1;3}$};
				\node[main node] (11) at (-2.25,-2.25)  {$v_{1;3}$};
				\node[main node] (13) at (-3,0)  {$u_{1;4}$};
				\node[main node] (15) at (-2.25,2.25)  {$v_{1;4}$};
				
				\path[every node/.style={font=\sffamily},line width=1pt]
				(1) edge [red] node [] {} (3)
				(5) edge [blue, dashed] node [] {} (3)
				(5) edge [red] node [] {} (7)
				(9) edge [blue, dashed] node [] {} (7)
				(9) edge [red] node [] {} (11)
				(13) edge [blue, dashed] node [] {} (11)
				(13) edge [red] node [] {} (15)
				(1) edge [blue, dashed] node [] {} (15);
				\end{tikzpicture}
			}
			\caption{Figure \ref{fig:4b}: $H(\alpha,2q)$ with $q=4$.}
			\label{fig:4b}
		\end{subfigure}
	\end{figure}

	\begin{figure}[H]
		\centering
		\begin{subfigure}[H]{0.20\textwidth}
			\centering
			\scalebox{0.5}{
				\begin{tikzpicture}[shorten >=1pt,auto,node distance=3cm,
				thick,main node/.style={circle,fill=blue!20,draw,font=\sffamily\LARGE\bfseries}]
				\node[main node] (1) at (-2,2)  {$u_1$};
				\node[main node] (2) at (-2,-2)  {$u_2$};
				\node[main node] (3) at (2,2)  {$v_1$};
				\node[main node] (4) at (2,-2)  {$v_2$};
				\node[main node] (7) at (0,-5)  {$\nn_1$};
				\path[every node/.style={font=\sffamily},line width=1pt]
				(1) edge node [] {} (3)
				(3) edge node [] {} (2)
				(2) edge node [] {} (4)
				(7) edge node [] {} (4);
				\node[draw=none,fill=none] (5) at (-2,3.5) {\fontsize{20}{22.4}\selectfont$\U{\alpha}$};
				\node[draw=none,fill=none] (6) at (2,3.5) {\fontsize{20}{22.4}\selectfont$\V{\alpha}$};
				\node[draw=none,fill=none] (6) at (-0.5,-3.5) {\fontsize{18}{22.4}\selectfont$\U{\alpha}\cap \V{\alpha}$};
				\draw[draw=black] ($(1.135)+(-0.3,0.3)$) rectangle  ($(2.-45)+(0.3,-0.3)$);
		\draw[draw=black] ($(3.135)+(-0.3,0.3)$) rectangle  ($(4.-45)+(0.3,-0.3)$);
				\draw[draw=black] ($(7.135)+(-0.3,0.3)$) rectangle  ($(7.-45)+(0.3,-0.3)$);
				\draw[draw=black] ($(7.135)+(-0.4,0.4)$) rectangle  ($(7.-45)+(0.4,-0.4)$);
				\end{tikzpicture}
				
			}
			\caption{Figure \ref{fig:5a}: $\alpha$.}
			\label{fig:5a}
		\end{subfigure}
		\qquad\qquad
		\begin{subfigure}[H]{0.60\textwidth}
			\centering
			\scalebox{0.60}{
				\begin{tikzpicture}[shorten >=1pt,auto,node distance=3cm,
				thick,main node/.style={circle,fill=blue!20,draw,font=\sffamily\LARGE\bfseries}]
				
				\node[main node] (1) at (0,6)  {$u_{1;1}$};
				\node[main node] (2) at (0,3)  {$u_{2;1}$};
				\node[main node] (3) at (4.5,4.5)  {$v_{1;1}$};
				\node[main node] (4) at (2.25,2.25)  {$v_{2;1}$};
				\node[main node] (5) at (6,0)  {$u_{1;2}$};
				\node[main node] (6) at (3,0)  {$u_{2;2}$};
				\node[main node] (7) at (4.5,-4.5)  {$v_{1;2}$};
				\node[main node] (8) at (2.25,-2.25)  {$v_{2;2}$};
				\node[main node] (9) at (0,-6)  {$u_{1;3}$};
				\node[main node] (10) at (0,-3)  {$u_{2;3}$};
				\node[main node] (11) at (-4.5,-4.5)  {$v_{1;3}$};
				\node[main node] (12) at (-2.25,-2.25)  {$v_{2;3}$};
				\node[main node] (13) at (-6,0)  {$u_{1;4}$};
				\node[main node] (14) at (-3,0)  {$u_{2;4}$};
				\node[main node] (15) at (-4.5,4.5)  {$v_{1;4}$};
				\node[main node] (16) at (-2.25,2.25)  {$v_{2;4}$};
				\node[main node] (17) at (0,0)  {$\nn_{1;0}$};
				\path[every node/.style={font=\sffamily},line width=1pt]
				(1) edge [red] node [] {} (3)
				(2) edge [red] node [] {} (3)
				(2) edge [red] node [] {} (4)
				(5) edge [blue, dashed] node [] {} (3)
				(6) edge [blue, dashed] node [] {} (3)
				(6) edge [blue, dashed] node [] {} (4)
				(5) edge [red] node [] {} (7)
				(6) edge [red] node [] {} (7)
				(6) edge [red] node [] {} (8)
				(9) edge [blue, dashed] node [] {} (7)
				(10) edge [blue, dashed] node [] {} (7)
				(10) edge [blue, dashed] node [] {} (8)
				(9) edge [red] node [] {} (11)
				(10) edge [red] node [] {} (11)
				(10) edge [red] node [] {} (12)
				(13) edge [blue, dashed] node [] {} (11)
				(14) edge [blue, dashed] node [] {} (11)
				(14) edge [blue, dashed] node [] {} (12)
				(13) edge [red] node [] {} (15)
				(14) edge [red] node [] {} (15)
				(14) edge [red] node [] {} (16)
				(1) edge [blue, dashed] node [] {} (15)
				(2) edge [blue, dashed] node [] {} (15)
				(2) edge [blue, dashed] node [] {} (16)
				(17) edge [purple, double distance between line centers=0.3em] node [] {} (4)
				(17) edge [purple, double distance between line centers=0.3em] node [] {} (8)
				(17) edge [purple,double distance between line centers=0.3em] node [] {} (12)
				(17) edge [purple,double, double distance between line centers=0.3em] node [] {} (16);
				\end{tikzpicture}
			}
			\caption{Figure \ref{fig:5b}: $H(\alpha,2q)$ with $q=4$.}
			\label{fig:5b}
		\end{subfigure}
	\end{figure}
	From Figures~\ref{fig:4b} and \ref{fig:5b}, one can easily verify that each vertex in $\U{\alpha}\cup \V{\alpha}\setminus (\U{\alpha}\cap \V{\alpha})$ is duplicated $q$ times in $H(\alpha,2q)$, and each vertex in $\U{\alpha}\cap \V{\alpha}$ appears as a single copy.
	Now let us consider an example where $\W{\alpha}\neq \emptyset$. 
	
	\begin{figure}[H]
		\centering
		\begin{subfigure}[H]{0.30\textwidth}
			\centering
			\scalebox{0.5}{
				\begin{tikzpicture}[shorten >=1pt,auto,node distance=3cm,
				thick,main node/.style={circle,fill=blue!20,draw,font=\sffamily\LARGE\bfseries}]
				\node[main node] (1) at (-2,0)  {$u_1$};
				\node[main node] (2) at (-2,-2)  {$u_2$};
				\node[main node] (3) at (2,0)  {$v_1$};
				\node[main node] (4) at (0,0)  {$w_1$};
				\path[every node/.style={font=\sffamily},line width=1pt]
				(1) edge node [] {} (4)
				(1) edge node [] {} (2)
				(2) edge node [] {} (4)
				(3) edge node [] {} (4);
				\node[draw=none,fill=none] (5) at (-2,1.5) {\fontsize{20}{22.4}\selectfont$\U{\alpha}$};
				\node[draw=none,fill=none] (6) at (2,1.5) {\fontsize{20}{22.4}\selectfont$\V{\alpha}$};
				\node[draw=none,fill=none] (7) at (0,1.5) {\fontsize{20}{22.4}\selectfont$\W{\alpha}$};
					\draw[draw=black] ($(1.135)+(-0.3,0.3)$) rectangle  ($(2.-45)+(0.3,-0.3)$);
		\draw[draw=black] ($(3.135)+(-0.3,0.3)$) rectangle  ($(3.-45)+(0.3,-0.3)$);
				\end{tikzpicture}
				
			}
			\caption{Figure \ref{fig:6a}: $\alpha$.}
			\label{fig:6a}
		\end{subfigure}
		\qquad
		\begin{subfigure}[H]{0.625\textwidth}
			\centering
			\scalebox{0.55}{
				\begin{tikzpicture}[shorten >=1pt,auto,node distance=3cm,
				thick,main node/.style={circle,fill=blue!20,draw,font=\sffamily\LARGE\bfseries}]
				
				\node[main node] (1) at (0,6)  {$u_{1;1}$};
				\node[main node] (2) at (0,2)  {$u_{2;1}$};
				\node[main node] (3) at (6,0)  {$v_{1;1}$};
				\node[main node] (4) at (3,3)  {$w_{1;1}$};
				\node[main node] (5) at (0,-6)  {$u_{1;2}$};
				\node[main node] (6) at (0,-2)  {$u_{2;2}$};
				\node[main node] (7) at (3,-3)  {$w'_{1;1}$};
				\node[main node] (8) at (-6,0)  {$v_{1;2}$};
				\node[main node] (9) at (-3,-3)  {$w_{1;2}$};
				\node[main node] (10) at (-3,3)  {$w'_{1;2}$};

				\path[every node/.style={font=\sffamily},line width=1pt]
				(1) edge [red] node {} (4)
				(2) edge [red] node {} (4)
				(3) edge [red] node {} (4)
				(3) edge [blue,dashed] node {} (7)
				(5) edge [blue,dashed] node {} (7)
				(6) edge [blue,dashed] node {} (7)
				(5) edge [red] node {} (9)
				(6) edge [red] node {} (9)
				(8) edge [red] node {} (9)
				(1) edge [blue,dashed] node {} (10)
				(2) edge [blue,dashed] node {} (10)
				(8) edge [blue,dashed] node {} (10)
			(1) edge [purple, double distance between line centers=0.3em] node {} (2)
			(5) edge [purple, double distance between line centers=0.3em] node {} (6);
				\end{tikzpicture}
			}
			\caption{Figure \ref{fig:6b}: $H(\alpha,2q)$ with $q=2$.}
			\label{fig:6b}
		\end{subfigure}
	\end{figure}
	
	From Figure~\ref{fig:6b}, one can easily verify that each vertex in $\W{\alpha}$ is duplicated $2q$ times in $H(\alpha,2q)$. 
To conclude our observations thus far, we obtain:
\begin{proposition} \label{prop:num}
    For a shape $\alpha=  (\U{\alpha},\V{\alpha},\W{\alpha},E(\alpha))$ and a $q\in\na$,
    \begin{align*}
        |V(\alpha,2q)|= q\cdot (|\U{\alpha}|+|\V{\alpha}| - 2\cdot |\U{\alpha}\cap \V{\alpha}|) + 2q\cdot |\W{\alpha}| + |\U{\alpha}\cap \V{\alpha}|\,.
    \end{align*}
\end{proposition}
	
Now we consider a map on $V(\alpha,2q)$ that assigns indices from $[n]$ to the vertices of a constraint graph.
In order for such map to be \emph{valid}, such map should not assign the same index to two different vertices within the same copy of $\alpha$ or $\alpha^T$.
Motivated by this, we make the following definition:
\begin{definition}[Piecewise injectivity]
		We say that a map $\phi:V(\alpha,2q) \to [n]$ is piecewise injective if $\phi$ is injective on each piece $V(\alpha_i)$ and each piece $V(\alpha^{\top}_i)$ for all $i \in [q]$. In other words, $\phi(u) \neq \phi(v)$ whenever $u,v \in V(\alpha_i)$ for some $i \in [q]$ or $u,v \in V(\alpha^\top_i)$ for some $i \in [q]$.
	\end{definition}
	With this definition of piecewise injectivity, we can express the trace power term as follows:
	\begin{proposition} \label{expansion}
		For all shapes $\alpha$ and all $q \in \mathbb{N}$,
		\[
		\E\left[\tr\left((M_{\alpha}M_{\alpha}^\top)^q\right)\right] = \sum_{\substack{\phi: V(\alpha,2q) \to [n]:\\
		\phi \text{ is piecewise injective}}}{\E[\chi_{\phi(E(\alpha,2q))}(G)]}\,.
		\]
	\end{proposition}
	We now observe that by symmetry, $\E[\chi_{\phi(E(\alpha,2q))}(G)]$ only depends on the set of pairs of vertices $\{(u,v): u,v \in V(\alpha,2q), \phi(u) = \phi(v)\}$ which are mapped to the same index. 
	We capture this set of pairs with a concept which we call \emph{constraint graphs}. 
	\begin{mdframed}[style=box]
	\begin{definition}[Constraint graphs] \label{def:constr}
		Given a set of vertices $V$ and a map $\phi:V \to [n]$, we construct the constraint graph $C(\phi)$ on $V$ associated to $\phi$ as follows:
		\begin{enumerate}
			\item We take $V(C(\phi)) = V$.
			\item For each pair of vertices $u,v \in V$ such that $\phi(u) = \phi(v)$, we add a constraint edge between $u$ and $v$.
			\item As long as there is a cycle, we delete one edge of this cycle (this choice is arbitrary). We do this until there are no cycles left.
		\end{enumerate}
		We say that two constraint graphs $C,C'$ on $V$ are equivalent (which we write as $C \equiv C'$) if for all $u,v \in V$, there is a path of constraint edges between $u$ and $v$ in $C$ if and only if there is a path of constraint edges between $u$ and $v$ in $C'$.
	\end{definition}
	\end{mdframed}
		\begin{remark}
		We delete all cycles from the constraint graph $C$ so that there are no redundant constraints and $|\phi(V)| = |V| - |E(C)|$.
	\end{remark}
	\noindent For our purpose, we specifically consider constraint graphs on $H(\alpha,2q)$ as follows:
	\begin{mdframed}[style=box]
	\begin{definition}[Constraint graphs on $H(\alpha,2q)$]
		We define $\mathcal{C}_{(\alpha,2q)} = \{C(\phi): \phi: V(\alpha,2q) \to [n] \text{ is piecewise injective}\}$ to be the set of all possible constraint graphs on $V(\alpha,2q)$ (up to equivalence) which come from a piecewise independent map $\phi: V(\alpha,2q) \to [n]$.
		
		Moreover, given a constraint graph $C \in \mathcal{C}_{(\alpha,2q)}$, we make the following definitions:
		\begin{enumerate}
			\item We define $N(C) = |\{\phi: V(\alpha,2q) \to [n]: \phi \text{ is piecewise injective}, C(\phi) \equiv C\}|$.
			\item We define $\val{C} = \E[\chi_{\phi(E(\alpha,2q))}(G)]$ where $\phi: V(\alpha,2q) \to [n]$ is any piecewise injective map such that $C(\phi)=C$. Note that this is well-defined due to the observation below Proposition~\ref{expansion}.
		\end{enumerate}
	\end{definition}
	\end{mdframed}
	With these definitions, we can simplify $\E\left[\tr\left((M_{\alpha}M_{\alpha}^\top)^q\right)\right]$ as follows:
	\begin{proposition} \label{diffrep}
		For all shapes $\alpha$ and all $q \in \mathbb{N}$,
		\[
		\E\left[\tr\left((M_{\alpha}M_{\alpha}^\top)^q\right)\right] = \sum_{C \in \mathcal{C}_{(\alpha,2q)}}{N(C)\val{C}}\,.
		\]
	\end{proposition}
	\noindent Thus, to analyze $\E\left[\tr\left((M_{\alpha}M_{\alpha}^\top)^q\right)\right]$ it is sufficient to analyze $N(C)$ and $\val{C}$ for all $C \in \mathcal{C}_{(\alpha,2q)}$.
	To that end, we make the following observation:
	\begin{proposition}\label{even}
		For every constraint graph $C \in \mathcal{C}_{(\alpha,2q)}$, $\val{C} = 1$ if every edge in $\phi(E(\alpha,2q))$ appears an even number of times and $\val{C} = 0$ otherwise (where $\phi: V(\alpha,2q) \to [n]$ is any piecewise injective map such that $C(\phi)=C$).
	\end{proposition}
	\begin{proof}
		Observe that for any multi-set of edges $E$, $\chi_E(G) = \chi_{E_{reduced}(G)}$ where $E_{reduced} = \{e: e \text{ appears in } E \text{ an odd number of times}\}$. Thus, $\E[\chi_E(G)] = 1$ if every edge in $E$ appears an even number of times and $\E[\chi_E(G)] = 0$ otherwise. Taking $E = \phi(E(\alpha,2q))$, the result follows.
	\end{proof}

	We include two examples of constraint graphs for the case when $\alpha$ consists of a single edge between $u_{1}$ and $v_{1}$ (i.e. the shape $\alpha$ from Example~\ref{ex:3} whose graph matrix is a $\pm 1$ Wigner matrix), and $q=4$. In both cases, we color equal edges with the same color and style.
	
		\begin{figure}[H]
			\centering
			\scalebox{0.65}{
				\begin{tikzpicture}[shorten >=1pt,auto,node distance=3cm,
				thick,main node/.style={circle,fill=blue!20,draw,font=\sffamily\LARGE\bfseries}]
				
				\node[main node] (1) at (0,4)  {$u_{1;1}$};
				\node[main node] (2) at (2.828,2.828)  {$v_{1;1}$};
				\node[main node] (3) at (4,0)  {$u_{1;2}$};
				\node[main node] (4) at (2.828,-2.828)  {$v_{1;2}$};
				\node[main node] (5) at (0,-4)  {$u_{1;3}$};
				\node[main node] (6) at (-2.828,-2.828)  {$v_{1;3}$};
				\node[main node] (7) at (-4,0)  {$u_{1;4}$};
				\node[main node] (8) at (-2.828,2.828)  {$v_{1;4}$};
				
				\path[every node/.style={font=\sffamily},line width=1pt]
			    (1) edge [red, solid] node [] {$\{\phi(u_{1;1}),\phi(v_{1;1})\}$} (2)
				(2) edge [red,solid] node [] {$\{\phi(v_{1;1}),\phi(u_{1;2})\}$} (3)
				(3) edge [blue,densely dotted] node [] {$\{\phi(u_{1;2}),\phi(v_{1;2})\}$} (4)
				(4) edge [olive,loosely dashed] node [] {$\{\phi(v_{1;2}),\phi(u_{1;3})\}$} (5)
				(5) edge [purple, densely dashdotdotted] node [] {$\{\phi(u_{1;3}),\phi(v_{1;3})\}$} (6)
				(6) edge [black, loosely dashdotdotted] node [] {$\{\phi(v_{1;3}),\phi(u_{1;4})\}$} (7)
				(7) edge [orange, densely dashed]  node [] {$\{\phi(u_{1;4}),\phi(v_{1;4})\}$} (8)
				(8) edge [blue,densely dotted]  node [] {$\{\phi(v_{1;4}),\phi(u_{1;1})\}$} (1);
				\path[every node/.style={font=\sffamily},line width=2pt]
				(1) edge node [] {} (3)
				(4) edge node [] {} (8)
				(2) edge node [] {} (6);
				\end{tikzpicture}
			}
			\caption{In this constraint graph $C$, $u_{1;1}=u_{1;2}$, $v_{1;1}=v_{1;3}$, and $v_{1;2}=v_{1;4}$. Here $\val{C}=0$ as we see edges that appear only once. More specifically, edges $\{v_{1;2},u_{1;3}\}$, $\{u_{1;3},v_{1;3}\}$, $\{v_{1;3},u_{1;4}\}$, and $\{u_{1;4},v_{1;4}\}$ appear only once.}\label{P10}
		\end{figure}
		
		\begin{figure}[H]
			\centering
			\scalebox{0.6}{
				\begin{tikzpicture}[shorten >=1pt,auto,node distance=3cm,
				thick,main node/.style={circle,fill=blue!20,draw,font=\sffamily\LARGE\bfseries}]
				
				\node[main node] (1) at (0,4)  {$u_{1;1}$};
				\node[main node] (2) at (2.828,2.828)  {$v_{1;1}$};
				\node[main node] (3) at (4,0)  {$u_{1;2}$};
				\node[main node] (4) at (2.828,-2.828)  {$v_{1;2}$};
				\node[main node] (5) at (0,-4)  {$u_{1;3}$};
				\node[main node] (6) at (-2.828,-2.828)  {$v_{1;3}$};
				\node[main node] (7) at (-4,0)  {$u_{1;4}$};
				\node[main node] (8) at (-2.828,2.828)  {$v_{1;4}$};
				
				\path[every node/.style={font=\sffamily},line width=1pt]
				(1) edge [red, solid] node [] {$\{\phi(u_{1;1}),\phi(v_{1;1})\}$} (2)
				(2) edge [blue,dotted] node [] {$\{\phi(v_{1;1}),\phi(u_{1;2})\}$} (3)
				(3) edge [black,dashed] node [] {$\{\phi(u_{1;2}),\phi(v_{1;2})\}$} (4)
				(4) edge [black,dashed] node [] {$\{\phi(v_{1;2}),\phi(u_{1;3})\}$} (5)
				(5) edge [orange,dashdotted] node [] {$\{\phi(u_{1;3}),\phi(v_{1;3})\}$} (6)
				(6) edge [orange,dashdotted] node [] {$\{\phi(v_{1;3}),\phi(u_{1;4})\}$} (7)
				(7) edge [blue, dotted] node [] {$\{\phi(u_{1;4}),\phi(v_{1;4})\}$} (8)
				(8) edge [red,solid] node [] {$\{\phi(v_{1;4}),\phi(u_{1;1})\}$} (1);
				\path[every node/.style={font=\sffamily},line width=2pt]
				(8) edge node [] {} (2)
				(3) edge node [] {} (5)
				(5) edge node [] {} (7);
				\end{tikzpicture}
			}
			\caption{In this constraint graph $C$, $v_{1;1}=v_{1;4}$ and $u_{1;2}=u_{1;3}=u_{1;4}$. Since every edge appears twice, we have $\val{C}=1$. }\label{P11}
		\end{figure}
	\begin{proposition} \label{numchoice}
		For any constraint graph $C \in \mathcal{C}_{(\alpha,2q)}$, $N(C) \leq n^{|V(\alpha,2q)| - |E(C)|}$.
	\end{proposition}
	\begin{proof}
		Observe that choosing a piecewise injective map $\phi: V(\alpha,2q) \to [n]$ with constraint graph $C$ is equivalent to choosing the distinct indices of $\phi(V(\alpha,2q))$. Since there are $|V(\alpha,2q)| - |E(C)|$ such indices, $N(C) = \frac{n!}{(n-|V(\alpha,2q)| + |E(C)|)!} \leq n^{|V(\alpha,2q)| - |E(C)|}$
	\end{proof}

	\section{Warm-up: The \texorpdfstring{$\pm{1}$}{+1} Wigner matrix}\label{pm1}
     
	As a warm-up, we consider the $\pm{1}$ Wigner matrix, i.e., $M_\alpha$ with the shape $\alpha$ from Example~\ref{ex:3}.
	The Wigner matrix and its norm have already been studied extensively; in particular,   sharp analyses from the prior works~\cite{bai1988necessary,soshnikov1999universality} demonstrate that with high probability, the norm of an $n \times n$ symmetric random $\pm{1}$ matrix is  $2(1+o(1))\sqrt{n}$. While our upper bound will not be as strong, it will illustrate the key ideas for the general case.
	\begin{theorem}\label{singular main}
		For all $\epsilon > 0$, the following probabilistic upper bound holds:
		\[
		\pP\left[\norm{M_\alpha} > \sqrt{8e\left\lceil{\log\left(\frac{n}{\epsilon}\right)}\right\rceil{n}}\right] < \epsilon\,.
		\]
	\end{theorem}
	\begin{proof}
		We prove this probabilistic norm bound using the trace power method. We have the following upper bound on $\E\left[\tr\left((M_{\alpha}M_{\alpha}^\top)^q\right)\right]$.
		\begin{lemma}\label{traceupper}
			For all $q \in \mathbb{N}$ such that $q \leq \frac{n}{4}$,
			\[
			\E\left[\tr\left((M_{\alpha}M_{\alpha}^\top)^q\right)\right] < 2^{2q}(2q)^{q}n^{q+1}\,.
			\]
		\end{lemma}
		\begin{proof}
			Recall that
			\[
			\E\left[\tr\left((M_{\alpha}M_{\alpha}^\top)^q\right)\right] = \sum_{C \in \mathcal{C}_{(\alpha,2q)}}{N(C)\val{C}}\,.
			\]
			Also recall from Proposition~\ref{numchoice} that $N(C) \leq n^{|V(\alpha,2q)| - |E(C)|}$ and that for all $C \in \mathcal{C}_{(\alpha,2q)}$, $\val{C} = 0$ or $1$.
			
			To prove our upper bound on $\E\left[\tr\left((M_{\alpha}M_{\alpha}^\top)^q\right)\right]$, we prove the following statements:
			\begin{enumerate}
				\item For any $C \in \mathcal{C}_{(\alpha,2q)}$ such that $\val{C} \neq 0$, $|E(C)| \geq q-1$.
				\item For all $k$, there are at most $2^{2q-1}(2q)^k$ constraint graphs $C \in \mathcal{C}_{(\alpha,2q)}$ which have exactly $k$ edges.
			\end{enumerate}
			Assuming these statements are true, since $q \leq \frac{n}{4}$,
			\[
			\E\left[\tr\left((M_{\alpha}M_{\alpha}^\top)^q\right)\right] \leq \sum_{k = q-1}^{2q-1}{n^{2q-k}2^{2q-1}(2q)^{k}} = 
			2^{2q-1}(2q)^{q-1}n^{q+1}\sum_{k' = 0}^{q}{\left(\frac{2q}{n}\right)^{k'}} < 2^{2q}(2q)^{q}n^{q+1}\,.
			\]
			Thus, we just need to prove these two statements.
			We begin with the first statement.
			\begin{definition}[Cycles] \label{def:cycle}
				For all $k \in \mathbb{N}$, we take $C_{2k}$ to be the cycle on $2k$ vertices, i.e. $C_{2k}$ is the multi-graph with vertices $V(C_{2k}) = \{\xx_i: i \in [2k]\}$ and edges $E(C_{2k}) = \{\{\xx_i,\xx_{i+1}\}~:~ i \in [2k]\}$ (where we take $v_{2k+1} = v_1$).
			\end{definition}
			Since the shape of interest is a single edge (see Example~\ref{ex:3}), the constraint graph $H(\alpha,2q)$ will be isomorphic to $C_{2q}$.
			Here and below, we will identify $H(\alpha,2q)$ with $C_{2q}$ to simplify notation.
			Note from Proposition~\ref{even} that when $\val{C}\neq 0$, it follows that each edge in $\phi(C_{2q})$ appears an even number of times, in particular, at least twice. 
			Hence, the first statement follows immediately from the following lemma. 
			\begin{lemma} \label{mink-1}
				For all $k \in \mathbb{N}$, for any map $\phi: V(C_{2k}) \to [n]$ such that each edge in $\phi(E(C_{2k}))$ appears at least twice, $|E(C(\phi))| \geq k-1$.
			\end{lemma}
			\begin{proof}
				We prove this lemma by induction. If $k = 1$, the result is trivial. If $k > 1$, assume that the lemma is true up to $k-1$ and let $\phi: V(C_{2k}) \to [n]$ be a map such that each edge in $\phi(E(C_{2k}))$ appears at least twice.
				
	For our inductive argument, we define the notion of unique vertices: we say a vertex $x$ in $C_{2k}$ is \emph{unique} if $x$ is the only vertex with the index $\phi(x)$.
				There are two cases to consider depending on the existence of a unique vertex.
		
		\begin{enumerate}
		    \item First, assume that $C(\phi)$ has no unique vertices.
		    In this case, $|E(C(\phi))| = \frac{V(C_{2k})}{2} = k>k-1$ and thus the result readily follows.
		    See Figure \ref{P12} for an illustration of this fact.
			\begin{figure}[H]
			\centering
			\scalebox{0.65}{
				\begin{tikzpicture}[shorten >=1pt,auto,node distance=3cm,
				thick,main node/.style={circle,fill=blue!20,draw,font=\sffamily\LARGE\bfseries}]
				
				\node[main node] (1) at (0,4)  {$x_1$};
				\node[main node] (2) at (2.828,2.828)  {$x_2$};
				\node[main node] (3) at (4,0)  {$x_3$};
				\node[main node] (4) at (2.828,-2.828)  {$x_4$};
				\node[main node] (5) at (0,-4)  {$x_5$};
				\node[main node] (6) at (-2.828,-2.828)  {$x_6$};
				\node[main node] (7) at (-4,0)  {$x_7$};
				\node[main node] (8) at (-2.828,2.828)  {$x_8$};
				
				\path[every node/.style={font=\sffamily},line width=1pt]
				(1) edge [red] node [] {} (2)
				(2) edge [blue,dashed] node [] {} (3)
				(3) edge [orange,dotted] node [] {} (4)
				(4) edge [black,dashdotted] node [] {} (5)
				(5) edge [red] node [] {} (6)
				(6) edge [blue,dashed] node [] {} (7)
				(7) edge [orange,dotted] node [] {} (8)
				(8) edge [black,dashdotted] node [] {} (1);
				\path[every node/.style={font=\sffamily},line width=2pt]
				(1) edge node [] {} (5)
				(3) edge node [] {} (7)
				(2) edge node [] {} (6)
				(4) edge node [] {} (8);
				\end{tikzpicture}
			}
			\caption{An illustration of a constraint graph $C_{8}$ without unique vertices. Edges of the same value are drawn with the same color and style. Notice that in this case, we have $|E(C_8)| = \frac{V(C_{8})}{2} = 4$ and hence the desired inequality holds.}\label{P12}
		    \end{figure}
				
\item Hence, we may assume that there is a unique vertex  $x_i \in V(C(\phi))$. 
We must have that $\phi(x_{i-1}) = \phi(x_{i+1})$ (where we take $x_{0} = x_{2k}$ and $x_{2k+1} = x_1$) because otherwise both edges $\phi(\{x_{i-1},x_i\})$ and $\phi(\{x_{i},x_{i+1}\})$ would  appear exactly once in $\phi(E(C_{2k}))$.
				
	We now contract the vertices $x_{i-1}$ and $x_{i+1}$ together and delete the vertex $x_i$ to obtain  the cycle $C'_{2k-2}$ of length $2k-2$.
	Moreover, we define the restricted index map $\phi': V(C'_{2k-2}) \to [n]$ so that the indices match the corresponding indices assigned by $\phi$ before contraction.
Observe that every edge in $\phi(E(C'_{2k-2}))$ still appears at least twice because the two deleted edges $\{x_{i-1},x_i\}$ and  $\{x_{i},x_{i+1}\}$  were the only two edges with the value $\{\phi(x_{i-1}),\phi(x_i)\} = \{\phi(x_{i}),\phi(x_{i+1})\} $.

Thus, by the inductive hypothesis, $|E(C(\phi'))| \geq k-2$.
Since the contraction got rid of one constraint edge between  $x_{i-1}$ and $x_{i+1}$, we have $|E(C(\phi))| = |E(C(\phi'))| + 1$, from which we obtain $|E(C(\phi))| \geq k-1$, as needed. An illustration of this process is shown in Figure \ref{P13}.
\begin{figure}[H]
        		\centering
        		\begin{subfigure}[H]{0.40\textwidth}
        			\centering
        			\scalebox{0.5}{
        				\begin{tikzpicture}[shorten >=1pt,auto,node distance=3cm,
        				thick,main node/.style={circle,fill=blue!20,draw,font=\sffamily\LARGE\bfseries}]

				\node[main node] (1) at (0,4)  {$x_1$};
				\node[main node] (2) at (2.828,2.828)  {$x_2$};
				\node[main node] (3) at (4,0)  {$x_3$};
				\node[main node] (4) at (2.828,-2.828)  {$x_4$};
				\node[main node] (5) at (0,-4)  {$x_5$};
				\node[main node] (6) at (-2.828,-2.828)  {$x_6$};
				\node[main node] (7) at (-4,0)  {$x_7$};
				\node[main node] (8) at (-2.828,2.828)  {$x_8$};
        				
        				\path[every node/.style={font=\sffamily},line width=1pt]
        				(1) edge [orange] node [] {} (2)
        				(2) edge [orange] node [] {} (3)
        				(3) edge [blue,dotted] node [] {} (4)
        				(4) edge [red,dashdotted] node [] {} (5)
        				(5) edge [red,dashdotted] node [] {} (6)
        				(6) edge [blue,dotted] node [] {} (7)
        				(7) edge [black,dashed] node [] {} (8)
        				(8) edge [black,dashed] node [] {} (1);
        				\path[every node/.style={font=\sffamily},line width=2pt]
        				(1) edge node [] {} (3)
        				(1) edge node [] {} (7)
        				(4) edge node [] {} (6);
        				
        				\draw [ultra thick, red] (0,-4) circle [radius=1];
        				\end{tikzpicture}

        			}
        			\caption{Figure \ref{fig:pm1-a}: $C(\phi)$.}
        				\label{fig:pm1-a}
        		\end{subfigure}
        	   		\begin{subfigure}[H]{0.05\textwidth}
        		    \centering
        		    \scalebox{1}{
        		        \begin{tikzpicture}[shorten >=1pt,auto,node distance=3cm,
        				thick]
                        \draw[->, very thick] (-1,0) -- (0,0);
        				\node[] (1) at (-0.5,1)  {Contraction};
                        \end{tikzpicture}
        		    
        		    }
        		\end{subfigure}
        		\qquad\qquad 
              		\begin{subfigure}[H]{0.40\textwidth}
        			\centering
        			\scalebox{0.5}{
        				\begin{tikzpicture}[shorten >=1pt,auto,node distance=3cm,
        				thick,main node/.style={circle,fill=blue!20,draw,font=\sffamily\LARGE\bfseries}, three/.style={circle split, fill=blue!20,draw,font=\sffamily\LARGE\bfseries}]
				
        				\node[main node] (1) at (0,4)  {$x_1$};
        				\node[main node] (2) at (3.46,2)  {$x_2$};
        				\node[main node] (3) at (3.46,-2)  {$x_3$};
        				\node[three] (4) at (0,-4)  {$x_4$ \nodepart{lower} $x_6$};
        				\node[main node] (7) at (-3.46,-2)  {$x_7$};
        				\node[main node] (8) at (-3.46,2)  {$x_8$};
        				
        				\path[every node/.style={font=\sffamily},line width=1pt]
        				(1) edge [orange] node [] {} (2)
        				(2) edge [orange] node [] {} (3)
        				(3) edge [blue,dotted] node [] {} (4)
        				(4) edge [blue,dotted] node [] {} (7)
        				(7) edge [black,dashed] node [] {} (8)
        				(8) edge [black,dashed] node [] {} (1);
        				\path[every node/.style={font=\sffamily},line width=2pt]
        				(1) edge node [] {} (3)
        				(1) edge node [] {} (7);
        				\end{tikzpicture}
        			}
        			\caption{Figure \ref{fig:pm1-b}: $C(\phi')$.}
        			\label{fig:pm1-b}
        		\end{subfigure}
        	\caption{An illustration of the contraction argument.
Figure~\ref{fig:pm1-a} illustrates a constraint graph with a unique vertex $x_5$.
This vertex is indeed unique as there are no constraint edge incident to this vertex. Figure~\ref{fig:pm1-b} illustrates the contraction process. After the contraction, $x_5$ disappears and the neighboring vertices, namely $x_4$ and $x_6$, are merged into one. 
The resulting graph still satisfies the property that every edge value appears at least twice, allowing us to use the inductive hypothesis.}
        	\label{P13}
        	\end{figure}

		\end{enumerate}		
		Combining the two cases above, the inductive step is completed, and hence Lemma~\ref{mink-1} follows.
	\end{proof}

			We now prove the second statement.
			\begin{proposition} \label{numconstr}
				For all $k \in \mathbb{N}$, there are at most $2^{|V(\alpha,2q)|-1}(2q)^k = 2^{2q-1}(2q)^k$ constraint graphs $C \in \mathcal{C}_{(\alpha,2q)}$ which have exactly $k$ constraint edges.
			\end{proposition}
			\begin{proof}
			To count the number of constraint graphs with exactly $k$ constraint edges, we encode each constraint graph as follows.
            Let $C$ be a given constraint graph with $k$ constraint edges.
            Since we have identified the constraint graph with cycles, we can naturally order the  vertices according to their indices, i.e., $x_1,x_2,\dots, x_{2q}$.
            Having ordered the vertices, define 
				\[
				V_{\mathrm{redundant}} := \{ x_i\in V(C_{2q}): \exists x_j \in V(C_{2q}): \phi(x_i) = \phi(x_j), ~~j<i\}\,,
				\]
				where $\phi: V(\alpha,2q) \to [n]$ is any piecewise injective map such that $C(\phi)=C$. Notice that $V_{\mathrm{redundant}}= k$ as there are $k$ constraint edges in $C$.
				
				Now we encode the constraint graph $C$ with the following data:
				\begin{enumerate}
					\item For each vertex $x_i \in V(C_{2q})$, is $x_i \in V_{\mathrm{redundant}}$? 
					\item For each vertex $x_i \in V_{\mathrm{redundant}}$, what is the vertex $x_j \in V(C_{2q})$ which comes before $x_i$ such that $\phi(x_i) = \phi(x_j)$? Here if there are multiple such vertices, we may pick any of them.
				\end{enumerate}
			Having defined the encoding scheme, one can easily observe that this encoding is valid, i.e., each encoding corresponds to a unique constraint graph $C$.
			
			Now, based on this encoding scheme, one can easily count the number of constraint graphs with $k$ constraint edges.
				Since the first vertex cannot be in $V_{\mathrm{redundant}}$, there are at most $2^{2q-1}$ possibilities for which vertices are in $V_{\mathrm{redundant}}$. For each of the $k$ vertices $x_i \in V_{\mathrm{redundant}}$, there are at most $2q$ choices for the vertex $x_j$ which comes before $x_i$ such that $\phi(x_i) = \phi(x_j)$. 
				
				Putting everything together, there are at most $2^{2q-1}(2q)^k$ constraint graphs $C \in \mathcal{C}_{(\alpha,2q)}$ which have exactly $k$ constraint edges, as needed.
			\end{proof}
		We have completed the proofs of the two statements, which completes the proof of Lemma \ref{traceupper}.
		\end{proof}
		Having established Lemma \ref{traceupper}, one can now prove Theorem \ref{singular main} using the trace power method. By the trace power method (Lemma \ref{basictracepower}), for all $\epsilon > 0$,
		\[
		\pP\left[\norm{M_{\alpha}} > \min_{q \in \mathbb{N}: q \leq \frac{n}{4}}{\left\{\sqrt[2q]{\frac{\B{2q}}{\epsilon}}\right\}}\right] < \epsilon
		\]
		where $\B{2q} = 2^{2q}(2q)^{q}n^{q+1}$. The derivative of
		\[
		\ln\left(\sqrt[2q]{\frac{\B{2q}}{\epsilon}}\right) = \ln\left(2\sqrt{2n}\sqrt{q}\sqrt[2q]{\frac{n}{\epsilon}}\right) = \ln(2\sqrt{2n}) + \frac{\ln(q)}{2} + \frac{\ln(\frac{n}{\epsilon})}{2q}
		\]
		with respect to $q$ is $\frac{1}{2q} - \frac{\ln(\frac{n}{\epsilon})}{2q^2}$ which is only $0$ at $q = \ln(\frac{n}{\epsilon})$. Since the second derivative is positive at $q = \ln(\frac{n}{\epsilon})$, this is a global minimum for our bound. However, we require that $q \in \mathbb{N}$, so we instead take $q = \lceil{\ln(\frac{n}{\epsilon})}\rceil$. 
		Taking $q = \lceil{\ln(\frac{n}{\epsilon})}\rceil$, \[
		\sqrt[2q]{\frac{\B{2q}}{\epsilon}} = 2\sqrt{2qn}\left(\frac{n}{\epsilon}\right)^{\frac{1}{\lceil{2\ln(\frac{n}{\epsilon})}\rceil}} \leq 
		2\sqrt{2e{\left\lceil{\ln\left(\frac{n}{\epsilon}\right)}\right\rceil}n}\,.
		\]
		Note that if $q = \lceil{\ln(\frac{n}{\epsilon})}\rceil \geq \frac{n}{4}$ then $2\sqrt{2qn} > n$ but $\norm{M_{\alpha}}$ is always at most $n$ so our bound holds in this case as well.
	\end{proof}
	
	\section{Modified techniques for graph matrices}\label{sec:modified}

     	In this section, we modify the techniques in Section~\ref{pm1} so that they can be applied to more general shapes $\alpha$.
	Due to technical reasons, for many of the matrices $M$ which we analyze, it is difficult to analyze $\E\left[\tr\left((M{M^\top})^q\right)\right]$ directly. Instead, we break $M$ up into pieces which are easier to analyze. 
To that end, we consider a vertex partitioning argument which we use to modify the trace power method.	

	\subsection{Vertex partitioning}
	
	\begin{definition}\label{def:alpha}
		Given a shape $\alpha$, we define an $\alpha$-partition of the indices to be a tuple $P = (I_v: v \in V(\alpha))$ such that 
		\begin{enumerate}
			\item $\forall v \in V(\alpha), I_v \subseteq [n]$.
			\item $\forall u \neq v \in V(\alpha), I_{u} \cap I_{v} = \emptyset$.
			\item $\cup_{v \in V(\alpha)}{I_{v}} = [n]$
		\end{enumerate}
		We define $\mathcal{P}_{\alpha}$ to be the set of all possible $\alpha$-partitions of the indices.
	\end{definition}
	\begin{proposition} \label{cardinality}
		$|\mathcal{P}_\alpha|= |V(\alpha)|^n$.
	\end{proposition}
	\begin{proof}
		The proposition follows from the fact that each element in $[n]$ can be assigned to $|V(\alpha)|$ different $I_v$'s.
	\end{proof}
	\begin{definition}
		Given an $\alpha$-partition $P$ of the indices, we say that a map $\phi:V(\alpha) \to [n]$ is $P$-respecting if $\forall v \in V(\alpha), \phi(v) \in I_v$.
	\end{definition}
	
	\begin{proposition} \label{respectcount}
		For any injective map $\phi: V(\alpha) \to [n]$, there are exactly $|V(\alpha)|^{n - |V(\alpha)|}$ $\alpha$-partitions $P \in \mathcal{P}_{\alpha}$ such that $\phi$ is $P$-respecting.
	\end{proposition}
	\begin{proof}
		Consider an integer $i\in[n]$. If $i=\phi(v)$ for some $v\in V(\alpha)$, then we have $i\in I_v$ by definition.
		Otherwise, i.e., $i\not\in \phi(V(\alpha))$, one can assign  $i$ to an arbitrary $I_v$, which amounts to $|V(\alpha)|$ different possibilities.
		Since there are $n-|V(\alpha)|$ elements outside the image of $\phi$, the conclusion follows.
	\end{proof}
	
	\begin{definition}
		We define $M_{\alpha,P}$ to be the matrix with entries \[
		M_{\alpha,P}(A,B) = \sum_{\phi: V(\alpha) \to [n]: \phi \textnormal{ is } P \text{-respecting, } \phi(\U{\alpha}) = A, \phi(\V{\alpha}) = B}{\chi_{\phi(E(\alpha))}(G)}\,.
		\]
	\end{definition}
	The following result is an immediate consequence of Propositions~\ref{cardinality} and~\ref{respectcount}:
	\begin{corollary} \label{repre}
		$M_{\alpha} = \frac{|V(\alpha)|^{|V(\alpha)|}}{|\mathcal{P}_\alpha|}\sum_{P \in \mathcal{P}_\alpha}{M_{\alpha,P}}$.
	\end{corollary}

	\subsection{Constraint graphs and vertex partitioning}\label{constraintgraph}
	We can still use constraint graphs to analyze the matrices $\M_{\alpha,P}$.
	To that end, we first consider an expansion similar to Proposition~\ref{expansion}. One difference is that due to the definition of $\M_{\alpha,P}$ the summation is over only $\phi$'s that \emph{respect} the partition $P$:
	\begin{definition}
		Given a shape $\alpha$, $q \in \mathbb{N}$, and $P \in \mathcal{P}_{\alpha}$, we say that a map $\phi: V(\alpha,2q) \to [n]$ is $P$-respecting if $\phi$ is $P$-respecting on each piece $V(\alpha_i)$ and each piece $V(\alpha^{\top}_i)$ for all $i \in [q]$. In other words, for all vertices $v \in V(\alpha)$ and all $i \in [q]$, letting $v_i \in V(\alpha_i)$ be the copy of $v$ in $V(\alpha_i)$ and letting $v'_i \in V(\alpha^\top_i)$ be the copy of $v$ in $V(\alpha^\top_i)$, $v_i, v'_i \in I_v$.
	\end{definition}
	\begin{remark}
		Note that if a map $\phi: V( \alpha,2q) \to [n]$ is $P$-respecting, then it is automatically piecewise injective. Also note that vertex partitioning does not affect $\val{C}$.
	\end{remark}
	The above definition yields the following  analogue of Proposition~\ref{expansion} for $\M_{\alpha, P}$: 
	\begin{proposition}
		For all shapes $\alpha$, all $P \in \mathcal{P}_{\alpha}$, and all $q \in \mathbb{N}$,
		\[
		\E\left[\tr\left((M_{\alpha,P}M_{\alpha,P}^\top)^q\right)\right] = \sum_{\phi: V(\alpha,2q) \to [n]:\phi \text{ is } P \text{-respecting}}{\E[\chi_{\phi(E(\alpha,2q))}(G)]}\,.
		\]
	\end{proposition}
	
	In fact, the fact that the summation is over $P$-respecting maps makes the analysis nicer because the constraint graphs corresponding to those maps are well-behaved:
	\begin{definition}
		We say that a constraint graph $C \in \mathcal{C}_{(\alpha,2q)}$ is well-behaved if for every edge $(u,w) \in E(C)$, $u,w$ are copies of the same vertex $v \in V(\alpha)$ (or its mirror in $V(\alpha^\top)$).
	\end{definition}
	Similarly to Proposition~\ref{diffrep},  we represent the summation in terms of constraint graphs by counting the number of maps $\phi$ corresponding to each constraint graph $C$.
	
	\begin{definition}
		Given a constraint graph $C \in \mathcal{C}_{(\alpha,2q)}$, we define 
		\[
		N_P(C) = |\{\phi: V( \alpha,2q) \to [n]: \phi \text{ is } P \text{-respecting}, C(\phi) \equiv C\}|\,.
		\]
	\end{definition}
	Since every constraint graph corresponding to a $P$-respecting map $\phi: V(\alpha,2q) \to [n]$ is well-behaved, the following holds:
	\begin{proposition}
		If $C \in \mathcal{C}_{(\alpha,2q)}$ is not well-behaved then for all $P \in \mathcal{P}_{\alpha}$, $N_{P}(C) = 0$.
	\end{proposition}
	Consequently, the summation is over only well-behaved constraints, which implies the following representation of the trace power:
	\begin{corollary} \label{expan2}
		For all shapes $\alpha$, all $P \in \mathcal{P}_{\alpha}$, and all $q \in \mathbb{N}$,
		\[
		\E\left[\tr\left((M_{\alpha,P}M_{\alpha,P}^\top)^q\right)\right] = \sum_{C \in \mathcal{C}_{(\alpha,2q)}~:~ C \text{ is well-behaved}}{N_{P}(C)\val{C}}\,.
		\]
	\end{corollary}
	\subsection{Vertex partitioning and the trace power method}
	One way to obtain a norm bound on $M_{\alpha}$ would be to use the norm bounds on $M_{\alpha,P}$ and then use a union bound to upper bound the probability that some norm bound fails. Unfortunately, there are too many matrices $M_{\alpha,P}$ for this naive analysis to work. However, with a careful modification of the trace power method we can obtain probabilistic norm bounds on $M_{\alpha}$ directly from upper bounds on $\E\left[\tr\left((M_{\alpha,P}M_{\alpha,P}^\top)^q\right)\right]$.
\begin{mdframed}[style=box]
	\begin{lemma}[Vertex partitioning lemma]\label{vpl}
		If $M$ is a matrix such that $M = \frac{1}{K}\sum_{i=1}^{K}{M_i}$ for some $K \in \mathbb{N}$ and matrices $\{M_i: i \in [K]\}$ and $B: \mathbb{N} \to \mathbb{R}$ is a function such that for all $i \in [K]$ and all $q \geq 3$, $\E\left[\tr\left((M_i {M_i}^\top)^q\right)\right] \leq B(2q)$
		then for all $\epsilon > 0$, $$\pP\left[\norm{M} > \min_{ q \geq 3}\left\{2\sqrt[2q]{\frac{B(2q)}{\epsilon}}\right\}\right] < \epsilon\,.$$
	\end{lemma}
\end{mdframed}
	\begin{proof}
		We can prove this lemma as follows:
		\begin{enumerate}
			\item Use the trace power method to obtain probabilistic norm bounds on the matrices $\{M_i: i \in [K]\}$.
			\item Combine these probabilistic norm bounds on the matrices $\{M_i: i \in [K]\}$ to obtain the probabilistic norm bound on $M$.
		\end{enumerate}
		We start with step 2 as this determines what probabilistic norm bounds we need on the matrices $\{M_i: i \in [K]\}$. For step 2, we use the following lemma.
		\begin{lemma}\label{combiningboundslemma}
			If $M$ is a matrix such that $M = \frac{1}{K}\sum_{i=1}^{K}{M_i}$ for some $K \in \mathbb{N}$ and matrices $\{M_i: i \in [K]\}$ and $B \in \mathbb{R}$ is a bound such that for all $i \in [K]$ and all $x \in [\frac{1}{2},K], \pP[\norm{M_i} > Bx] \leq \frac{\epsilon}{32x^3}$
			then $\pP[\norm{M} \geq B] < \epsilon$.
		\end{lemma}
		\begin{proof}[Proof of Lemma~\ref{combiningboundslemma}]
			To prove this, we make the following observation. 
			\begin{prp} \label{int:prp}
				For all $j \in [0,\log_2(K)+1]$, $\pP\left[|\{{i \in [K]:\norm{M_i} > 2^{j-1}}B\}| > \frac{K}{2^{2j+1}}\right] \leq \frac{\epsilon}{2^{j+1}}$.
			\end{prp}
			\begin{proof}[Proof of Proposition~\ref{int:prp}]
				We prove this by contradiction. If $\pP\left[|\{{i \in [K]:\norm{M_i} > 2^{j-1}}B\}| > \frac{K}{2^{2j+1}}\right] > \frac{\epsilon}{2^{j+1}}$, then 
				$\E\left[|\{{i \in [K]:\norm{M_i} > 2^{j-1}}B\}|\right] > \frac{K}{2^{2j+1}}\cdot\frac{\epsilon}{2^{j+1}} =  \frac{{\epsilon}K}{2^{3j+2}}$. However, plugging in $x = 2^{j-1}$ we have that 
				\[
				\E\left[|\{{i \in [K]:\norm{M_i} > 2^{j-1}}B\}|\right] \leq K\max_{i}\left\{\pP[\norm{M_i} > Bx]\right\} \leq \frac{{\epsilon}K}{32x^3} = \frac{{\epsilon}K}{2^{3j+2}}
				\] 
				which gives a contradiction.
			\end{proof}
\noindent			Using this proposition, with probability at least $1 - \sum_{j=0}^{\lfloor{\log_2(K) + 1}\rfloor}{\frac{\epsilon}{2^{j+1}}} > 1 - \epsilon$,
			\[
			\forall j \in [0,\log_2(K) + 1],~~ |\{i: ||M_i|| > 2^{j-1}B\}| \leq \frac{K}{2^{2j+1}}\,.
			\]
		Writing the above statement in a different way, we have
		\begin{enumerate}
		    \item $\forall j \in [0,\lceil{\log_2(K)}\rceil],~~ |\{i: 2^{j-1}B < \norm{M_i} \leq 2^{j}B\}| \leq \frac{K}{2^{2j+1}}$ and
		    \item there are  no matrices $M_i$ such that $\norm{M_i} > 2^{\lceil{\log_2(K)}\rceil}B$.
		 \end{enumerate}
		    This implies that 
			\[
			\norm{M} \leq \sum_{j=0}^{\lceil{\log_2(K)}\rceil}{\frac{1}{2^{2j+1}} \cdot 2^{j}B} = \sum_{j=0}^{\lceil{\log_2(K)}\rceil}{\frac{B}{2^{j+1}}} <  B\,.
			\]
			Thus, $\pP[\norm{M} \geq B] < \epsilon$, as needed.
		\end{proof}
		With Lemma \ref{combiningboundslemma} in mind, we now obtain probabilistic norm bounds on the matrices $\{M_i: i \in [K]\}$ which hold with probability at least $1 - \frac{\epsilon}{32x^3}$. 
		
		Assume that $q \geq 3$, $\epsilon > 0$, and $x \geq \frac{1}{2}$. Applying the trace power method (Lemma \ref{basictracepower}) with $\epsilon' = \frac{\epsilon}{32x^3}$, we obtain:
		\[
		\pP\left[\norm{M_i} > \sqrt[2q]{\frac{B(2q)}{\epsilon'}}\right] < \epsilon'\,.
		\]
		Now take $B = 2\sqrt[2q]{\frac{B(2q)}{\epsilon}}$. Since $q \geq 3$, 
		$\sqrt[2q]{\frac{B(2q)}{\epsilon'}} = \sqrt[2q]{32x^3\frac{B(2q)}{\epsilon}} < Bx$. Thus, for all $i \in [K]$ and $x \geq \frac{1}{2}$,
		$\pP\left[\norm{M_i} > Bx\right] < \frac{\epsilon}{32x^3}$.
		
		Plugging these bounds into Lemma \ref{combiningboundslemma}, we obtain $\pP\left[\norm{M} \geq 2\sqrt[2q]{\frac{B(2q)}{\epsilon}}\right] < \epsilon$. Since this holds for all $q \geq 3$ and $\epsilon > 0$, the result follows after taking the minimum over $q\geq 3$.
	\end{proof}
	\begin{remark}[Alternative vertex partitioning arguments]
		An alternative way to obtain a probabilistic norm bound on $M$ is as follows. Instead of expressing $M$ as the average of $K$ matrices, we ``cover'' $M$ by $m = \OO{\log(n)}$ matrices $\{M_i: i \in [m]\}$. More precisely, we take each $M_i$ to be a matrix of the form $M_{\alpha,P}$ except that we zero out any entry which was already included in a previous matrix. Thus, $M = \sum_{i=1}^{m}{M_i}$. We can then use the fact that $\norm{M} \leq m\cdot \max_{i \in [m]}{\norm{M_i}}$ and use a union bound to bound the probability that any of the norm bounds on the matrices $\{M_i: i \in [m]\}$ fail.
		Indeed, this argument was employed in prior works (e.g. see the beginning of Section 4 in \cite{barak2016noisy}).
		
		While this alternate approach is simpler than the given approach in some ways, the bound it gives is worse by a logarithmic factor.
	\end{remark}

	\section{Proving norm bounds on graph matrices (Theorem~\ref{thm:mainresult})} \label{pf:main}

	In this section, we prove our main theorem (Theorem~\ref{thm:mainresult}), building on the proof for symmetric $\pm 1$ Wigner matrices in Section~\ref{pm1}. 
	We first give the formal statement of Theorem~\ref{thm:mainresult}:
	\begin{mdframed}[style=box2]
			\begin{theorem}[Formal statement of Theorem~\ref{thm:mainresult}] \label{thm:formal:mainresult}
		For a given shape $\alpha = (\U{\alpha},\V{\alpha},\W{\alpha},E(\alpha))$, let $\smin$ be a minimum vertex separator between  $\U{\alpha}$ and $\V{\alpha}$, and  $\iso\subset \W{\alpha}$ be the subset of isolated middle vertices.
		Based on this, taking 	$\ccal:=|\W{\alpha}| + |\smin| -|\iso| -|\U{\alpha} \cap \V{\alpha}|$ and $\dal:= |V(\alpha)|-|\iso|-|\U{\alpha}\cap \V{\alpha}|$, the following bounds hold:
			\begin{enumerate}
			    \item (Probabilistic bounds for $\ccal>0$) If $\ccal>0$, for all $\epsilon>0$, we have
			\begin{align*} 
	 \pP\left[\norm{\M_\alpha} \geq   \left( 6e \left\lceil \frac{1}{3\ccal}\log(n^{|\smin| }/\epsilon)\right\rceil\right)^{\frac {1}{2}\ccal}\cdot  (2\dal)^{\dal}  n^{\frac{1}{2}(|V(\alpha)|+|\iso|-|\smin|)}\right] \leq \eps \,.
			\end{align*}
			\item (Trivial deterministic bounds) $\norm{\M_\alpha} \leq n^{\frac{1}{2}(|V(\alpha)|+|\W{\alpha}|-|\U{\alpha}\cap \V{\alpha}|)}$.
			\end{enumerate}
			\end{theorem} 
		\end{mdframed}
We will first restrict ourselves to the special case where the left and right vertices are disjoint, i.e., $\U{\alpha}\cap \V{\alpha}=\emptyset$ and there are no isolated middle vertices, i.e., $\iso =\emptyset$.
This special case delivers the key ideas of the proof while avoiding unnecessary complications.
Later in this section, we will discuss how we can relax the restrictions and prove the full results in Theorem~\ref{thm:formal:mainresult}.

	\subsection{Proof for the special case:  \texorpdfstring{$\U{\alpha}\cap \V{\alpha}=\emptyset$}{+1} and \texorpdfstring{$\iso =\emptyset$}{+1}} \label{sec:simple}
 We begin with proving the trivial deterministic bound:
\begin{proof}[Proof of the trivial deterministic bound]
Since we are considering the case when $U_{\alpha} \cap V_{\alpha} = \emptyset$, we need to show that $\norm{M_{\alpha}} \leq n^{\frac{1}{2}(|V(\alpha)|+|\W{\alpha}|)}$.
To prove this, we will use a simple fact that an $n_1\times n_2$ matrix with the magnitude of each entry bounded by $c$ has spectral norm at most $c \cdot \sqrt{n_1n_2}$.
Indeed, this fact is an easy consequence of the inequality $\norm{M}\leq \norm{M}_F$ for any matrix $M$ where $\norm{M}_F$ is the Frobenius norm of $M$.
Now it is straightforward from the definition of graph matrices (Definition~\ref{def:graphmatrix}) that there are at most $n^{|\U{\alpha}|}$ rows and at most $n^{|\V{\alpha}|}$ columns in $\M_\alpha$.
Moreover, each entry of $\M_\alpha$ is the summation of at most $n^{|\W{\alpha}|}$ $\pm{1}$ random variables, so the magnitude of each entry is bounded by $n^{|\W{\alpha}| }$.
	Hence, from the simple fact above, the trivial deterministic bound follows:
$\norm{\M_\alpha} \leq n^{|\W{\alpha}|}\cdot \sqrt{n^{|\U{\alpha}|}\cdot n^{|\V{\alpha}|}} = n^{\frac{1}{2} (|V(\alpha)| +|\W{\alpha}|)}$.
	\end{proof}

Now we prove the probabilistic bound.
Before the proof, we make one notation remark:
\begin{itemize}
    \item  Here and below, we simply write $s:=|\smin|$.
\end{itemize}Note that for this special case the probabilistic bound reads:
\begin{align} \label{bound:int}
    \pP\left[  \norm{\M_\alpha} \geq  \left(6e \left\lceil \frac{1}{3(|\W{\alpha}| + s)}\ln(n^{s}/\epsilon)\right\rceil\right)^{\frac{1}{2}(|\W{\alpha}| + s )}\cdot (2|V(\alpha)|)^{|V(\alpha)|}  n^{\frac{1}{2}(|V(\alpha)|-s)} \right] \leq \epsilon
\end{align}

The proof of \eqref{bound:int} bears great similarities to the proof for the $\pm1$ random matrix.
One main distinction is that we now decompose $\M_\alpha$ into \emph{well-behaved} pieces and conquer each piece separately.
More formally, we decompose $\M_\alpha$ as follows in light of Corollary~\ref{repre}: 
\begin{align}
    M_{\alpha} = \frac{|V(\alpha)|^{|V(\alpha)|}}{|\mathcal{P}_\alpha|}\sum_{P \in \mathcal{P}_\alpha}{M_{\alpha,P}}\,. \label{decomp}
\end{align}

 Now, similar to the $\pm1$ random matrix case, the key step is to establish upper bounds on the terms $\E\left[\tr\left((M_{\alpha,P}M_{\alpha,P}^\top)^q\right)\right]$.
 Indeed, one can turn upper bounds on the expected trace power into a probabilistic upper bound on $\norm{\M_\alpha}$ via the vertex partitioning lemma (Lemma~\ref{vpl}).
As we shall see later, it turns out the following analogue of Lemma~\ref{traceupper} suffices:
\begin{lemma}\label{traceupper2}
	For all $q \in \mathbb{N}$ such that $q \leq \frac{n}{4}$,
	\begin{align*}
			  \E\left[\tr\left((M_{\alpha,P}M_{\alpha,P}^\top)^q\right)\right]  \leq 2^{2q|V(\alpha)|}(2q)^{q(|\W{\alpha}|+s)}n^{q(|V(\alpha)|-s) + s}\,.
		\end{align*}
\end{lemma}
\begin{proof}

Paralleling the proof of Lemma~\ref{traceupper}, the following statements are main ingredients to proving Lemma~\ref{traceupper2}:
\begin{enumerate}
  	\item For any well-behaved $C \in \mathcal{C}_{(\alpha,2q)}$ such that $\val{C} \neq 0$, $|E(C)| \geq  q|\W{\alpha}| + (q-1)s $. 
    \item 	For all $k$, there are at most  $2^{|V(\alpha,2q)| - 1}(2q)^{k}$ well-behaved constraint graphs $C \in \mathcal{C}_{(\alpha,2q)}$ which have exactly $k$ constraint edges. 
\end{enumerate}
Assuming these statements are true, the proof of Lemma~\ref{traceupper2} is immediate.
To see this, first recall that by  Corollary~\ref{expan2},
	\begin{align*}
		\E\left[\tr\left((M_{\alpha,P}M_{\alpha,P}^\top)^q\right)\right] &= \sum_{C \in \mathcal{C}_{(\alpha,2q)}:~ C \text{ is well-behaved}}{N_{P}(C)\val{C}}\\
		&\leq \sum_{\substack{C \in \mathcal{C}_{(\alpha,2q)}:~ C \text{ is well-behaved}\\
		\val{C}\neq 0}} n^{|V(\alpha,2q)|-|E(C)|}\,,
			\end{align*}
where the inequality follows from Proposition~\ref{numchoice} together with the simple fact  $N_P(C) \leq N(C)$.
Now, we rearrange the summation in terms of $|E(C)|$.
Due to the first statement,  $|E(C)|$ should be greater than or equal to  $q|\W{\alpha}| + (q-1)s$ and at most $|V(\alpha,2q)|-1$.
Moreover, the summand for  $|E(C)|=k$ is at most $2^{|V(\alpha,2q)| - 1}(2q)^{k}\cdot n^{|V(\alpha,2q)|-k}$ due to the second statement.
Thus, the trace power term can be upper bounded as follows:
\begin{align*}
	\E\left[\tr\left((M_{\alpha,P}M_{\alpha,P}^\top)^q\right)\right] 
			&\leq \sum_{k =q|\W{\alpha}| + (q-1)s}^{|V(\alpha,2q)|-1}{2^{|V(\alpha,2q)| - 1}(2q)^{k}n^{|V(\alpha,2q)|-k}}\\
		&\leq  2^{|V(\alpha,2q)| - 1}n^{|V(\alpha,2q)|} \left(\frac{2q}{n}\right)^{q|\W{\alpha}| + (q-1)s}  \cdot \sum_{k =0}^{\infty}{\left(\frac{2q}{n}\right)^{k}}\\
	&\leq 2^{|V(\alpha,2q)|} n^{|V(\alpha,2q)|} \left(\frac{2q}{n}\right)^{q|\W{\alpha}| + (q-1)s}  \,,
	\end{align*}
where the last inequality is due to the assumption $q\leq \frac{n}{4}$. 
It is then straightforward to check that the last upper bound is upper bounded by the bound in Lemma~\ref{traceupper2} using Proposition~\ref{prop:num} which says that $|V(\alpha,2q)| = q|\U{\alpha} \cup \V{\alpha}|   + 2q|\W{\alpha}|$.
		
Thus, we just need to prove these two statements.
We begin with the first statement.
For the first statement, we need to recall classical results from graph theory: K\"{o}nig's Theorem~\cite{koniggrafok} and Menger's Theorem~\cite{Menger}.
For completeness, we first state them here:
  
	\begin{proposition}[K\"{o}nig's Theorem]\label{Konigstheorem}
	Given a graph $G$, a vertex cover of $G$ is a set of vertices $V \subseteq V(G)$ such that all edges of $G$ are incident with at least one vertex in $V$.
		If $G$ is a bipartite graph with partite sets $U$ and $V$ then the minimal size of a vertex cover of $G$ is equal to the maximal size of a matching between $U$ and $V$. 
	\end{proposition} 
	\begin{proposition}[Menger's Theorem]\label{Menger}
	If $G$ is a graph and $U,V \subseteq V(G)$, we say $S$ is  a vertex separator between $U$ and $V$ if  all paths from $U$ to $V$ intersect $S$.
		If $G$ is a graph and $U,V \subseteq V(G)$ then the minimum size of a vertex separator between $U$ and $V$ is equal to the maximal number of vertex disjoint paths between $U$ and $V$.
	\end{proposition}
Now we are ready to prove the first statement:

\begin{lemma} \label{mingeneral}
				 For any well-behaved $C \in \mathcal{C}_{(\alpha,2q)}$ such that $\val{C} \neq 0$, $|E(C)| \geq  q|\W{\alpha}| + (q-1)s $.  
				As a special case, if $\alpha$ is bipartite, then $|E(C)|\geq (q-1)s$.
			\end{lemma}
			\begin{proof}
			 As a warmup, we begin with the case where $\alpha$ is bipartite. In this case, $s$ is the minimum size of a vertex cover of $\alpha$.

				By K\"{o}nig's Theorem (Proposition~\ref{Konigstheorem}), there is a matching of size $s$ between $\U{\alpha}$ and $\V{\alpha}$ in $\alpha$. 
				Choose an arbitrary edge $\{u,v\}$ from the matching and consider the cycle $T=C_{2q}$  (see Definition~\ref{def:cycle}) formed by the  copies of $u$ and $v$ in $C$.
			Then, there is no  constraint edge between $T$ and $C\setminus C_{2q}$ because the constraint graph is well-behaved.
				This implies that edges in $\phi(E(T))$ do not appear in $\phi(E(\alpha,2q)\setminus E(T))$.
				Hence, each edge in $\phi(E(T))$ appears at least twice within $\phi(E(T))$.
				Applying Lemma~\ref{mink-1} to $T$, there should be at least $q-1$ constraint edges between vertices of $T$.
				
				Since $\{u,v\}$ was arbitrarily chosen,  one can apply this  argument to all of the edges in the perfect matching, which gives the desired lower bound of $(q-1)s$.
				
				For general $\alpha$, we can use similar ideas. By Menger's Theorem (Proposition~\ref{Menger}), there are $s$ vertex-disjoint paths $P_1,\ldots,P_{s}$ from $\U{\alpha}$ to $\V{\alpha}$ such that for each of these paths $P_i$:
				\begin{enumerate}
					\item The length of $P_i$ is at least $1$.
					\item $P_i$ only intersects $\U{\alpha}$ in its first vertex and $P_i$ only intersects $\V{\alpha}$ in its last vertex. In other words, there is no subpath of $P_i$ from $\U{\alpha}$ to $\V{\alpha}$.
				\end{enumerate}
				Let $l_1,\ldots,l_{s}$ be the lengths of these paths. Letting $T_1,\ldots,T_{s}$ be the vertex disjoint cycles where $T_i$ consists of all of the copies of the path $P_i$, $T_i$ has length $2q{l_i}$ and following similar logic, there must be at least $q{l_i} - 1$ constraint edges within $T_i$. Thus, there must be at least $q\left(\sum_{i=1}^{s}{l_i}\right) - s$ constraint edges within these cycles.
				
				In addition, there are $2q\left(|\W{\alpha}| - \sum_{i=1}^{s}{(l_i - 1)}\right)$ vertices in $V(\alpha,2q)$ which are copies of vertices in $\W{\alpha}$ but are not part of any of these cycles. 
				Observe that each of these vertices must be incident to at least one constraint edge as otherwise there would be an edge in $\phi(V(\alpha))$ which only appears once. Thus, there must be at least $q\left(|\W{\alpha}| - \sum_{i=1}^{s}{(l_i - 1)}\right)$ additional constraint edges.
				
				Putting everything together, $|E(C)| \geq q|\W{\alpha}| + (q-1)s$.
			\end{proof} 
 	 
	We now prove the second statement.
	\begin{proposition}\label{numconstr2}
				For all $k$, there are at most  $2^{|V(\alpha,2q)| - 1}(2q)^{k}$ well-behaved constraint graphs $C \in \mathcal{C}_{(\alpha,2q)}$ which have exactly $k$ constraint edges.
			\end{proposition}
			\begin{proof}
				The proof  is very similar to the proof of Lemma~\ref{numconstr}. To specify the constraint graph, it is sufficient to specify the following data:
				\begin{enumerate}
					\item The set of vertices $V_{\mathrm{redundant}} \subseteq V(\alpha,2q)$ which are equal to a previously seen vertex.
					\item For each vertex $v \in V_{\mathrm{redundant}}$, which previously seen vertex $u \in V(\alpha,2q)$ is $v$ equal to?
				\end{enumerate}
				The first vertex we consider cannot be in $V_{\mathrm{redundant}}$ and each vertex we consider afterwards is either in $V_{\mathrm{redundant}}$ or it is not. Thus, there are at most $2^{|V(\alpha,2q)|-1}$ choices for which vertices are in $V_{\mathrm{redundant}}$.
				
				Since the constraint graph $C$ must be well-behaved, for each vertex $v \in V_{redundant}$, there are at most $2q$ choices for which previous vertex $v$ is equal to. Since $|V_{\mathrm{redundant}}| = k$, there are at most $(2q)^{k}$ choices for the previous vertices.
				
				Putting everything together, there are at most $2^{|V(\alpha,2q)| - 1}(2q)^{k}$ possibilities, as needed.
			\end{proof}

		This completes the proof of Lemma~\ref{traceupper2}.  				\end{proof}

	Having established the upper bounds on the expected trace power terms, one can turn them into a probabilistic norm bound on $\norm{\M_\alpha}$ using the vertex partitioning lemma (Lemma~\ref{vpl}).
	In particular, due to \eqref{decomp} and Lemma~\ref{traceupper2}, one can apply the vertex partitioning lemma (Lemma~\ref{vpl}) with \begin{enumerate}
	    \item $M\leftarrow \frac{1}{|V(\alpha)|^{|V(\alpha)|}}\M_\alpha$,
	    \item $\{M_i\}\leftarrow \{\M_{\alpha,P}\}$, and 
	    \item $B(2q) \leftarrow 2^{2q|V(\alpha)|}(2q)^{q(|\W{\alpha}|+s)}n^{q(|V(\alpha)|-s) + s}$
	\end{enumerate}
	This gives the following for all $\epsilon>0$:
$$\pP\left[\norm{\M_\alpha} >  (2|V(\alpha)|)^{|V(\alpha)|}   n^{\frac{1}{2}(|V(\alpha)|-s)}2^{\frac{1}{2}(|\W{\alpha}|+s)}\min_{ 3\leq q \leq \frac{n}{4}}\left\{ q^{\frac{1}{2}(|\W{\alpha}|+s)} (n^{s}/\epsilon)^{1/2q}  \right\}\right] < \epsilon\,.$$

Now let us estimate the norm bound by computing the minimum value.
A similar derivative calculation concludes that the upper bound acheives the minimum at $q^\star:=\frac{1}{|\W{\alpha}|+s}\ln(n^s/\epsilon)$.
To ensure that our choice of $q$ is an integer at least $3$, let us opt for a slightly suboptimal choice  $q=3\lceil q^\star/3 \rceil$; this will surely guarantee $q\geq 3$ since $3\lceil a/3 \rceil \geq 3$ for any $a>0$.

	Let us first assume that $ 3\lceil q^\star/3 \rceil\leq n/4$, then choosing $q=3\lceil q^\star/3 \rceil$ in the minimum term, we obtain:
		\begin{itemize}
			\item $q^{\frac{1}{2}(|\W{\alpha}|+s)} = \left( 3\left\lceil\frac{1}{3(|\W{\alpha}|+s)}\ln(n^s/\epsilon)\right\rceil\right)^{\frac{1}{2}(|\W{\alpha}|+s)} $.
			\item $(n^{s}/\epsilon)^{1/2q}
			\leq (n^{s}/\epsilon)^{\frac{|\W{\alpha}|+s}{2\ln(n^s/\epsilon)}}  = e^{\frac{1}{2}(|\W{\alpha}|+s)}$.
		\end{itemize}
		Using these, we recover \eqref{bound:int}.
	 
	Now we are left with the corner  case $3\lceil q^\star/3 \rceil =3\left\lceil\frac{1}{3( |\W{\alpha}|+s)}\ln(n^s/\epsilon)\right\rceil > \frac{n}{4}$, i.e., $\epsilon$ is very small.
	In this case, we  show that the probabilistic bound trivially holds since it becomes larger than the trivial upper bound on $\norm{\M_\alpha}$:
	\begin{align*}
	&\left(6e \left\lceil \frac{1}{3(|\W{\alpha}| + s)}\ln(n^{s}/\epsilon)\right\rceil\right)^{\frac{1}{2}(|\W{\alpha}| + s )}\cdot (2|V(\alpha)|)^{|V(\alpha)|}  n^{\frac{1}{2}(|V(\alpha)|-s)}\\
			 	&> \left(\frac{en}{2}\right)^{\frac{1}{2}(|\W{\alpha}| + s )}\cdot (2|V(\alpha)|)^{|V(\alpha)|}  n^{\frac{1}{2}(|V(\alpha)|-s)}\\
			 	&= \left(\frac{e}{2}\right)^{\frac{1}{2}(|\W{\alpha}| + s )}\cdot (2|V(\alpha)|)^{|V(\alpha)|}  n^{\frac{1}{2}(|V(\alpha)|+|\W{\alpha}|)}\geq n^{\frac{1}{2}(|V(\alpha)|+|\W{\alpha}|)}\,.
			\end{align*}

\subsection{The case with nonempty left-right intersection: \texorpdfstring{$\U{\alpha}\cap \V{\alpha}\neq\emptyset$}{+1}}
\label{reduc:1}

Now we extend the argument from the previous section to the case where $\U{\alpha}\cap \V{\alpha}\neq \emptyset$.
The key observation is that in that case $M_\alpha$ is a block-diagonal matrix whose blocks \emph{behave} like the graph matrices corresponding to the shape without the intersection.

To be specific, consider a submatrix determined by fixing indices corresponding to the intersection $\U{\alpha}\cap \V{\alpha}$, i.e, fixing $\sigma(\U{\alpha}\cap \V{\alpha})$.
Then, one can readily notice from  the definition of $M_\alpha$ (Definition~\ref{def:graphmatrix}) that $M_\alpha$ is the matrix having such submatrices as diagonal blocks.
Note that the total number of blocks is 
$$n(n-1)\cdots (n+1-|\U{\alpha}\cap \V{\alpha}|) \leq n^{|\U{\alpha} \cap \V{\alpha}|}\,.$$
Moreover, one can observe that regardless of what we choose for the $\sigma(\U{\alpha}\cap \V{\alpha})$, the distribution of the resulting submatrix is the same.

Since the spectral norm of a block-diagonal matrix is equal to the maximum of those of individual blocks, one can bound the norm of each block and use the union bound to bound $\norm{\M_\alpha}$.
On the other hand, since each block has the same distribution, it is enough to prove the norm bound for an arbitrarily chosen block.
Let $\LL{\alpha}$ be a diagonal block corresponds to an arbitrarily fixed $\sigma(\U{\alpha}\cap \V{\alpha})$.
Then, in order to prove Theorem~\ref{thm:formal:mainresult} for this case, the union bound implies that it suffices to show:
\begin{align}
		 \pP\left[\norm{\LL{\alpha}} >   \left( 6e \left\lceil \frac{1}{3\ccal}\log(n^{s }/\epsilon)\right\rceil\right)^{\frac {1}{2}\ccal}\cdot  (2\dal)^{\dal}  n^{\frac{1}{2}(|V(\alpha)|-s)}\right] < \frac{\epsilon}{n^{|\U{\alpha}\cap \V{\alpha}|}} \,. \label{upper:inter}
\end{align}
Now we apply the same arguments as Section~\ref{sec:simple} to prove \eqref{upper:inter}.
Since the arguments overlap significantly, we do not repeat them here.
Below, we instead provide a high-level intuition as to why \eqref{upper:inter} is true.

The key principle is that $\U{\alpha}\cap \V{\alpha}$ can be \emph{effectively} ignored throughout the arguments in Section~\ref{sec:simple}.
To see this more precisely, first note that in the constraint graph, vertices in $ \U{\alpha}\cap \V{\alpha}$ have only a single copy of themselves.
Hence for any edge in $\alpha$ incident to a vertex $v\in \U{\alpha}\cap \V{\alpha}$, its copies in the constraint graph will automatically have the copy of $v$ as a common end.
Among the edges in $\alpha$ incident to a vertex in $\U{\alpha}\cap \V{\alpha}$, one can easily conclude that the only type of edges that might affect the value of $\E\left[\tr\left((M_{\alpha,P}M_{\alpha,P}^\top)^q\right)\right]$ are the edges between $\U{\alpha}\cap \V{\alpha}$ and $\W{\alpha}$.
However, even the copies of these edges  do not contribute to $\E\left[\tr\left((M_{\alpha,P}M_{\alpha,P}^\top)^q\right)\right]$ because for a constraint graph $C$ with $\val{C}\neq 0$, the copies of the endpoint in $W_{\alpha}$ must be paired up with each other. Thus, the product of the edge variables $\chi_e$ corresponding to the copies of these edges will be equal to $1$.

Consequently, the norm bound one can prove is equal to the result obtained in Section~\ref{sec:simple} applied to the shape $\alpha'$ which is obtained from $\alpha$ by removing the intersection $\U{\alpha}\cap \V{\alpha}$.
To write the result explicitly,  \eqref{bound:int} applied to $\alpha'$ yields
\begin{align*}  
    \pP\left[  \norm{\M_{\alpha'}} \geq  \left(6e \left\lceil \frac{1}{3(|\W{\alpha'}| + s')}\ln(n^{s'}/\epsilon')\right\rceil\right)^{\frac{1}{2}(|\W{\alpha'}| + s' )}\cdot (2|V(\alpha')|)^{|V(\alpha')|}  n^{\frac{1}{2}(|V(\alpha')|-s')}\right] \leq \epsilon'\,,
\end{align*}
where $s'$ is the  size of a minimum vertex separator between  $\U{\alpha'}$ and $\V{\alpha'}$.
One can easily verify:
\begin{itemize}
    \item $s' = s-|\U{\alpha}\cap \V{\alpha}|$ ;
    \item $|\W{\alpha'}|+s' = |\W{\alpha}|+s-|\U{\alpha}\cap \V{\alpha}| = \ccal$; 
    \item $|V(\alpha')| = |V(\alpha)|- |\U{\alpha}\cap \V{\alpha}| = \dal$; and 
    \item $|V(\alpha')|-s' = |V(\alpha)|-s$.
\end{itemize}
Consequently, the above bound becomes
\begin{align*} 
    \pP\left[  \norm{\M_{\alpha'}} \geq   \left(6e \left\lceil \frac{1}{3\ccal }\ln(n^{s-|\U{\alpha}\cap \V{\alpha}|}/\epsilon')\right\rceil\right)^{\frac{1}{2}\ccal}\cdot (2\dal)^{\dal}  n^{\frac{1}{2}(|V(\alpha)|-s)} \right] \leq \epsilon'\,,
\end{align*}
which recovers \eqref{upper:inter} by making the substitution $\epsilon' \leftarrow \frac{\epsilon}{n^{|\U{\alpha}\cap \V{\alpha}|}}$.

\subsection{The case with isolated middle vertices: \texorpdfstring{$\iso \neq \emptyset$}{+1}}

\label{reduc:2}
Lastly, we extend our proof to the case $\W{\alpha}\neq \emptyset$.
In particular, we reduce the case with isolated vertices to the case without them via the following proposition:
		\begin{proposition}  \label{prop:isolated}
			Given a shape  $\alpha = (\U{\alpha},\V{\alpha},\W{\alpha},E(\alpha))$, let $\iso \subset \W{\alpha}$ be the set of isolated vertices in $\W{\alpha}$, and let $\tnoniso:= |V(\alpha)|-|\iso|$.
			Let $\beta$ be the shape obtained from $\alpha$ by removing the isolated middle vertices, i.e. $\beta = (\U{\alpha}, \V{\alpha}, \W{\alpha}\setminus \iso, E(\alpha))$.
			Then, $$\M_\alpha = (n-\tnoniso)(n-\tnoniso-1)\cdots (n-|V(\alpha)|+1) \cdot \M_{\beta}\,.$$
			Consequently, we have $\norm{\M_\alpha} \leq n^{|\iso|}\norm{\M_{\beta}} $.
		\end{proposition}	
		\begin{proof}
			Recall from the definition of graph matrices (Definition~\ref{def:graphmatrix}) that 
			\begin{align*}
			M_{\alpha}(A,B) = \sum_{\substack{\sigma: V(\alpha) \to [n]:\\\sigma \text{ is injective},\\ \sigma(\U{\alpha}) = A,~ \sigma(\V{\alpha}) = B}}{\chi_{\sigma(E(\alpha))}}\,.
			\end{align*}
			Now observe that $\chi_{\sigma(E(\alpha))}$ is equal to $\chi_{\sigma'(E(\alpha))}$ as long as $\sigma(x)=\sigma'(x)$ for $x\in V(\alpha)\setminus \iso$.
			Thus, we can collect identical Fourier characters in the above summation. After fixing the values of $\sigma$ restricted to the subset $ V(\alpha)\setminus \iso$, there are $(n-\tnoniso)(n-\tnoniso-1)\cdots (n-|V(\alpha)|+1)$ different possible assignments for $\sigma(\iso)$.
			Since each such assignment results in the same Fourier character, we have that
			\begin{align*}
			M_{\alpha}(A,B) &= (n-\tnoniso)(n-\tnoniso-1)\cdots (n-|V(\alpha)|+1)\cdot \sum_{\substack{\sigma: V(\alpha)\setminus \iso \to [n]:\\\sigma \text{ is injective},\\ \sigma(\U{\alpha}) = A,~ \sigma(\V{\alpha}) = B}}{\chi_{\sigma(E(\alpha))}} \\
			&= (n-\tnoniso)(n-\tnoniso-1)\cdots (n-|V(\alpha)|+1) \cdot M_{\beta}(A,B) \,.
			\end{align*}
		\end{proof}
Using Proposition~\ref{prop:isolated}, it is sufficient to prove Theorem~\ref{thm:formal:mainresult} for the case where there are no isolated vertices, which we did in the previous subsection.

	\section{Generalized graph matrices}\label{gendef}

	So far, we have considered a restrictive setting where the index monomials are multilinear and the input distribution is a random graph on the ground set $[n]$.
	In this section, we will relax these restrictions and come up with the notion of generalized graph matrices.
	We will then extend the techniques developed in Section~\ref{sec:modified} to this generalized setting.

    \subsection{Motivation for generalizing graph matrices}
    We begin by giving several examples to illustrate why generalizing the concept of graph matrices is useful.
    \begin{example}
    Let's say that instead of a random graph, the input is a $m \times n$ matrix $X$ with random $\pm{1}$ entries. In this case, we represent the entry $X_{ij}$ of $X$ by an edge between a vertex $i$ and a vertex $j$. However, here it is important to know which index is the row index and which index is the column index. To handle this, we consider two different types of vertices, row vertices and column vertices. Row vertices take values in $[m]$ while column vertices take values in $[n]$.
    \end{example}
    \begin{example}\label{hyperedgeexample}
    Let's say that the input is an $n \times n \times n$ tensor $T$ such that
    \begin{enumerate}
        \item $T$ is symmetric (i.e. $\forall i,j,k \in [n], T_{ijk} = T_{ikj} = T_{jik} = T_{jki} = T_{kij} = T_{kji}$).
        \item Any entry of $T$ with a repeated index is $0$ (i.e $\forall i \in [n], T_{iii} = 0$ and $\forall i,j \in [n], T_{iij} = 0$).
        \item All entries $T_{ijk}$ where $i < j < k$ are independently chosen to be $-1$ or $1$.
    \end{enumerate}
    In this case, we need $3$ indices to describe each entry of the input rather than $2$. To handle this, we represent each entry $T_{ijk}$ where $i,j,k$ are distinct by a hyperedge between vertices $i$, $j$ and $k$.
    \end{example}
    \begin{example}\label{gaussianinputexample}
    Let's say that the input is a symmetric matrix $X$ with random Gaussian off-diagonal entries and zeros on the diagonal. In this case, we need to represent higher degree polynomials of the input entries like $X_{ij}^2$. To handle this, we take an orthonormal basis $h_1,h_2,\ldots$ for our input distribution, decompose each polynomial in terms of this basis, and represent $h_{k}(X_{ij})$ by an edge between vertices $i$ and $j$ with label $k$.
    
    For Gaussian inputs, we take the orthonormal basis to be the normalized Hermite polynomials $h_1(x) = x$, $h_2(x) = \frac{x^2 - 1}{\sqrt{2!}}$, $h_3(x) = \frac{x^3 - 3x}{\sqrt{3!}}$, etc. To represent $X_{ij}^2$, we write $X_{ij}^2 = \sqrt{2}h_2(X_{ij}) + 1$, which then can be represented as $\sqrt{2}$ times an edge with label $2$ between vertices $i$ and $j$ plus a non-edge (which represents $1$).
    \end{example}
    \begin{remark}
    In Examples \ref{hyperedgeexample} and \ref{gaussianinputexample}, we restricted which entries of the matrix or tensor were nonzero so that every input entry was represented by the same kind of edge/hyperedge. If we had nonzero $T_{iii}$, $T_{iij}$, and/or $X_{ii}$ entries, we would need to represent these entries with different kinds of hyperedges.
    \end{remark}
    \begin{example}
    When analyzing the Sum-of-Squares hierarchy, we need to analyze a moment matrix $M$ which is indexed by monomials $p$ and $q$. When the variables of $p$ and $q$ are Boolean variables or $\pm{1}$ variables, we can assume that $p$ and $q$ are multilinear. However, we cannot make this assumption in general. For example, we may need to consider the entry $M_{pq}$ where $p = x^2_1{x_2}$ and $q = x^2_2{x_3}$. 
    
    To handle this, we define non-multilinear matrix indices where we specify the degree of each element of the matrix index.
    \end{example}
    Based on these examples, we introduce the following generalizations for graph matrices.
    \begin{enumerate}
        \item Different types of vertices to handle different types of indices.
        \item Hyperedges to handle input entries which are described by more than $2$ indices.
        \item Labeled edges/hyperedges to handle different input distributions.
        \item Non-multilinear matrix indices.
    \end{enumerate}
	\subsection{Generalized definitions} \label{def:gen}
\subsubsection{Generalized matrix indices and index shapes}
	We first generalize the index monomials.
	We consider the index monomials that (i) are not necessarily multilinear and (ii) consist of variables of multiple types.
	For instance, we can now have $x_1^2 x_2x_3^4y_1y_2^5z_1^3$ as a index monomial (here there are $3$ different types of variables, $x,y,z$).
	We assume throughout this section that there are $t$ different types of variables, and for each $i\in [t]$, let us denote by $x_{i;1},x_{i;2},x_{i;3},\cdots$ the variables of type $i$.
	Now we generalize the definition of ground set (Definition~\ref{simp:ground}) as follows:	\begin{definition}[Ground sets]
		For each $i=1,2,\dots, t$, we denote by $\N{i}$  the $i$-th ground set for the indices of  variables of type $i$.
		In other words, we consider monomials consisting of variables $\bigcup_{i=1}^t\{x_{i;j}\}_{j\in \N{i}}$.
		Unless stated otherwise, we consider $\N{i}=[n_i]$ for $n_i\in\na$, $i=1,2,\dots, t$.
	\end{definition}

		In order to represent the generalized index monomials, we also need to generalize Definition~\ref{simp:indices}.
	More specifically, since an index monomial is  no longer multilinear, it \emph{cannot} be represented by a single matrix index.
	To cope with this, we decompose a given monomial into pieces so that each piece consists of variables of the same type to the same power.
	To illustrate, let us consider the following example:
	\begin{example}[Decomposing monomials according to types and powers] Consider the case $t=3$, i.e., there are three types of variables. 
	Let us denote the three different variables by $x,y,z$ for simplicity.
	Moreover assume that $N_1=[5], N_2=[3], N_3=[2]$ and consider a monomial $x_1^4z_1^8y_1x_2^4y_2x_3^3y_3z_2^4x_4^4x_5^3$
	Then, one can decompose the monomial according to types and powers as follows:
	\begin{enumerate}
	    \item Decompose them according to  types: $x_1^4x_2^4x_3^3x_4^4x_5^3$,  
	   $y_1y_2y_3$, $z_1^8z_2^4$.
	    \item Further decompose them according to powers: $x_3^3x_5^3$, $x_1^4x_2^4x_4^4$, $y_1y_2y_3$, $z_2^4$, $z_1^8$.
	\end{enumerate}
	\end{example}
\noindent Having decomposed an index monomial into pieces, one can represent each piece with a tuple.
Formally, we consider the following definition about generalized matrix index pieces; again due to a technical reason addressed in Remark~\ref{rmk:why} we adopt tuples over subsets for representations.
	\begin{definition}[Generalized matrix index pieces] \label{def:gmp}
	For a given collection of ground sets 
	$\{\N{i}\}$, we define a generalized matrix index piece $A = ((a_{1},\ldots,a_{m}), i, p)$ for some $m\in \na$, where 
	\begin{enumerate} 
			\item[(i)]  $i\in[t]$ denotes the type of variables,
			\item[(ii)]  $p\in \mathbb{N}$ denotes the power of variables, and 
			\item[(iii)] $a_1,\dots,a_{m}\in \N{i}$ are distinct  elements representing variable indices.
		\end{enumerate}
	In other words, the generalized matrix index piece $A$ corresponds to the monomial $p_{A} := \prod_{j = 1}^{m}{x^{p}_{i;a_j}}$.
		In addition, we make the following definitions about generalized matrix index pieces:
		\begin{enumerate}
		    \item We denote by $|A|$ the size of the generalized matrix index piece, i.e., $|A|=m$.
			\item We define $V(A)$ to be the set of vertices $\{(a_{j},i): j=1,2,\dots, m\}$.
			\item We say that $A_1$ and $A_2$ are disjoint if $V(A_1) \cap V(A_2) = \emptyset$.
		\end{enumerate}
	\end{definition}
	\noindent Let us illustrate generalized matrix index pieces through an example:
	\begin{example}[Examples of generalized matrix index pieces]
Let $t=3$, i.e., there are $3$ different types of variables, and let $\N{1}=[4]$, $\N{2}=[5]$, and $\N{3}=[3]$.
To simplify notations, let us denote the three types of variables by $x$, $y$, and $z$.
For example, $A_1 = ((1,2,4),1,1)$ is the generalized matrix index piece corresponding to the multi-type monomial $x_1x_2x_4$.
In this example, $V(A_1) =\{(1,1),(2,1),(4,1) \}$ and $|A|=3$.
Next, $A_2 = ((2,5,3,4) ,2, 3)$ corresponds to $y_2^3y_3^3y_4^3y_5^3$.
In this example, $V(A_2) =\{(2,2),(3,2),(4,2), (5,2) \}$ and $|A|=4$.
Lastly, $A_3 = ((2,1),3,2)$ corresponds to $z_1^2z_2^2$.
In this example, $V(A_2) =\{(1,3),(2,3)\}$ and $|A|=2$.
	
	\end{example}
	
	\noindent With generalized matrix index pieces, we can generalize Definition~\ref{simp:indices} as follows:
	
	\begin{definition}[Generalized matrix indices]
		For a given collection of ground sets $\{\N{i}\}$, we define a generalized matrix index $A = \{A_i=\uu{a_{i;1},\dots,a_{i;|A_i|} }{t_i}{p_i}\}$ to be a set of disjoint generalized matrix index pieces where for all $i < j$, either (i) $t_i < t_j$ or (ii) $t_i = t_j$ and $p_i < p_j$.
		We make the following definitions about generalized matrix indices:
		\begin{enumerate}
			\item We associate the monomial $p_{A} = \prod_{A_i \in A}{p_{A_i}}$ with $A$.
			\item We define the size $|A|$ to be a $t$-tuple $(s_1,s_2,\dots, s_t)$, where $s_i$ is the sum of the sizes of the generalized matrix index pieces whose types are $i$.
			\item We define $V(A)$ to be the set of vertices $\cup_{A_i \in A}{V(A_i)}$.
		\end{enumerate}
	\end{definition}

		Next, we consider generalized index shapes. 
Paralleling Definition~\ref{def:gmp}, we first define generalized index shape pieces:
	\begin{definition}[Generalized  index shape pieces]  \label{def:gisp}
	For a given collection of ground sets  $\{\N{i}\}$, we define a generalized index shape piece to be a tuple of distinct variables $U = ((u_{1},\ldots,u_{m}), i, p)$ for some $m\in \na$, where 
	\begin{enumerate} 
			\item[(i)]  $i\in[t]$ specifies the type of variables,
			\item[(ii)]  $p\in \mathbb{N}$ specifies the power of variables, and 
			\item[(iii)] $u_1,\dots,u_{m}\in \N{i}$ are distinct variables.
		\end{enumerate}
 	In addition, we make the following definitions about generalized index shape pieces:
		\begin{enumerate}
		    \item We denote by $|U|$ the size of the generalized index shape piece, i.e., $|U|=m$.
			\item We define $V(U)$ to be the set of vertices $\{(u_{j},i): j=1,2,\dots, m\}$.
			\item We say that two generalized index shape pieces $U_1$ and $U_2$ are disjoint if $V(U_1) \cap V(U_2) = \emptyset$.
			\item If $U$ and $V$ are two generalized index shape pieces, we say that $U = V$ if $U$ and $V$ have the same type and power and $|U| = |V|$.
		\end{enumerate}
	\end{definition}
 
	\noindent Based on Definition~\ref{def:gisp}, we can generalize index shapes as follows:
	
	\begin{definition}[Generalized index shapes] We define a generalized index shape $U=\{U_i=\uu{u_{i;1},\dots,u_{i;|U_i|} }{t_i}{p_i}\}$ to be a set of generalized index shape pieces such that for all $i < j$, either (i) $t_i < t_j$ or (ii) $t_i = t_j$ and $p_i < p_j$. We make the following definitions about generalized index shapes:
		\begin{enumerate}
			\item We define the size $|U|$ to be a $t$-tuple $(s_1,s_2,\dots, s_t)$, where $s_i$ is the sum of the sizes of the generalized index shape pieces whose types are $i$.
			\item 	We define $V(U)$ to be the set of vertices $\cup_{U_i\in U}V(U_i)$.
			\item 	If $U = \{U_i\}$ and $V = \{V_i\}$ are two generalized index shapes, we say that $U = V$ if $U$ and $V$ have the same number of pieces and for all $i$, $U_i = V_i$.
		\end{enumerate}
	\end{definition}
	
	Finally, we make the following definition about the relation between generalized matrix indices and generalized index shapes.
	\begin{definition}
		We say a generalized matrix index $A = \{A_i\}$ has generalized index shape $U = \{U_i\}$ if $A$ and $U$ have the same number of pieces and for all $i$, $|A_i| = |U_i|$ and $A_i$ and $U_i$ have the same type and power. In other words, there is an assignment of values to the unspecified indices of $U$ which results in $A$. 
	\end{definition}
	
\subsubsection{Generalized input distributions}
	Next, we generalize the input distribution.
	Since the matrix indices have multiple types, the generalized input distribution is now over vertices of multiple types.
	Moreover, we generalize random graphs to random \emph{hypergraphs}, i.e., some subsets of nodes (not necessarily of size $2$) are assigned independent random variables.
	The distribution of these random variables is denoted by $\Omega$, which is not necessarily the $\{\pm 1\}$ Bernoulli distribution. In order to do Fourier analysis on the input, we need to use an orthonormal basis for $\Omega$.
	\begin{definition}[Orthonormal basis for $\Omega$]
		Given an input distribution $\Omega$, we define the polynomials $\{h_i\}$ to be the ones found through the Gram-Schmidt process such that 
		\begin{enumerate}
			\item $\forall i, E_{\Omega}[h^2_i(x)] = 1$
			\item $\forall i \neq j, E_{\Omega}[h_i(x)h_j(x)] = 0$
			\item For all $i$, the leading coefficient of $h_i(x)$ is positive.
		\end{enumerate}
		If $\Omega$ has finite support, then we will only have the polynomials $\{h_i: i \in [k-1]\}$ where $k$ is the size of the support of $\Omega$.
	\end{definition}
	\begin{example}[Fourier analysis for the normal distribution] \label{ex:fourier}
		If $\Omega$ is the normal distribution, then the polynomials $\{h_i\}$ are the Hermite polynomials with the appropriate normalization so that for all $i$, $E_{\Omega}[h^2_i(x)] = 1$. For example, $h_0(x) = 1$, $h_1(x) = x$, $h_2(x) = \frac{x^2 - 1}{\sqrt{2!}}$, $h_3(x) = \frac{x^3 - 3x}{\sqrt{3!}}$, and $h_4(x) = \frac{x^4 - 6x^2+3}{\sqrt{4!}}$.
	\end{example}
\noindent	Given this orthonormal basis for $\Omega$, we assign a label to each hyperedge:
\begin{definition}[Labels of hyperedges]
We assign a label $l_e$ to each hyperedge $e$ to denote which element $h_{l_e}$ of the orthonormal basis to apply to the random variable associated to $e$. 
More formally, if the hyperedge $e$ is assigned a label $l_e$ then this hyperedge corresponds to $h_{l_e}(x_e)$ where $x_e$ is the random variable associated with $e$.
\end{definition}
	With these edge labels, we consider the following Fourier characters. 
	\begin{definition}[Fourier characters]
		Given a multi-set of labeled hyperedges $E$, we define the Fourier character $\chi_E$ to be 
		\[
		\chi_E = \prod_{e:~e\in E}{h_{l_e}{(x_e)}}\,.
		\]
	\end{definition}
	
	\subsubsection{Generalized ribbons and shapes}
	
	Having established the above definitions, we are ready to consider generalized ribbons:

	\begin{mdframed}[style=box]\begin{definition}[Generalized ribbons] \label{def:gribbon}
			A generalized ribbon $R = (A_R,B_R,C_R,E(R))$  consists of two generalized matrix indices $A_R$ and $B_R$ (possibly $V(A_R)\cap V(B_R)\neq\emptyset$), an additional set of indices with types $C_R$ such that $C_R \cap (V(A_R) \cup V(B_R)) = \emptyset$, and a set $E(R)$ consisting of labeled subsets from $V(A_R) \cup V(B_R)\cup  C_R$. We make the following definitions about generalized ribbons:
		\begin{enumerate}
			\item (Graphical representation of a generalized ribbon) We identify $R$ with its graphical representation defined as a labeled hypergraph with vertices $V(R) = V(A_R) \cup V(B_R) \cup C_R$ and labelled hyperedges $E(R)$.
			Here we regard $V(A_R)$ and $V(B_R)$ as distinguished sets, and refer to them as the left and right vertices of $R$, respectively. We refer to $C_R$ as the middle vertices of $R$.
				
			\item (Fourier character of a generalized ribbon) We define the Fourier character $\chi_{R}$ to be $\chi_{E(R)}$.
		
			\item (Matrix associated to a generalized ribbon) Given a generalized ribbon $R$ where $|A_R|=(s_1,s_2,\dots , s_t )$ and $|B_R|=(r_1,r_2,\dots, r_t)$ and a random input $X$ with a random variable $x_e$ for each hyperedge $e \in E(R)$, we define   $M_R$ to be a $\prod_{i=1}^t\frac{n_i!}{(n_i-s_i)!}\times \prod_{i=1}^t\frac{n_i!}{(n_i-r_i)!}$ matrix such that $M_R(A_R,B_R) = \chi_{R}(X)$ and $M_R(A,B) = 0$ if $A \neq A_R$ or $B \neq B_R$.
		\end{enumerate}
	\end{definition}
	\end{mdframed}

	\noindent	We also define generalized shapes as follows:
			\begin{mdframed}[style=box]
	\begin{definition}[Generalized shapes] \label{def:gshape}
		A generalized shape $\alpha = (\U{\alpha},\V{\alpha},\W{\alpha},E(\alpha))$ consists of generalized index shapes $\U{\alpha}$ and $\V{\alpha}$ (possibly having common variables), an additional set $\W{\alpha}$  of variables with types distinct from  the ones in $V(\U{\alpha})\cup V(\V{\alpha})$, and  $E(\alpha)$ consisting of labeled hyperedges from $V(\U{\alpha}) \cup V(\V{\alpha})\cup \W{\alpha}$.
		
			Furthermore, we identify $\alpha$ with its graphical representation defined as a labeled hypergraph graph with vertices $V(\alpha)=V(\U{\alpha}) \cup V(\V{\alpha}) \cup \W{\alpha}$ and labeled hyperedges $E(\alpha)$.
		Here we regard $V(\U{\alpha})$ and $V(\V{\alpha})$ as distinguished sets, and refer to them as the left and right vertices of $\alpha$, respectively. We refer to $\W{\alpha}$ as the middle vertices of $\alpha$.
	\end{definition}
	\end{mdframed}
 There are a few more ingredients required for generalized graph matrices.
	\begin{definition}[Type-respecting realizations]
		Given a realization $\sigma:U\to \cup_{i}{N_i}$,  it is said to be type-respecting if for each $i$, vertices of type $i$ are mapped to the $i$-th ground set $\N{i}$.
	\end{definition}
	
	\begin{definition}[Realizations of generalized shapes]
		Given a generalized shape $\alpha$ and a type-respecting realization $\sigma: V(\alpha) \to \cup_{i}{N_i}$, we define   the ribbon $\sigma(\alpha):=(\sigma(\U{\alpha}), \sigma(\V{\alpha}), \sigma(\W{\alpha}), \sigma(E(\alpha)))$.
	\end{definition}
	\begin{definition}[Shapes of generalized ribbons]
	Given a generalized shape $\alpha$ and a generalized ribbon $R$, we say $R$ has shape $\alpha$ if there exists a type-respecting realization $\sigma$ such that $\sigma(\alpha)=R$.
	\end{definition}
	
\noindent With all the definitions above, we are now ready to define generalized graph matrices.

\begin{mdframed}[style=box]
	\begin{definition}[Generalized graph matrices] \label{def:gengraphmatrix}
		Given a generalized shape $\alpha$, the generalized  graph matrix $M_{\alpha}$ is defined as follows:
		for generalized matrix indices $A$ and $B$ such that $|A|=|\U{\alpha}|$ and $|B|= |\V{\alpha}|$,
		\begin{align} \label{def:g1}
		M_{\alpha}(A,B) = \sum_{\substack{\sigma \text{ is a type-respecting realization of } \alpha,\\~\sigma(\U{\alpha}) = A, \sigma(\V{\alpha})=B}}{\chi_{\sigma(E(\alpha))}}\,.
		\end{align}
	\end{definition}
\end{mdframed}	
	\begin{remark}
	Paralleling Remark~\ref{rmk:altdef1},  an alternative way to define the generalized graph matrix is to consider
		\begin{align} \label{def:g2}
		M_{\alpha} = \sum_{R \text{ is a generalized ribbon of shape } \alpha}{M_{R}}\,.
		\end{align}
		Again, similarly to Remark~\ref{rmk:altdef1}, the two definitions \eqref{def:g1} and \eqref{def:g2} only differ by a multiplicative constant.
	\end{remark}
	
	\subsection{Generalized techniques}
	In this section, we generalize  the techniques in Section~\ref{sec:modified} so that they can handle this generalized setting. 
	In particular, we will generalize the vertex partitioning technique as well as constraint graphs.
	
	\subsubsection{Generalized vertex partitioning} 
	Paralleling the vertex partition argument in Section~\ref{sec:modified}, we consider an $\alpha$-partition of the ground set.  
	Here the main difference is that 
	for each $i\in[t]$, we partition the $i$-th ground set according to the vertices of type $i$ in $\alpha$. 
	Formally, we make the following definitions about $\alpha$-partitions.
	\begin{definition}[$\alpha$-partitions]
		Given a generalized shape $\alpha$ and for each $i\in[t]$,  let $V_i(\alpha)$ be the set of vertices of type $i$, and  let $P_i$ be a tuple  $(I_v: v \in V_i(\alpha))$ such that $\cup_{v \in V_i(\alpha)}{I_{v}}$ is a partition of $\N{i}$.
		We define an $\alpha$-partition of the ground sets to be a tuple of partitions $(P_i~:~i\in[\type])$.
		We define $\mathcal{P}_{\alpha}$ to be the set of all possible $\alpha$-partitions of the ground set.
	\end{definition}
\noindent	Similarly to Proposition~\ref{cardinality}, we can count the number of possible $\alpha$-partitions:
	\begin{proposition} \label{g:cardinality}
		$|\mathcal{P}_\alpha|=\prod_{i=1}^t |V_i(\alpha)|^{|\N{i}|}$.
	\end{proposition}
\noindent Given this definition of $\alpha$-partitions, the notion of $P$-respecting map can be analogously defined:
	\begin{definition}
		Given an $\alpha$-partition $P$ of the indices, we say that a type-respecting and injective map $\phi:V(\alpha) \to \N{}$ is $P$-respecting if $\forall v \in V(\alpha)$, $\phi(v) \in I_v$.
	\end{definition}
\noindent	The following is a direct generalization of Proposition~\ref{respectcount}:
	\begin{proposition} \label{g:respectcount}
		For any type-respecting and injective map $\phi:V(\alpha) \to \N{}$, there are exactly $\prod_{i=1}^{t}|V_i(\alpha)|^{|\N{i}| - |V_i(\alpha)|}$ $\alpha$-partitions $P \in \mathcal{P}_{\alpha}$ such that $\phi$ is $P$-respecting.
	\end{proposition}
\noindent	Having established all of this, we can similarly define $M_{\alpha,P}$ as follows:
	\begin{definition}
		Given a random input $X$, We define $M_{\alpha,P}$ to be the matrix with entries \[
		M_{\alpha,P}(A,B) := \sum_{\substack{\phi:V(\alpha) \to \N{} : \phi \text{ is } P \text{-respecting}, \\ \phi(\U{\alpha}) = A, \phi(\V{\alpha}) = B}}{\chi_{\phi(E(\alpha))}(X)}\,.
		\]
	\end{definition}
\noindent	The following result is an immediate consequence of Propositions~\ref{g:cardinality} and~\ref{g:respectcount}:
	\begin{corollary}[Decomposition induced by an $\alpha$-partition] \label{cor:vpa}
		$M_{\alpha} = \frac{\prod_{i=1}^t{|V_i(\alpha)|^{|V_i(\alpha)|}}}{|\mathcal{P}_\alpha|}\sum_{P \in \mathcal{P}_\alpha}{M_{\alpha,P}}$.
	\end{corollary}

	\subsubsection{Generalized constraint graphs}
	Having defined $\alpha$-partitions and $M_{\alpha,P}$, the constraint graphs can be defined akin to Section~\ref{constraintgraph}.
	Actually, the definitions and discussions in Section~\ref{constraintgraph} go through in the exactly same manner except that $\phi$ is now  a type-respecting map from $H(\alpha,2q)$ to $\cup_{i}{N_i}$.
	Moreover, we can define the notion of $P$-respecting maps and well-behaved constraint graphs analogously.
	Therefore, we obtain the following generalization of Corollary~\ref{expan2}:
	\begin{proposition}\label{expan3}
		For all generalized shapes $\alpha$ and all $P \in \mathcal{P}_{\alpha}$, and all $q \in \mathbb{N}$,
		\[
		\E\left[\tr\left((M_{\alpha,P}M_{\alpha,P}^\top)^q\right)\right] = \sum_{C \in \mathcal{C}_{(\alpha,2q)}: C \text{ is well-behaved}}{N_{P}(C)\val{C}}\,,
		\]
		where  $N_P(C) := |\{\phi: V(H(\alpha,2q)) \to  \N{}: \phi \text{ is } P \text{-respecting}, ~C(\phi) \equiv C\}|$.
	\end{proposition}
	
As we will soon see in Section~\ref{pf:generalgeneral}, the proof of norm bounds for generalized graph matrices will follow by bounding the quantities $N_P(C)$ and $\val{C}$. 
	    However, since the constraint graphs now admit different types of vertices, counting $N_P(C)$ will require some more delicate combinatorial arguments.
		Moreover, since the input distribution no longer satisfies $\Omega=\{\pm 1\}$, we cannot simply assert that  $\val{C}\leq 1$ as before.
		Bounding $\val{C}$ under some mild conditions on $\Omega$ will be the focus of the next subsection.
	
	\subsection{Handling general input distributions}
	\label{sec:handle}
To come up with a bound on $\val{C}$, we first make the following definition about bound functions for a given input distribution.
	
	\begin{definition}[Bound function of input distribution]
		For any input distribution $\Omega$, we say $B_{\Omega}: \mathbb{N} \to \mathbb{R}$ is a bound function of $\Omega$ if it satisfies:
		\begin{enumerate}
			\item It is non-decreasing.
			\item $\forall j \in \mathbb{N}, \left|\E_{\Omega}[x^j]\right| \leq B_{\Omega}(j)^{j}$.	
	\end{enumerate}  	\end{definition}
\noindent Note that there could be multiple bound functions for a single input distribution. 
		For a better bound on $\val{C}$, it is suggested to take one that takes minimal values.
		This will be clearer when we state our norm bounds for generalized graph matrices.

Moreover, to upper bound each polynomial in the orthonormal basis $\{h_i\}$ of $\Omega$, we consider the following polynomials:
	\begin{definition}
		Given an input distribution $\Omega$, if $h_k(x) = \sum_{j=0}^{k}{c_{kj}x^{j}}$ is the degree $k$ polynomial in the orthonormal basis for $\Omega$ then we define $h^{+}_k(x)$ to be $h^{+}_{k}(x) = \sum_{j=0}^{k}{|c_{kj}|x^{j}}$.
	\end{definition}
\noindent With these definitions of bound functions and ${h^+_k(x)}$, we obtain the following upper bounds:
	\begin{lemma} \label{upper:dist}
		For all $k,r \in \mathbb{N}$, $\E_{\Omega}[h_k(x)^{r}] \leq h^{+}_{k}(B_{\Omega}(kr))^{r}$.
	\end{lemma}
	\begin{proof}
		Writing $h_k(x)^{r} = \sum_{j=0}^{kr}{{a_j}x^j}$ and $h^{+}_k(x)^{r} = \sum_{j=0}^{kr}{{b_j}x^j}$, we can easily deduce that  $|a_j| \leq b_j$ for all $0\leq j \leq kr$.
		Thus, we have 
		\begin{align*}
		\E_{\Omega}[h_k(x)^{r}] &= \sum_{j=0}^{kr}{{a_j}\E_{\Omega}[x^j]}\leq \sum_{j=0}^{kr}{{b_j}\left|\E_{\Omega}[x^j]\right|} \leq \sum_{j=0}^{kr}{{b_j}B_{\Omega}(kr)^{j}} = h^{+}_k(B_{\Omega}(kr))^{r}\,,
		\end{align*}
	where the second inequality follows from the definition of the bound function: for $0\leq j \leq kr$,  $\left|\E_{\Omega}[x^j]\right| \leq B_{\Omega}(j)^j\leq B_{\Omega}(kr)^j$.
	\end{proof}
\noindent Thanks to the above upper bound, we obtain the following upper bound on $\val{C}$:
	\begin{corollary}\label{cor:valCbound}
		For a generalized shape $\alpha$ with labels $\{l_e\}_{e\in E(\alpha)}$ of the hyperedges, let $C \in \mathcal{C}_{(\alpha,2q)}$ be a well-behaved constraint graph. Then, the following holds: 
		\[
		\val{C} \leq \prod_{e \in E(\alpha)}{h^{+}_{l_{e}}(B_{\Omega}(2q{l_e}))^{2q}}.
		\]
	\end{corollary}
	\begin{proof}
		Recall that $\val{C} = \E_X[\chi_{\phi(E(\alpha,2q))}(X)]$.
			The bound is trivial when $\val{C}=0$, so consider the case where $\val{C}\neq 0$.
			Then, each hyperedge in the multi-set $\phi(E(\alpha,2q))$ appears an even number of times.
			Moreover, since $C$ is well-behaved, two hyperedges in $E(\alpha,2q)$ correspond to the same hyperedge under the map $\phi$ only if they are copies to each other.
			Hence, due to independence between hyperedge variables, the expectation value becomes $\prod_{e\in E(\alpha)} \ex_X[\chi_{\phi(C_e)}(X)]$, where $C_e\subset H(\alpha,2q)$ is the subset consisting of $2q$ copies of $e$.

			Fix a hyperedge $e\in E(\alpha)$. 
			Let $r_1,\dots, r_p$ be the multiplicities of hyperedges appearing in the multi-set $\phi(C_e)$, i.e., $2q$ elements of $\phi(C_e)$ can be partitioned into $p$ subsets so that each subset consists of the same hyperedge, and different subsets contain different elements.
			Recall that since $\val{C}\neq 0$, each $r_i$ is even for $i=1,2,\dots, p$. 
			Then, we have 
			$$ \ex_X[\chi_{\phi(C_e)}(X)] = \prod_{i=1}^p \ex_X[h_{l_e}(x)^r_i]  \leq \prod_{i=1}^ph^+_{l_e}(B_\Omega(l_e r_i))^{r_i}\leq \prod_{i=1}^p h^+_{l_e}(B_\Omega(2ql_e ))^{r_i} =h^+_{l_e}(B_\Omega(2ql_e ))^{2q} \,,$$
			where the first inequality is due to Lemma~\ref{upper:dist}, the second inequality is due to the facts that $B_\Omega(\cdot)$ is increasing and $r_i\leq 2q$ $\forall i=1,\dots, p$, and the last equality is due to $\sum_{i=1}^pr_i = 2q$.
		
	Taking a product over all $e\in E(\alpha)$, we get 
			$$\val{C} = \prod_{e\in E(\alpha)} \ex_X[\chi_{\phi(C_e)}(X)] \leq \prod_{e\in E(\alpha)} h^+_{l_e}(B_\Omega(2ql_e ))^{2q}\,,$$
			as desired. 
	\end{proof}

	\section{Norm bounds on generalized graph matrices}\label{pf:generalgeneral}

In this section, we state and  prove norm bounds for generalized graph matrices.
		Recall that in the case of single-type vertex graph matrices (Section~\ref{simpledef}), the norm bound depends on the minimal size of a vertex separator of $U$ and $V$ of $\alpha$.
		In the general case of multiple types, a similar story goes through:
		viewing the shape $\alpha$ as some type of graph with weighted vertices, it turns out that the combinatorial quantity that captures the norm bound is the weight of the \emph{minimum weight seperator} between $U$ and $V$.
		To formalize this, we make the following definitions about weights and minimum weight separators of shapes:
	\begin{definition}[Weights in shapes and constraint graphs]
		Let $n=\max_{i\in[t] }|\N{i}|$, $\alpha = (\U{\alpha},\V{\alpha},\W{\alpha},E(\alpha))$ be a generalized shape, and   $C \in \mathcal{C}_{(\alpha,2q)}$ be a well-behaved constraint graph.
		We make the following definitions about weights of vertices in $\alpha$ or $C$:
		\begin{enumerate}
			\item For each vertex $v$ in  $\alpha$ (or in $C$), we define the weight $\w{v}$ to be $\log_n|\N{i}|$, where $i$ is the type of the vertex $v$ (or the type of the original vertex).
			\item For a subset $S$ of $V(\alpha)$ (or $V(C)$), we define $\w{S}$ to be $\sum_{v\in S}\w{v}$. 
			\item For each constraint edge $e$ in $E(\alpha)$ (or in $E(C)$), we define the weight $\w{e}$ to be the weight of either end of $e$. (Note that this is well-defined since $C$ is well-behaved.)
			\item For a subset of constraint edges $D$ of $V(\alpha)$ (or $V(C)$), we define $w(D)$ to be $\sum_{e\in D}\w{e}$. 
		\end{enumerate} 
		
	\end{definition}

	\begin{definition}[Paths in generalized shapes]
		Let $\alpha = (\U{\alpha},\V{\alpha},\W{\alpha},E(\alpha))$ be a generalized shape.
		We say a sequence of vertices $v_0,v_1,\dots, v_\ell$ is a path from $\U{\alpha}$ to $\V{\alpha}$ if it satisfies the following:
		\begin{enumerate}
			\item $v_0\in \U{\alpha}$ and $v_\ell\in \V{\alpha}$.
			\item  For each $i=0, 1,2,\dots, \ell-1$, $v_i$ and $v_i+1$ are contained in some hyperedge of $\alpha$.
		\end{enumerate} 
	\end{definition}

	\begin{definition}[Minimum weight separators]
		Let $\alpha = (\U{\alpha},\V{\alpha},\W{\alpha},E(\alpha))$ be a generalized shape. 
		We define a minimum weight separator  $\smin$ to be a subset of vertices that satisfies the followings: 
		\begin{enumerate}
			\item $\smin$ separates $\U{\alpha}$ from $\V{\alpha}$. 
			In other words, each path from $\U{\alpha}$ to $\V{\alpha}$ intersects $\smin$.
			\item For any subset $S$ that separates $\U{\alpha}$ from $\V{\alpha}$, we have $\w{\smin} \leq \w{S}$.
		\end{enumerate}
	\end{definition}	
With the above definitions, we are now ready to state our norm bounds for generalized graph matrices.
We first state and prove the most general form of the norm bounds which can handle any input distributions.
Later, we will see how the norm bounds can be simplified for common input distributions such as the Bernoulli distribution and the normal distribution.

\begin{mdframed}[style=box2]
	\begin{theorem}[Bounds on generalized graph matrices]
		\label{boundgeneral}
	For a given generalized shape $\alpha=(\U{\alpha},\V{\alpha},\W{\alpha},E(\alpha))$, let $\smin$ be a minimum weight vertex separator between  $\U{\alpha}$ and $\V{\alpha}$, and let $W_{iso}$ be the set of isolated vertices in $W_{\alpha}$.
		Further, for each $i\in[t]$, let $m_i$ be the number of vertices of type $i$ excluding those in $\U{\alpha}\cap \V{\alpha}$.
		Adopting notations from Section~\ref{def:gen} and \ref{sec:handle}, the following norm bound holds with probability at least $1-\epsilon$:
\begin{align*}
\norm{M_{\alpha}} &\leq 2\left(\prod_{i=1}^t m_i^{m_i}\right) \cdot  
n^{\frac{w(V(\alpha)) - w(\smin) + w(W_{iso})}{2}} \cdot \\
&\min_{ q \geq 3}\left\{(2q)^{|V(\alpha) \setminus (U_{\alpha} \cap V_{\alpha})|} 
\left(\prod_{e \in E(\alpha)}{h^{+}_{l_{e}}(B_{\Omega}(2q{l_e}))}\right)
 \left(\frac{n^{ \w{\smin} }}{\epsilon} \right)^{\frac{1}{2q}}
\right\}.
		\end{align*}
	\end{theorem}
\end{mdframed}

	\begin{remark} As discussed in the introduction, Theorem~\ref{boundgeneral} can be used to give alternative proofs of several technical lemmas in the literature ~\cite{barak2012hypercontractivity, hopkins2016fast,ge2015decomposing}. Moreover, our proofs based on graph matrices are mechanical, whereas the original analyses often required clever arguments. For details, see Section~\ref{sec:appl}. 
	
\end{remark}

Note that the current form of Theorem~\ref{boundgeneral} is not user-friendly since the bound is written as the minimum over $q\geq 3$.
Below, we write the norm bound for $q=3\left\lceil\frac{\log \left(n^{\w{\smin}}/\epsilon \right)}{3|V(\alpha) \setminus (U_{\alpha} \cap V_{\alpha})|}  \right\rceil$, which will turn out to be a nearly optimal choice for many cases. 
In the subsequent subsections, we will further tune the choice of $q$ for specific distributions. 
\begin{mdframed}[style=box]
    \begin{corollary}[User-friendly version of Theorem~\ref{boundgeneral}]
   Under the setting of Theorem~\ref{boundgeneral}, the following probabilistic upper bound holds:	with probability at least $1-\epsilon$,
   \[
   \norm{M_\alpha} \leq  2\left(\prod_{i=1}^t m_i^{m_i}\right) \cdot \left(\prod_{e \in E(\alpha)}{h^{+}_{l_{e}}(B_{\Omega}(2q{l_e}))}\right)
   (eq)^{|V(\alpha) \setminus (U_{\alpha} \cap V_{\alpha})|} 
n^{\frac{w(\alpha) - w(\smin) + w(W_{iso})}{2}}\,,
   \]
   where 
   \[
   q = 3\left\lceil\frac{\log \left(n^{\w{\smin}}/\epsilon \right)}{3|V(\alpha) \setminus (U_{\alpha} \cap V_{\alpha})|}  \right\rceil\,.
   \]
\end{corollary}
	\end{mdframed}

	\subsection{Proof of Theorem~\ref{boundgeneral}}
	
We first note that the arguments in Section~\ref{reduc:1} and \ref{reduc:2} work similarly here. Similar to Proposition \ref{prop:isolated}, each isolated vertex $v$ gives a multiplicative factor of at most $n^{\w{v}}$, so it is sufficient to consider the case where there are no isolated vertices. When $U_{\alpha} \cap V_{\alpha}$ is non-empty, we can still partition $M_{\alpha}$ into blocks based on the values of the vertices in $U_{\alpha} \cap V_{\alpha}$ and analyze each block separately (see Section \ref{reduc:1}). To handle hyperedges $e$ with vertices in $U_{\alpha} \cap V_{\alpha}$, we can simply delete these vertices from $e$ (viewing $e$ as a hyperedge on its remaining vertices) because in each block of $M_{\alpha}$, vertices in $U_{\alpha} \cap V_{\alpha}$ are fixed. Thus, it suffices to prove Theorem~\ref{boundgeneral} for the case where $\U{\alpha}\cap \V{\alpha} = \emptyset$ and $\iso = \emptyset$. Hence, we assume throughout the proof that $\U{\alpha}\cap \V{\alpha} = \emptyset$ and $\iso = \emptyset$.

Following the previous proofs, the main step of the proof is to prove upper bounds on the trace power terms.
Akin to Section~\ref{pf:main}, Theorem~\ref{boundgeneral} will then be an immediate consequence of the vertex partitioning lemma (Lemma~\ref{vpl}).

	\begin{lemma} 	\label{traceupper4}
		Let $\smin$ be a minimum weight vertex separator between  $\U{\alpha}$ and $\V{\alpha}$. Then, for all $P\in \mathcal{P}_\alpha$ and $q \in \mathbb{N}$, \[	\E\left[\tr\left((M_{\alpha,P}M_{\alpha,P}^\top)^q\right)\right] \leq 
		\left(\prod_{e \in E(\alpha)}{h^{+}_{l_{e}}(B_{\Omega}(2q{l_e}))^{2q}}\right){(2q)^{2q|V(\alpha)|}}n^{q(w(\alpha) - \w{\smin}) + \w{\smin}}\,.
		\]
	\end{lemma}

	\begin{proof}[Proof of Lemma~\ref{traceupper4}]
		Recall from Proposition~\ref{expan3} that
		\[
		\E\left[\tr\left((M_{\alpha,P}M_{\alpha,P}^\top)^q\right)\right] = \sum_{C \in \mathcal{C}_{(\alpha,2q)}: C \text{ is well-behaved}}{N_{P}(C)\val{C}}\,.
		\]
	By Corollary \ref{cor:valCbound}, for any well-behaved constraint graph $C \in \mathcal{C}_{(\alpha,2q)}$, 
		\[
		\val{C} \leq \prod_{e \in E(\alpha)}{h^{+}_{l_{e}}(B_{\Omega}(2q{l_e}))^{2q}}
		\]

		Now, let us upper bound $N_P(C)$. 
		First,  we generalize  Proposition~\ref{numchoice} to cover the case with vertices of multiple types:

		\begin{proposition} \label{numchoice2}
			Let $C \in \mathcal{C}_{(\alpha,2q)}$ be a well-behaved constraint graph. 
			Then, for 	 $n=\max_{i\in[t] }|\N{i}|$,  we have  $N_P(C) \leq n^{\w{V(C)} -\w{E(C)}}$.
		\end{proposition}
		\begin{proof}
			Observe that choosing a type-respecting, $P$-respecting and injective map $\phi: V(\alpha,2q) \to \N{} $ with constraint graph $C$ is equivalent to choosing the distinct indices of $\phi(V(\alpha,2q))$.
			For each vertex $v\in V(\alpha)$, let $C_v$ be the subset of $C$ that consists of the copies of $v$, and let $E_v$ be the set of constraint edges between $C_v$.
			Note that $E(C) = \cup_{v\in V(\alpha)} E_v$ due to well-behavedness of $C$.
			
			Then, it is straightforward to see that the number of different choices for $\phi(C_v)$ is at most $|\N{i}|^{|C_v|-|E_v|} = n^{\w{C_v}-\w{E_v}}$, where $i$ is the type of $v$. 
			Taking a product of these inequalities for all $v\in V(\alpha)$, we obtain the result.
		\end{proof}

		Thus, Lemma~\ref{traceupper4} follows directly from the following two statements:
		
		\begin{enumerate}
			\item For any well-behaved $C \in \mathcal{C}_{(\alpha,2q)}$ such that $\val{C} \neq 0$, the total weight of the constraint edges is at least $q\cdot \w{\W{\alpha}}+ (q-1)\cdot \w{\smin}$.
			\item There are at most $(2q)^{2q|V(\alpha)|}$ well-behaved constraint graphs $C \in \mathcal{C}_{(\alpha,2q)}$.
		\end{enumerate} 
		
		Now we prove the two statements. We begin with the second statement.
		
		\begin{lemma} \label{firststate}
			There are at most $(2q)^{2q|V(\alpha)|}$ well behaved constraint graphs $C \in \mathcal{C}_{(\alpha,2q)}$.
		\end{lemma}
		\begin{proof}
			To specify a constraint graph, it is sufficient to take each $v \in V(\alpha,2q)$ and specify whether $v$ is equal to a later vertex and if so, what is the next vertex $w \in V(\alpha,2q)$ such that $v = w$. Since $|V(\alpha,2q)| \leq 2q|V(\alpha)|$ and there are at most $2q$ possibilities for each vertex ($2q-1$ possible vertices $w$ and the possibility that $v$ is not equal to a later vertex), the result follows.
		\end{proof}
		\begin{remark}
			A careful reader might notice that the proof of the above lemma is based on a much cruder counting than those appearing in the proof of  Lemmas~\ref{numconstr} and \ref{numconstr2}.
				Note that the counting arguments in the proof of  Lemmas~\ref{numconstr} and \ref{numconstr2} do not readily carry over to this case due to the presence of weights in the generalized constraint graph.
				It would be interesting to investigate whether we can do a better counting and come up with a better bound.
				On the other hand, for the purpose of our norm bounds, this slack only results in an additional poly-logarithmic factor in the norm bounds.
		\end{remark}
		We now prove the first statement. The proof bears some resemblances to that of Lemma~\ref{mingeneral}.
		However, more sophisticated arguments are needed as we are in the weighted case.
		We remark that there is an alternative proof to this statement which we present in Appendix~\ref{altproof}.
		\begin{lemma}\label{mingeneral2}
			For any well-behaved $C \in \mathcal{C}_{(\alpha,2q)}$ such that $\val{C} \neq 0$, the total weight of the constraint edges is at least $q\cdot \w{\W{\alpha}}+ (q-1)\cdot  \w{\smin } $.
		\end{lemma}
		\begin{proof}
			First, note that in the constraint graph, every vertex that is a copy of a vertex in $\W{\alpha}$ must be incident to at least one constraint edge as we have $\val{C} =0$ otherwise.
			Since there are $2q$ copies of $\W{\alpha}$, such constraint edges give a total weight of $\frac{2q\cdot \w{\W{\alpha}}}{2} = q\cdot \w{\W{\alpha}}$. Thus, to prove Lemma~\ref{mingeneral2}, we must show that there must be ``additional'' constraint edges with total weight at least $(q-1)\cdot \w{\smin}$.

			For each vertex $v$ in the shape $\alpha$,  let $\wei{v}$ be the number of additional constraint edges between the copies of $v$. 
			Note that $0\leq \wei{v} \leq q-1$ for any $v\in V(\alpha)$ and also, $\sum_{v\in V(\alpha)}\wei{v}$ is equal to the total number of constraint edges due to $C$ being well-behaved. 
			Hence, the Lemma~\ref{mingeneral2} statement can be rephrased as 
			\begin{align} \label{wts}
			    \sum_{v\in V(\alpha)} \w{v}\wei{v}  \geq  (q-1)\cdot \w{\smin}\,.
			\end{align}

		Now consider any path $P$ in $\alpha$ from $\U{\alpha}$ to  $\V{\alpha}$ whose vertices lie in $\W{\alpha}$ except for the endpoints.
	   Then we can repeat the  argument from the proof of  Lemma~\ref{mingeneral} to show $\sum_{v\in P} \wei{v} \geq q-1$:
	    \begin{itemize}
	        \item[] Let $l$ be the length of $P$. Consider the cycle $T$ which consists of the $2q$ copies of $P$ in the constraint graph. Then, $T$ has length $2ql$.  Applying Lemma~\ref{mink-1} to $T$, we conclude that there must be at least $ql-1$ constraint edges within $T$.
	        However, not all of them are ``additional'' constraint edges that we mentioned earlier.
            As we noted above, each vertex that is a copy of a vertex in $\W{\alpha}$ must be incident to a constraint edge.
	        Therefore, in terms of the number of ``additional'' constraint edges, there are $ql-1 -q(l-1) = q-1$ of them.
	    \end{itemize}

In order to show \eqref{wts}, we formulate an integer programming problem where we minimize $\sum_{v\in V(\alpha)} \w{v}\wei{v}$ subject to the constraints that $\wei{v}\in\{0,1,\dots, q-1\} $ and $\sum_{v\in P} \wei{v} \geq q-1$ for any path from $\U{\alpha}$ to $\V{\alpha}$.
			Then it suffices to show that the objective value of this  integer programming problem is lower bounded by  $(q-1)\cdot  \w{\smin}$.
			To simplify our consideration, define $\weii{v}:=\wei{v}/(q-1)$ and relax the constraints $\weii{v}\in\{0,\frac{1}{q-1},\dots, 1\} $ to $\weii{v}\in [0,1]$.
			The resulting linear program is 
			\begin{align} \label{lp}
			\min_{\weii{v}\in[0,1]}\left\{ \sum_{v\in V(\alpha)} w(v)\weii{v}~:~ \sum_{v\in P} \weii{v} \geq 1~\text{for any path from $\U{\alpha}$ and $\V{\alpha}$} \right\}\,.
			\end{align}
			Now it suffices to show that the objective value of the above linear program is lower bounded by  $ \w{\smin}$.
			\begin{proposition}
				There is an integral optimal solution to the above linear program.
			\end{proposition}
			\begin{proof}
				Let $\{\weii{v}^*\}$ be  an optimal solution to the program.
				We will show that there is an integral solution whose objective value is less than or equal to that of $\{\weii{v}^*\}$. 
				For any path, let us define the weight of the path to be the sum of  the weights of vertices  along the path.
				For each vertex $v$ in $\alpha$, let $d(v)$ be the minimum weight among the paths from $\U{\alpha}$ to $v$.
				For instance, $d(u)=x^*_u$ for any $u\in \U{\alpha}$.
				\begin{claim}\label{claim:int}
					For each $v \in \V{\alpha}$, $d(v)=1$.
				\end{claim}
				\begin{proof}[Proof of Claim~\ref{claim:int}]
					Due to the fact that $\{\weii{v}^*\}$ satisfies the constraint, we have $d(v)\geq 1$ for any $v \in \V{\alpha}$.
					Now suppose to the contrary that there is a vertex $w\in \V{\alpha}$ such that $d(w)>1$, say $d(w)=1+\epsilon$ for some constant $\epsilon>0$.
					Then, there exists a path $P$ from $u\in \U{\alpha}$ to $w$ with weight $1+\epsilon$.
					Let $w'$ be the last vertex on $P$ with nonzero $\weii{}^*$-value, i.e., $\weii{w'}^*>0$.
					Then, it should be the case that $d(w')=1+\epsilon$.  
					Note that slightly decreasing $\weii{w'}^*$ does not conflict with any constraints from the program; the only case when it could possibly affect is when $w'$ was on a path with weight exactly equal to $1$, which contradicts the fact that  $d(w')=1+\epsilon$.
					Hence, one can slightly decrease the objective value by decreasing $\weii{w'}^*$, which  contradicts the optimality of $\{\weii{v}^*\}$.
				\end{proof}
				Now, for each $t\in [0,1]$, let $S_t:= \left\{  u ~:~ t \in [d(u)-\weii{u}^*,d(u)] \right\}$.
				Note that $S_t$ is a  vertex separator between $\U{\alpha}$ and $\V{\alpha}$ for each $t\in[0,1]$.
				Indeed, for any path $P$, let $w$ be the first vertex along the path satisfying $d(w)\geq t$ as we read it from $\U{\alpha}$ to $\V{\alpha}$.
				If it is the starting vertex of $P$ in $\U{\alpha}$, then $t \in [0, d(w)]$, which implies  $w\in S_t$.
				Otherwise, let $w'$ be the vertex right before $w$ on $P$.
				Then, $d(w')\geq d(w)-\weii{w}^*$, otherwise there is a path to $w$ with weight smaller than $d(w)$.
				Since $d(w')< t$ by definition, we have $d(w)-\weii{w}^* \leq d(w') < t \leq d(w)$, again implying that $w\in S_t$.
				
				Hence, letting $\weii{v}^{(t)}:= \mathbb{I}_{\{ v \in S_t\}}$, we conclude that $\{\weii{v}^{(t)}\}$ is a feasible solution for each $t\in[0,1]$.  Consider a uniform random variable $T$ in $[0,1]$.
				By the construction, we have $\pP(\weii{v}^{(T)}=1)= \weii{v}^*$, which implies $\ex_{t\sim T} \left[\sum_{v\in V(\alpha)} w(v) \weii{v}^{(t)}\right] = \sum_{v\in V(\alpha)} w(v)\pP(\weii{v}^{(T)}=1) = \sum_{v\in V(\alpha)} w(v) \weii{v}^*$.
				Therefore, there must exist $t\in[0,1]$ such that the objective value of $\{\weii{v}^{(t)}\}$ is at least that of  $\{\weii{v}^*\}$, implying the existence of an integral optimal solution. 
			\end{proof}
	Now note that for any integral feasible  solution $\{\weii{v}\}$, the set $\{v\in V(\alpha)~:~\weii{v}=1 \}$ is a seperator between $\U{\alpha}$ and $\V{\alpha}$.
	This is because any path from $\U{\alpha}$ to $\V{\alpha}$ must intersect the set due to the constraint.
	Therefore, the objective value of the linear program \eqref{lp} is lower bounded by $\w{\smin}$, which  completes the proof of Lemma~\ref{mingeneral2}.
		\end{proof}
		
		Since the two statements (Lemma~\ref{firststate} and \ref{mingeneral2}) are proved,  Lemma~\ref{traceupper4} follows.
	\end{proof}
	Having established Lemma~\ref{traceupper4}, Theorem~\ref{boundgeneral} follows from the vertex partitioning lemma (Lemma~\ref{vpl}) together with Corollary~\ref{cor:vpa}:
	more precisely, we apply the vertex partitioning lemma (Lemma~\ref{vpl}) with $M\leftarrow \frac{1}{\prod_{i=1}^t{|V_i(\alpha)|^{|V_i(\alpha)|}}}\M_\alpha$, $\{M_i\}\leftarrow \{\M_{\alpha,P}\}$, and  
	$$	B(2q)\leftarrow\left(\prod_{e \in E(\alpha)}{h^{+}_{l_{e}}(B_{\Omega}(2q{l_e}))^{2q}}\right){(2q)^{2q|V(\alpha)|}}n^{q(w(\alpha) - \w{\smin}) + \w{\smin}}\,.$$
	This completes the proof of Theorem~\ref{boundgeneral}.

	\subsection{Explicit bounds for special cases}
	We now give more explicit and user-friendly bounds for two  common input distributions: (i) the $\pm 1$ Bernoulli distribution and (ii) the normal distribution.
	\subsubsection{The \texorpdfstring{$\pm{1}$}{+1} Bernoulli distribution case}
	For the $\pm 1$ Bernoulli distributions, it can be easily verified that $B_{\Omega} = 1$.
    Moreover, since the support is of size $2$, we have only two elements $h_0$ and $h_1$ in the orthonormal basis.
Here, one can easily check that  $h^{+}_{1}(1) = 1$. 
    Using these explicit calculations of the input distribution related quantities, one can state a refined version of   Theorem~\ref{boundgeneral} as follows.
    
	\begin{mdframed}[style=box]
		\begin{corollary}[Refined bounds for Bernoulli distributions] \label{cor:bern}
		 Under the setting of Theorem~\ref{boundgeneral}, assume further that the input distribution $\Omega$ takes value between $+1$ and $-1$, i.e. it is a $\pm 1$-Bernoulli distribution.
		 Then, the following probabilistic upper bound holds	with probability at least $1-\epsilon$:
		\begin{align*}
\norm{M_{\alpha}} &\leq 2\left(\prod_{i=1}^t m_i^{m_i}\right) \cdot  \left(6e\left\lceil  \frac{\log ( \frac{n^{\w{\smin}}}{\epsilon})}{6|\V{\alpha}\setminus (\U{\alpha}\cap \V{\alpha})|} \right \rceil\right)^{|\V{\alpha}\setminus (\U{\alpha}\cap \V{\alpha})|}
n^{\frac{w(V(\alpha)) - w(\smin) + w(W_{iso})}{2}} \,.
		\end{align*}
	\end{corollary}	
	\end{mdframed}
	\begin{proof}
		Due to the above observations that  $B_{\Omega} = 1$ and  $h^{+}_{1}(1) = 1$,  the term we need to minimize over $q\geq 3$ is equal to  $$(2q)^{|\V{\alpha}\setminus (\U{\alpha}\cap \V{\alpha})|} \left(\frac{n^{\w{\smin}}}{\epsilon}\right)^{\frac{1}{2q}}\,.$$
		Taking the natural logarithm of the above gives 
		$|\V{\alpha}\setminus (\U{\alpha}\cap \V{\alpha})|\cdot \log(2q) + \frac{1}{2q}\cdot \log ( \frac{n^{\w{\smin}}}{\epsilon})$, which after taking the derivative  with respect to $q$ becomes  $ |\V{\alpha}\setminus (\U{\alpha}\cap \V{\alpha})|\frac{1}{q} - \frac{1}{2q^2}\cdot \log ( \frac{n^{\w{\smin}}}{\epsilon})$.
        This derivative is equal to  $0$ when $q = \frac{\log ( \frac{n^{\w{\smin}}}{\epsilon})}{2|\V{\alpha}\setminus (\U{\alpha}\cap \V{\alpha})|}$. 
        This suggests that we take   $q = 3\left\lceil  \frac{\log ( \frac{n^{\w{\smin}}}{\epsilon})}{6|\V{\alpha}\setminus (\U{\alpha}\cap \V{\alpha})|} \right \rceil \geq 3$ (again the extra multiplication by $3$ is for ensuring $q\geq 3$), which gives Corollary~\ref{cor:bern}.
		\end{proof}
		
		\subsubsection{The normal distribution case}
	
		Again, to come up with a refined bound for  normal distributions, we need to first compute $B_{\Omega}$ and $h^{+}_k$ for $\Omega = N(0,1)$.
        The following lemma characterizes these quantities:		
	
		\begin{lemma}\label{Gaussianbounds}
			If $\Omega = N(0,1)$ is the normal distribution then we can take $B_{\Omega}(j) = \sqrt{j}$ and we have that for all $k \in \mathbb{N}$ and all $x \in \mathbb{R}$, 
			$h^{+}_{k}(x) \leq \frac{1}{\sqrt{k!}}(x^2 + k)^{\frac{k}{2}} \leq \left(\frac{e}{k}(x^2 + k)\right)^{\frac{k}{2}}$.
		\end{lemma}
		\begin{proof}
			For the first part, recall that $E_{\Omega}[x^j] = \frac{(2j)!}{j!2^{j}} = \prod_{i=0}^{\frac{j}{2}-1}{(2i+1)}$ if $j$ is even and $E_{\Omega}[x^j] = 0$ if $j$ is odd. Thus, for all $j \in \mathbb{N}$, 
			$\left|E_{\Omega}[x^j]\right| \leq j^{\frac{j}{2}} = (\sqrt{j})^{j}$.
			
			For the second part, recall that the orthonormal basis for $\Omega = N(0,1)$ is 
			\[
			h_k(x) = \frac{1}{\sqrt{k!}}{\sum_{j=0}^{\lfloor\frac{k}{2}\rfloor}{(-1)^{j}\frac{k!}{(k-2j)!(2j)!} \cdot \frac{(2j)!}{j!2^{j}}x^{k-2j}}}\,.
			\]
			Thus,
			\begin{align*}
			h^{+}_k(x) &= \frac{1}{\sqrt{k!}}{\sum_{j=0}^{\lfloor\frac{k}{2}\rfloor}{\frac{k!}{(k-2j)!j!2^{j}}x^{k-2j}}} \\
			&\leq \frac{1}{\sqrt{k!}}\sum_{j=0}^{\lfloor\frac{k}{2}\rfloor}{\frac{(\lfloor\frac{k}{2}\rfloor)!}{(\lfloor{\frac{k}{2}}\rfloor-j)!j!}k^{j}(x^2)^{\frac{k}{2}-j}} \\
			&= \frac{1}{\sqrt{k!}}(x^2 + k)^{\lfloor\frac{k}{2}\rfloor}(x^2)^{\frac{k}{2} - \lfloor\frac{k}{2}\rfloor} \leq \frac{1}{\sqrt{k!}}(x^2 + k)^{\frac{k}{2}}\,. 
			\end{align*}
		\end{proof}  
        \noindent Using  Lemma~\ref{Gaussianbounds}, we can  refine Theorem~\ref{boundgeneral}  as follows:        
        \begin{mdframed}[style=box]
		\begin{corollary}[Refined bounds for normal distributions] \label{cor:normal}
		 Under the setting of Theorem~\ref{boundgeneral}, assume further that the input distribution $\Omega$  is a normal distribution (Gaussian distribution).
		 Let $l(\alpha):= \sum_{e\in E(\alpha)} l_e$.
		 Then, the following probabilistic upper bound holds	with probability at least $1-\epsilon$:
		\begin{align*}
\norm{M_{\alpha}} \leq& 2\left(\prod_{i=1}^t m_i^{m_i}\right) \cdot 
n^{\frac{w(V(\alpha)) - w(\smin) + w(W_{iso})}{2}}\cdot  \\
& \left(6e \left\lceil{\frac{\log (\frac{n^{\w{\smin}}}{\epsilon})}{6(|V(\alpha) \setminus (U_{\alpha} \cap V_{\alpha})| + l(\alpha))}}\right \rceil\right)^{l(\alpha)+|V(\alpha) \setminus (U_{\alpha} \cap V_{\alpha})|}
		\end{align*}
	\end{corollary}	
	\end{mdframed}
        
        	\begin{proof}
		By Lemma \ref{Gaussianbounds}, we then have the following for each $k$:
		\begin{align*}
		   h^{+}_{k}(B_{\Omega}(2qk)) &\leq \left(\frac{e}{k}(B_{\Omega}(2qk)^2+k) \right)^{k/2}\\
		   &=\left(\frac{e}{k}(2qk+k) \right)^{k/2}\\
		   &= e^{k/2} (2q+1)^{k/2}\,.
		\end{align*}
		Plugging this into Theorem \ref{boundgeneral}, the term we need to minimize over $q$ becomes (using the notation $l(\alpha):= \sum_{e\in E(\alpha)} l_e$):
        \begin{align*}
         &(2q)^{|V(\alpha) \setminus (U_{\alpha} \cap V_{\alpha})|} 
\left(\prod_{e \in E(\alpha)} e^{l_e/2}(2q+1)^{l_e/2}\right)
 \left(\frac{n^{ \w{\smin} }}{\epsilon} \right)^{\frac{1}{2q}}\\
 =&(2q)^{|V(\alpha) \setminus (U_{\alpha} \cap V_{\alpha})|}   e^{l(\alpha)/2} (2q+1)^{l(\alpha)/2}
 \left(\frac{n^{ \w{\smin} }}{\epsilon} \right)^{\frac{1}{2q}}\,.
        \end{align*}
        Using the elementary inequality that $x+1\leq x^2$ for $x\geq 2$, the above term then becomes:
         \begin{align*}
         (2q)^{l(\alpha)+|V(\alpha) \setminus (U_{\alpha} \cap V_{\alpha})|}   e^{l(\alpha)/2} 
 \left(\frac{n^{ \w{\smin} }}{\epsilon} \right)^{\frac{1}{2q}}\,.
        \end{align*} 
		Now choosing  $q = 3\left\lceil{\frac{\log (\frac{n^{\w{\smin}}}{\epsilon})}{6(|V(\alpha) \setminus (U_{\alpha} \cap V_{\alpha})| + l(\alpha))}}\right \rceil \geq 3$ (again the extra multiplication by $3$ is for ensuring $q\geq 3$), we recover Corollary~\ref{cor:normal}.
	\end{proof}

	\section{Applications of graph matrix norm bounds} \label{sec:appl}

It was already shown in prior works~\cite{Norms,DBLP:journals/corr/DeshpandeM15, DBLP:journals/corr/HopkinsKP15,DBLP:journals/corr/RaghavendraS15, FinalPlantedClique,PAPSOS} that graph matrices are very useful in proving lower bounds on the Sum-of-Squares hierarchy.
In this section, we describe other applications of graph matrices: 
it turns out one can also capture various \emph{upper bound} analyses for the Sum-of-Squares hierarchy and other spectral algorithms with generalized graph matrices.
Here we highlight that our proofs based on graph matrices are \emph{mechanical}, whereas the original  analyses often required  clever arguments.

	\subsection{Warm-up: How can graph matrices be applied?}
	In many prior works on the sum of squares hierarchy and other spectral algorithms~\cite{barak2012hypercontractivity,hopkins2016fast,ge2015decomposing,barak2016noisy,potechin2017exact,raghavendra2017strongly}, the  main step of the analysis was showing some upper bounds on the maximum of a homogeneous polynomial over the unit sphere, an important theoretical question of independent interest extensively studied in \cite{Bhattiprolu2016SumofSquaresCF,bhattiprolu2017weak}.
	To recover such upper bounds using graph matrices, our approach will hinge upon matrix representations of homogeneous polynomials: 
	\begin{definition}
		For a degree $d$ (assume $d$ is even) homogeneous polynomial $f$ on $n$ variables $x_1,\dots, x_n$, we say a $n^{d/2}\times n^{d/2}$ matrix $M$ is a \emph{matrix representation} of $f$ if $f(x_1,\dots, x_n) = (x^{\otimes d/2})^\top Mx^{\otimes d/2}$ for all $x\in \re^n$.
		Moreover, we say a matrix representation $M$ is a \emph{graph matrix representation} if $M$ is, in addition, a graph matrix.
		If $M$ is a graph matrix  corresponding to a shape $\alpha$, i.e., $M=M_{\alpha}$, then we say $\alpha$ is a \emph{shape representation} of $f$.
	\end{definition}
	For an even number $d$, let $f$ be a degree-$d$ homogeneous polynomial that we want to prove an upper bound for.
	Assuming that $f$ has a graph matrix representation, one can easily find an upper bound on $\max_{\norm{x}=1}f(x)$.
	In particular,  if $\alpha$ is a shape representing $f$,  we have $f(x) = (x^{\otimes d/2})^\top M_{\alpha}x^{\otimes d/2} \leq \norm{M_{\alpha}}$, which implies $\max_{\norm{x}=1}f(x) \leq \norm{M_{\alpha}}$.
	We can then invoke  our main theorems to upper bound $\norm{M_{\alpha}}$.
Therefore, to upper bound $f$, one needs to find a shape representation of $f$.
We illustrate how one can find a shape representation through the following example. 
	(This  example will reappear in Sec.~\ref{fast:sparse}.)
One remark on notation is that we will use $u_1,u_2,\dots$ to denote the vertices in $\U{}$,   $v_1,v_2,\dots$ to denote the vertices in $\V{}$,  $\nn_1,\nn_2,\dots$ to denote the vertices in $\U{}\cap \V{}$, and  $w_1,w_2,\dots$  to denote the vertices in $W$.	
	\begin{example}[An example adapted from \cite{hopkins2016fast}] \label{ex1}
		For integers $n,m$ such that $n \leq m$, 
		let $\{b_{i,j}\}_{i\in[n],~j\in[m]}$ be  random variables drawn i.i.d. from $\Nc(0,1 )$.
		Let $f_{\text{Ex}1}(x_1,\dots ,x_n)$ be a degree-$2$ homogeneous polynomial defined as 
		\begin{align*}
		f_{\text{Ex}1}(x_1,\dots ,x_n):= \sum_{\substack{j_1\neq j_2 \in [n]}} x_{j_1} x_{j_2}\left[\sum_{\substack{i\in [m],~j_3\in [n]\\
				j_3\notin\{j_1,j_2\}}} b_{i,j_1}b_{i,j_2} \left(   b_{i,j_3}^2 -1 \right)  \right]\,.
		\end{align*} 
	\end{example}
	Let us find a shape $\alpha_{\text{Ex}1}$ that represents $f_{\text{Ex}1}$. 
	The  principal observations are as follows:
	\begin{mdframed}
	\begin{itemize}
	    \item Since $j_1,j_2,j_3$ are distinct, there have to be four vertices in  $\alpha_{\text{Ex}1}$, each representing a different index among $i,j_1,j_2,j_3$. 
	    \item Since the variable of each summand is equal to $x_{j_1}x_{j_2}$,
the graph matrix representation will have dimension $n\times n$.
Moreover,  for the two vertices corresponding to $j_1$ and $j_2$, one of them should in $\U{}$ and the other is in $\V{}$.
By symmetry, one may assume that $u_1\in\U{}$ and $v_1\in\V{}$.
\item The other two vertices corresponding to $i$ and $j_3$ belong to $\W{}$. Let us denote those vertices by $w_1$ and $w_2$, respectively.
	\end{itemize}
	\end{mdframed}
	Accounting for the fact that $i\in[m]$ and $j_1,j_2,j_3\in[n]$,  we follow the setting in Section~\ref{pf:generalgeneral} and assign  weight $\log_N (m)$ to $w_1$ and assign weight $\log_N (n)$ to each of $u_1,v_1,w_2$, where $N:=\max\{n,m\}$.
	In view of Example~\ref{ex:fourier},  random variables $b_{i,j_1}$, $b_{i,j_2}$ and $  (b_{i,j_3}^2 -1)$ are (constant multiples) of the orthonormal basis of $\Nc(0,1 )$, and hence, one can represent those random variables by edges $\{w_1,u_1\}$, $\{w_1,v_1\}$, and $\{w_1,w_2\}$ with labels  $l_{\{w_1,u_1\}}=1$, $l_{\{w_1,v_1\}} = 1$, and $l_{\{w_1,w_2\}}=2$, respectively.
	Hence, one can see that $f_{\text{Ex}1}(x)=x^\top M_{\alpha_{\text{Ex}1}}x$, 
	where $\alpha_{\text{Ex}1}$ is a shape on four vertices $u_1,v_1,w_1,w_2$ with weights $\w{u_1}=\w{v_1}=\w{w_2}=\log_N (n)$, $\w{w_1}=\log_N(m)$ and three edges $\{u_1,w_1\}, \{v_1,w_1\}, \{w_1,w_2\}$ with labels $l_{\{w_1,u_1\}}=1$, $l_{\{w_1,v_1\}} = 1$, and $l_{\{w_1,w_2\}}=2$.
	See Figure~\ref{warmup1}. 
	Here and below, the edges are of label $1$ unless specified otherwise.
	
	\begin{figure}[H]
		\centering
		\scalebox{0.55}{
			\begin{tikzpicture}[one/.style={circle,fill=blue!25,draw,font=\sffamily\LARGE\bfseries},two/.style={diamond,fill=red!5,draw,font=\sffamily\LARGE\bfseries}]
			\node[one] (u1) at (-4,2)  {$u_1=j_1$};
			\node[two] (w1) at (0,2)  {$w_1=i$};
			\node[one] (v1) at (4,2)  {$v_1=j_2$};
			\node[one] (w2) at (0,-2) {$w_2=j_3$};
			\path[every node/.style={font=\sffamily\LARGE\bfseries},line width=1pt]
			(u1) edge  (w1) 
			(w1) edge    (v1)
			(w1) edge node[right]{$2$}  (w2);
				\draw[draw=black] ($(v1.135)+(-0.5,0.5)$) rectangle  ($(v1.-45)+(0.5,-0.5)$);
		\draw[draw=black] ($(u1.135)+(-0.5,0.5)$) rectangle  ($(u1.-45)+(0.5,-0.5)$);
			\node[draw=none,fill=none] (8) at (-4,4.5) {\fontsize{20}{22.4}\selectfont$\U{}$};
			\node[draw=none,fill=none] (10) at (4,4.5) {\fontsize{20}{22.4}\selectfont$\V{}$};
			\node[draw=none,fill=none] (10) at (0,4.5) {\fontsize{20}{22.4}\selectfont$\W{}$};
			\end{tikzpicture}
		}
		
		\caption{an illustration of how shape $\alpha_{\text{Ex}1}$ represents $f_1$. We depict the nodes with weight $\log_N(n)$ as circles and a node with weight $\log_N(m)$ as a diamond. Here node $u_1$ represents index $j_1$, $v_1$ represents $j_2$, $w_1$ represents $i$, and $w_2$ represents $j_3$.
			Consequently, the shape represents the random variable $b_{i,j_1}$ with the edge $\{w_1,u_1\}$, $b_{i,j_2}$ with the edge $\{w_1,v_1\}$, and $ (b_{i,j_3}^2 -1) $ with the edge $\{w_1,w_2\}$. 
		The edges have labels $l_{\{w_1,u_1\}}=1$, $l_{\{w_1,v_1\}} = 1$ and $l_{\{w_1,w_2\}}=2$.
		Here and below, the edges are of label $1$ unless specified otherwise.
		}
		\label{warmup1}
	\end{figure}
	
	Hence, we have $\max_{\norm{x}=1}f_{\text{Ex}1}(x) \leq \norm{M_{\alpha_{\text{Ex}1}}}$. 
	To obtain a norm bound, due to Corollary~\ref{cor:normal}, it suffices to compute a minimum weight separator of the shape $\alpha$.  
	Since we are in the regime $n \leq m$,  $\{u_1\}$ is a minimum weight separator, so $\max_{\norm{x}=1}f_{\text{Ex}1}(x)\leq \norm{M_{\alpha_{\text{Ex}1}}} = \tO{n\sqrt{m}}$.  \qed

 Example~\ref{ex1} illustrates how we can find a graph matrix representation for a given homogeneous polynomial. 
	However, in some cases, there could be more than one shape representing a single polynomial, as one can see from the following example. (This example will  reappear in Sec.~\ref{overtensor}).

	\begin{example}[An example excerpted from \cite{ge2015decomposing}] \label{ex2}
		For integers $n,m$ such that $n\leq m \leq c\cdot n^{1.5}$ for some constant $c$,  
		let $\{b_{i,j}\}_{i\in[m],~j\in[n]}$ be  i.i.d. Rademacher random variables.
		Let $f_{\text{Ex}2}(x_1,\dots ,x_n)$ be a degree-$4$ homogeneous polynomial defined as 
		\begin{align*}
		f_{\text{Ex}2}(x_1,\dots ,x_n):=\sum_{\substack{j_1,j_2,j_3,j_4\in[n]:\\ j_1,j_2,j_3,j_4\text{ are distinct}}} \left[ \sum_{\substack{i_1\neq i_2\in [m],~ j_5\in[n]\\j_5\notin  \{ j_1,j_2,j_3,j_4 \}   } } b_{i_1,j_5} b_{i_2,j_5} b_{i_1,j_1} b_{i_1,j_2} b_{i_2,j_3} b_{i_2,j_4}
		x_{j_1}x_{j_2}x_{j_3}x_{j_4} \right]\,.
		\end{align*} 
	\end{example}
	Following similar reasoning to Example~\ref{ex1}, one can come up with a graph matrix representation.
However, one distinction is that for this case, there are multiple shapes representing $f_{\text{Ex}2}$.
	We can easily verify that the following two shapes both represent $f_{\text{Ex}2}$ (here since the input distribution is Rademacher, we do not need to consider labels of edges): 
	\begin{itemize}
		\item $\alpha_{\text{Ex}2}$ is a graph on $7$ vertices $u_1, u_2,  v_1,v_2, w_1,w_2,w_3$ with weights $\w{u_1}=\w{u_2}=\w{v_1}=\w{v_2}=\w{w_2}=\log_N(n)$, $\w{w_1}=\w{w_3}=\log_N(m)$ and edges $\{\{u_1,w_1\}$, $\{u_2,w_3\}$, $\{w_1,w_2\}$, $\{w_2,w_3\}$, $\{w_1,v_1\}$, $\{w_3,v_2\}  \}$. See Figure~\ref{case1a:t}.
		\item  $\beta_{\text{Ex}2}$ is a graph on $7$ vertices $u_1, u_2,  v_1,v_2, w_1,w_2,w_3$ with weights $\w{u_1}=\w{u_2}=\w{v_1}=\w{v_2}=\w{w_2}=\logn{n}$, $\w{w_1}=\w{w_3}=\logn{m}$ and edges $\{\{u_1,w_1\}$, $\{u_2,w_1\}$, $\{w_1,w_2\}$, $\{w_2,w_3\}$, $\{w_3,v_1\}$, $\{w_3,v_2\}  \}$. See Figure~\ref{case1b:t}.
	\end{itemize}
	
	\begin{figure}[H]
		\centering
		\begin{subfigure}[H]{0.24\textwidth}
			\centering
			\scalebox{0.5}{
				\begin{tikzpicture}[one/.style={circle,fill=blue!25,draw,font=\sffamily\LARGE\bfseries},two/.style={diamond,fill=red!5,draw,font=\sffamily\LARGE\bfseries}]
				\node[one] (u1) at (-3,2)  {$u_1$};
				\node[one] (u2) at (-3,-2)  {$u_2$};
				\node[two] (w1) at (0,2)  {$w_1$};
				\node[one] (w2) at (0,0)  {$w_2$};
				\node[two] (w3) at (0,-2)  {$w_3$};
				\node[one] (v1) at (3,2)  {$v_1$};
				\node[one] (v2) at (3,-2)  {$v_2$};
				\path[every node/.style={font=\sffamily},line width=1pt]
				(u1) edge  (w1)
				(u2) edge  (w3)
				(w1) edge  (w2)
				(w2) edge  (w3)
				(w1) edge (v1)
				(w3) edge (v2);
						\draw[draw=black] ($(v1.135)+(-0.3,0.3)$) rectangle  ($(v2.-45)+(0.3,-0.3)$);
		\draw[draw=black] ($(u1.135)+(-0.3,0.3)$) rectangle  ($(u2.-45)+(0.3,-0.3)$);
			\node[draw=none,fill=none] (8) at (-3,3.5) {\fontsize{20}{22.4}\selectfont$\U{}$};
			\node[draw=none,fill=none] (10) at (3,3.5) {\fontsize{20}{22.4}\selectfont$\V{}$};
			\node[draw=none,fill=none] (10) at (0,3.5) {\fontsize{20}{22.4}\selectfont$\W{}$};
				\end{tikzpicture}
				
			}
			\caption{Figure~\ref{case1a:t}: $\alpha_{\text{Ex}2}$.}
			\label{case1a:t}
		\end{subfigure}
		\quad\quad\quad\quad
		\begin{subfigure}[H]{0.24\textwidth}
			\centering
			\scalebox{0.5}{
				\begin{tikzpicture}[one/.style={circle,fill=blue!25,draw,font=\sffamily\LARGE\bfseries},two/.style={diamond,fill=red!5,draw,font=\sffamily\LARGE\bfseries}]
				\node[one] (u1) at (-3.5,2)  {$u_1$};
				\node[one] (u2) at (-3.5,-2)  {$u_2$};
				\node[two] (w1) at (-2,0)  {$w_1$};
				\node[one] (w2) at (0,0)  {$w_2$};
				\node[two] (w3) at (2,0)  {$w_3$};
				\node[one] (v1) at (3.5,2)  {$v_1$};
				\node[one] (v2) at (3.5,-2)  {$v_2$};
			
				\path[every node/.style={font=\sffamily},line width=1pt]
				(u1) edge  (w1)
				(u2) edge  (w1)
				(w1) edge  (w2)
				(w2) edge  (w3)
				(w3) edge (v1)
				(w3) edge (v2);
						\draw[draw=black] ($(v1.135)+(-0.3,0.3)$) rectangle  ($(v2.-45)+(0.3,-0.3)$);
		\draw[draw=black] ($(u1.135)+(-0.3,0.3)$) rectangle  ($(u2.-45)+(0.3,-0.3)$);
			\node[draw=none,fill=none] (8) at (-3.5,3.5) {\fontsize{20}{22.4}\selectfont$\U{}$};
			\node[draw=none,fill=none] (10) at (3.5,3.5) {\fontsize{20}{22.4}\selectfont$\V{}$};
			\node[draw=none,fill=none] (10) at (0,3.5) {\fontsize{20}{22.4}\selectfont$\W{}$};
				\end{tikzpicture}
			}
			\caption{Figure~\ref{case1b:t}: $\beta_{\text{Ex}2}$.}
			\label{case1b:t}
		\end{subfigure}
		\caption{Two shapes $\alpha_{\text{Ex}2}$ and $\beta_{\text{Ex}2}$ representing $f_{\text{Ex}2}$. As before, nodes with weight $\logn{n}$ are depicted as circles and nodes with weight $\logn{m}$ are depicted as diamonds.}
		\label{case1:t}
	\end{figure}
	
	However, note that each representation gives a different upper bound.
	More specifically, one can check that the set of two vertices $\{u_1,u_2\}$ is a minimum weight separator for  $\alpha_{\text{Ex}2}$, while a single vertex $\{w_2\}$ is a minimum weight separator for $\beta_{\text{Ex}2}$.
	It follows from Corollary~\ref{cor:bern} that $\norm{M_{\alpha_{\text{Ex}2}}}=\tO{mn^{1.5}}$ and $\norm{M_{\beta_{\text{Ex}2}}} = \tO{mn^2}$.
	Hence,  it is preferable to choose $\alpha_{\text{Ex}2}$ as a shape representation of $f_{\text{Ex}2}$ as it provides a sharper norm bound.
	Indeed, one can easily verify that  $\alpha_{\text{Ex}2}$ is the shape that gives the best upper bound. \qed
	
	The main takeaway from Example~\ref{ex2} is that when there are multiple shapes representing a single polynomial, we should choose the shape that gives the best upper bound on $\max_{\norm{x}=1}f(x)$.
	This observation motivates us to consider the following definition:
	\begin{definition}[Optimal shape]
		For a homogeneous polynomial $f$ having a graph matrix representation,
		we say shape $\alpha$ is an \emph{optimal shape} for $f$
		if the graph matrix $M_\alpha$ corresponding to $\alpha$ gives the best possible upper bound (graph matrix norm bound) on $\max_{\norm{x}=1}f(x)$ up to polylogarithmic factors.
	\end{definition}
	
	Having established all of this, we finally consider a general case where $f$ does not directly admit a graph matrix representation.
	As we shall see in the following subsections, this is the case for many applications.
	The general principle for this case is as follows:
	\begin{mdframed}
	\begin{itemize}
	    \item We decompose $f$ into a sum of homogeneous polynomials $f_1+\cdots +f_k$ so that each $f_i$ admits a graph matrix representation.
	\end{itemize}
	\end{mdframed}
With such decomposition, one can obtain  an upper bound on each term $\sup_{\norm{x}=1}f_i(x)$ (for $i=1,2,\dots, k$) using the graph matrix norm bounds. 
	For the rest of this section, we  will recover various upper bound analyses in the literature~\cite{barak2012hypercontractivity,hopkins2016fast,ge2015decomposing} using this strategy.

	We begin with the results from \cite{barak2012hypercontractivity}, namely bounding the $2\to 4$ norm of a random operator.
	The original approach therein relies on some advanced large deviation inequalities based on Orlicz norms, whereas our graph matrix approach is fairly mechanical.
	While recovering the result, we  also illustrate several principles that we will use for other applications in later subsections.

	\subsection{Bound on the $2\to 4$ norm of a random operator}
	\label{fournorm} 
	Our first application is the upper bound on the $2\to 4$ norm of a random operator due to~\cite{barak2012hypercontractivity}. Here there are two dimension parameters $n,m\in \na$.
	Given these parameters, let $A$ be an $m\times n$ matrix with each entry drawn independently from $\Nc(0,1)$.
	Our goal is to prove an upper bound on $\|A\|_{2\to 4}$ defined as $\max_{\left\| x \right\|_2=1} \norm{Ax}_4$.
		\begin{theorem}[{\cite[Theorem 7]{barak2012hypercontractivity}}] \label{baraketal7}
		With high probability, the following upper bound holds: $\|A\|_{2\to 4} \leq 3  m+\tO{ \max\left\{n\sqrt{m},~ n^2\right\}}$.
	\end{theorem}
	\noindent The original proof in~\cite{barak2012hypercontractivity} requires some advanced large deviation inequalities based on Orlicz norms. Here, we take a more direct approach and recover the same bound using graph matrices.

	Let us first represent the $2\to 4$ norm as the maximum of a homogeneous polynomial:
	\begin{align*}
	\max_{\norm{x}_2=1}\norm{Ax}_4^4  = \max_{\norm{x}_2=1}   \sum_{i\in[m]}\sum_{j_1,j_2,j_3,j_4\in [n]} \left[A_{i,j_1}A_{i,j_2}A_{i,j_3}A_{i,j_4}x_{j_1}x_{j_2}x_{j_3}x_{j_4}\right]=:\max_{\norm{x}_2=1}  f(x) \,. 
	\end{align*}
	With this polynomial representation, our analysis will follow the two steps: (i) we  first analyze the dominant terms in the polynomial and then (ii) show that the remaining terms are negligible compared to the dominant terms.
	Let us begin with the analysis of the dominant terms.
	\subsubsection{Analysis of dominant terms}
Generally, the dominant terms are terms where all indices are distinct and the terms which have nonzero expected value:  
\begin{itemize}
        \item ({\bf Distinct indices}:) First, consider the case where all four indices  are distinct.
		One can easily verify that this case can be represented by shape  which consists of vertices $\{u_1,u_2,v_1,v_2,w_1\}$, weights $\w{u_1}=\w{u_2}=\w{v_1}=\w{v_2}= \logn{n}$, $\w{w_1}=\logn{m}$ (here as before, $N:=\max\{n,m\}$),  and edges $\{  \{u_1,w_1\}$,  $\{u_2,w_1\}$, $\{v_1,w_1\}$, $\{v_2,w_1\}  \}$. See Figure~\ref{case1:h}.
	  \begin{figure}[H]
			\centering
			\scalebox{0.5}{
				\begin{tikzpicture}[one/.style={circle,fill=blue!25,draw,font=\sffamily\LARGE\bfseries},two/.style={diamond,fill=red!5,draw,font=\sffamily\LARGE\bfseries}]
				\node[one] (a) at (-2,1)  {$u_1$};
				\node[one] (b) at (-2,-1)  {$u_2$};
				\node[two] (c) at (0,0)  {$w_1$};
				\node[one] (d) at (2,1)  {$v_1$};
				\node[one] (e) at (2,-1)  {$v_2$};
				\node[draw=none,fill=none] (f) at (-2,2) {\fontsize{20}{22.4}\selectfont$U$};
				\node[draw=none,fill=none] (g) at (0,2) {\fontsize{20}{22.4}\selectfont$W$};
				\node[draw=none,fill=none] (h) at (2,2) {\fontsize{20}{22.4}\selectfont$V$};
				\path[every node/.style={font=\sffamily},line width=1pt]
				(a) edge  (c)
				(b) edge  (c)
				(d) edge  (c)
				(e) edge  (c);
				\draw[draw=black] ($(a.135)+(-0.2,0.2)$) rectangle  ($(b.-45)+(0.2,-0.2)$);	
				\draw[draw=black] ($(d.north west)+(-0.2,0.2)$) rectangle ($(e.south east)+(0.2,-0.2)$);
				\end{tikzpicture}
			}
			\caption{A shape representing the case of distinct indices.}
			\label{case1:h}
		\end{figure}
	Depending on the relative sizes between $m$ and $n$, either $\{w_1\}$ or $\{u_1,u_2\}$ is a minimum weight separator.
	Hence, the norm bound for this case reads 
	$\tO{\max\{n^2,~n\sqrt{m} \}}$.
	    \item ({\bf Terms with non-zero expected values}:) Now let us consider the terms which have nonzero expected values. 
	    Note that the terms with non-zero expected values occur only when the shape representation does not contain any edges.
	    Having noticed this, one can easily conclude that the terms with non-zero expected values appear in the following cases: (i) there are two pairs of equal indices or (ii) all four indices are the same.

	     Let us first find  the term with non-zero expected value  for the case (i). 
	     One can notice that the term corresponding to the zero-expected value is represented by the shape consisting of vertices $\{\nn_1,\nn_2,w_1\}$ (recall that we use $\nn$ to denote vertices in the left/right intersection), weights $\w{\nn_1}=\w{\nn_2} =\logn{n}$, $\w{w_1}=\logn{m}$ without edges. See Figure~\ref{case34:h}.
	     	\begin{figure}[H] 
			\centering
			\scalebox{0.5}{
				\begin{tikzpicture}[one/.style={circle,fill=blue!25,draw,font=\sffamily\LARGE\bfseries},two/.style={diamond,fill=red!5,draw,font=\sffamily\LARGE\bfseries}]
				\node[two] (c) at (0,-0.5)  {$w_1$};
				\node[one] (e) at (-1,1.5)  {$\nn_1$}; 
				\node[one] (f) at (1,1.5)  {$\nn_2$}; 
			\node[draw=none,fill=none] (g) at (3,-0.5) {\fontsize{20}{22.4}\selectfont$W$};
		\node[draw=none,fill=none] (h) at (3,1.5) {\fontsize{20}{22.4}\selectfont$U\cap V$};
				\draw[draw=black] ($(e.north west)+(-0.2,0.2)$) rectangle ($(f.south east)+(0.2,-0.2)$);
				\draw[draw=black] ($(e.north west)+(-0.3,0.3)$) rectangle ($(f.south east)+(0.3,-0.3)$);
				\end{tikzpicture}
			}
			\caption{  An optimal shape representing the case when there are two pairs of equal indices.}
			\label{case34:h}
		\end{figure}
		 Since the minimum weight separator is $\{\nn_1,\nn_2\}$, the norm bound reads $\tO{m}$ due to the fact that the middle vertex is isolated.  
		 We will not use this result as our goal is to exactly recover the coefficient appearing in front of the $m$ term.
		 Instead, let us explicitly compute the term:  since there are $\binom{4}{2}/2 = 3$ possible ways of obtaining such equality pattern,  the term for this case is equal to
	    \begin{align} \label{remain:1}
	         3  \sum_{i\in[m]}\sum_{\substack{j_1,j_2 \in [n]\\ j_1,j_2\text{ are distinct}}} \E[A^2_{i,j_1}] \E[A^2_{i,j_2}] x_{j_1}^2x_{j_2}^2=  3  \sum_{i\in[m]}\sum_{\substack{j_1,j_2 \in [n]\\ j_1,j_2\text{ are distinct}}}  x_{j_1}^2x_{j_2}^2\,.
	    \end{align} 
	    
	    Next, let us find the term with non-zero expected value for the case (ii). One can verify that the nonzero expectation term for this case is represented by the shape which consists of vertices $\{\nn_1, w_1\}$, weights $\w{\nn_1} =\logn{n}$, $\w{w_1}=\logn{m}$ without  edges. See Figure~\ref{case52:h}.
	    	\begin{figure}[H] 
			\centering
			\scalebox{0.5}{
				\begin{tikzpicture}[one/.style={circle,fill=blue!25,draw,font=\sffamily\LARGE\bfseries},two/.style={diamond,fill=red!5,draw,font=\sffamily\LARGE\bfseries}]
				\node[two] (c) at (1,-0.5)  {$w_1$};
				\node[one] (e) at (1,1.5)  {$\nn_1$};  
				\node[draw=none,fill=none] (g) at (3,-0.5){\fontsize{20}{22.4}\selectfont$W$};
			\node[draw=none,fill=none] (h) at (3,1.5) {\fontsize{20}{22.4}\selectfont$U\cap V$}; 
					 		\draw[draw=black] ($(e.north west)+(-0.2,0.2)$) rectangle ($(e.south east)+(0.2,-0.2)$);
					 		\draw[draw=black] ($(e.north west)+(-0.3,0.3)$) rectangle ($(e.south east)+(0.3,-0.3)$);
		\end{tikzpicture}
				
			}
			\caption{ An optimal shape representing the case where all four indices are equal.}
			\label{case52:h}
		\end{figure}
	Again the graph matrix norm bound reads $\tO{m}$, and hence, we explicitly represent the term as before:
	     \begin{align} \label{remain:2}
	        \sum_{i\in[m]}\sum_{j_1 \in [n] } \E[A^4_{i,j_1}] x_{j_1}^4 = 3 \sum_{i\in[m]}\sum_{j_1 \in [n] }  x_{j_1}^4\,.
	    \end{align} 
	    Combining  \eqref{remain:1} with \eqref{remain:2}, we obtain the term $ 3 \sum_{i\in[m]} \left(\sum_{j_1}x_{j_1}^2\right)^2  = 3 m$ since $\norm{x}_2=1$. 
	    This exactly recovers the desired $3 m$ in Theorem~\ref{baraketal7}.
	\end{itemize}
Having analyzed the dominant terms, now it suffices to show that the remaining terms are negligible compared to the dominant terms.

	\subsubsection{Analysis of negligible terms}
 Below, we take two different approaches for this. The first approach will be more brute-force in nature: it will consider all possible cases and show---using graph matrices---that each case is dominated by the dominant terms.
 After presenting the brute-force approach, we will present another approach which deals with cases succinctly via some counting arguments. 
 As we move on to other applications where the shape representations become more complicated than Figure~\ref{case1:h}, such counting arguments will be effective in simplifying our arguments.

{\bf 1) First approach: brute-force case studies.}
Let us first list all possible equality patterns between the indices:
	\begin{enumerate} 
	    \item $|\{j_1,j_2,j_3,j_4\}|=4$, i.e., $j_1$, $j_2$, $j_3$, and $j_4$ are all distinct. 
	    \item $|\{j_1,j_2,j_3,j_4\}|=3$, i.e., there is a single pair of indices that are equal.
	    \item $|\{j_1,j_2,j_3,j_4\}|=2$  with each equivalence class of size $2$, i.e., there are two pairs of equal indices.
	    \item $|\{j_1,j_2,j_3,j_4\}|=2$  with one equivalence class of size $3$, i.e., there are three indices that are equal.
	    \item $|\{j_1,j_2,j_3,j_4\}|=1$, i.e., all indices are equal.
	\end{enumerate}
	Note that we already have covered Case 1, and the terms corresponding to the nonzero expectations in Case 3 and 5 are also analyzed.
	We now represent each of the remaining cases by graph matrices and bound each term with our norm bounds: 
		\begin{enumerate}  
	    \item[2:] The second case corresponds to the polynomial
	    \begin{align*}
	        f_{2}(x):= 6\sum_{i\in[m]}\sum_{\substack{j_1,j_2,j_3\in [n]\\ j_1,j_2,j_3\text{ are distinct}}} \left[A_{i,j_1}^2A_{i,j_2}A_{i,j_3}x_{j_1}^2x_{j_2}x_{j_3}\right]\,.
	    \end{align*}
	    To represent each random variable with a member of the orthonormal basis of $\Nc(0,1)$ (see Example~\ref{ex:fourier}), we further decompose this polynomial as follows: $f_2=\sqrt{2}\cdot f_{2\text{-}1}+f_{2\text{-}2}$, where
	    \begin{align*}
	        f_{2\text{-}1}(x)&:= 6\sum_{i\in[m]}\sum_{\substack{j_1,j_2,j_3\in [n]\\ j_1,j_2,j_3\text{ are distinct}}} \left[\frac{A_{i,j_1}^2-1}{\sqrt{2}}A_{i,j_2}A_{i,j_3}x_{j_1}^2x_{j_2}x_{j_3}\right]~~\text{and}\\
	        f_{2\text{-}1}(x)&:= 6\sum_{i\in[m]}\sum_{\substack{j_1,j_2,j_3\in [n]\\ j_1,j_2,j_3\text{ are distinct}}} \left[ A_{i,j_2}A_{i,j_3}x_{j_1}^2x_{j_2}x_{j_3}\right]\,.
	    \end{align*}
	    One can easily verify that $\alpha_{2\text{-}1}$ and $\alpha_{2\text{-}2}$ in Figures~\ref{case21:h} and \ref{case22:h} are optimal shapes for  $f_{2\text{-}1}$ and $f_{2\text{-}2}$, respectively.
	    \begin{figure}[H]
	    \begin{subfigure}[H]{0.5\textwidth}
			\centering
			\scalebox{0.5}{
				\begin{tikzpicture}[one/.style={circle,fill=blue!25,draw,font=\sffamily\LARGE\bfseries},two/.style={diamond,fill=red!5,draw,font=\sffamily\LARGE\bfseries}]
				\node[one] (a) at (-2,1.5)  {$u_1$}; 
				\node[two] (c) at (0,-0.5)  {$w_1$};
				\node[one] (d) at (2,1.5)  {$v_1$}; 
				\node[one] (e) at (0,1.5)  {$\nn_1$}; 
				\node[draw=none,fill=none] (f) at (-2,2.5) {\fontsize{15}{22.4}\selectfont$U\setminus V$};
				\node[draw=none,fill=none] (g) at (-1.5,-0.5) {\fontsize{20}{22.4}\selectfont$W$};
				\node[draw=none,fill=none] (h) at (2,2.5) {\fontsize{15}{22.4}\selectfont$V\setminus U$};
					\node[draw=none,fill=none] (h) at (0,2.5) {\fontsize{15}{22.4}\selectfont$U\cap V$};
				\path[every node/.style={font=\sffamily\LARGE},line width=1pt]
			(a) edge  (c)
			(e) edge node[right]{$2$} (c)
				(d) edge  (c);
				\draw[draw=black] ($(a.135)+(-0.2,0.2)$) rectangle  ($(a.-45)+(0.2,-0.2)$);	
				\draw[draw=black] ($(d.north west)+(-0.2,0.2)$) rectangle ($(d.south east)+(0.2,-0.2)$);
				\draw[draw=black] ($(e.north west)+(-0.2,0.2)$) rectangle ($(e.south east)+(0.2,-0.2)$);
				\draw[draw=black] ($(e.north west)+(-0.3,0.3)$) rectangle ($(e.south east)+(0.3,-0.3)$);
				\end{tikzpicture}
			}
			\caption{Figure \ref{case21:h}: $\alpha_{2\text{-}1}$.}
			\label{case21:h}
		\end{subfigure}
		\begin{subfigure}[H]{0.5\textwidth}
			\centering
			\scalebox{0.5}{
				\begin{tikzpicture}[one/.style={circle,fill=blue!25,draw,font=\sffamily\LARGE\bfseries},two/.style={diamond,fill=red!5,draw,font=\sffamily\LARGE\bfseries}]
				\node[one] (a) at (-2,1.5)  {$u_1$}; 
				\node[two] (c) at (0,-0.5)  {$w_1$};
				\node[one] (d) at (2,1.5)  {$v_1$}; 
				\node[one] (e) at (0,1.5)  {$\nn_1$}; 
				\node[draw=none,fill=none] (f) at (-2,2.5) {\fontsize{15}{22.4}\selectfont$U\setminus V$};
				\node[draw=none,fill=none] (g) at (-1.5,-0.5) {\fontsize{20}{22.4}\selectfont$W$};
				\node[draw=none,fill=none] (h) at (2,2.5) {\fontsize{15}{22.4}\selectfont$V\setminus U$};
					\node[draw=none,fill=none] (h) at (0,2.5) {\fontsize{15}{22.4}\selectfont$U\cap V$};
				\path[every node/.style={font=\sffamily},line width=1pt]
				(a) edge  (c)
				(d) edge  (c);
				\draw[draw=black] ($(a.135)+(-0.2,0.2)$) rectangle  ($(a.-45)+(0.2,-0.2)$);	
				\draw[draw=black] ($(d.north west)+(-0.2,0.2)$) rectangle ($(d.south east)+(0.2,-0.2)$);
				\draw[draw=black] ($(e.north west)+(-0.2,0.2)$) rectangle ($(e.south east)+(0.2,-0.2)$);
				\draw[draw=black] ($(e.north west)+(-0.3,0.3)$) rectangle ($(e.south east)+(0.3,-0.3)$);
				\end{tikzpicture}
			}
			\caption{Figure \ref{case22:h}: $\alpha_{2\text{-}2}$.}
			\label{case22:h}
		\end{subfigure}
			\caption{Shapes for Case 2. We depict nodes with weight $\logn{n}$ as circles and a node with weight $\logn{m}$ as a diamond. Again, the edges are of label $1$ unless specified otherwise.}
	    \end{figure}
	    Notice the edge of label $2$ in $\alpha_{2\text{-}1}$ due to appearance of the random variable $(A_{i,j_1}^2-1)/\sqrt{2}$.
	    For both  $\alpha_{2\text{-}1}$ and $\alpha_{2\text{-}2}$, depending on the relative sizes between $m$ and $n$, either $\{\nn_1,u_1\}$ or $\{\nn_1,w_1\}$  is a minimum weight separator.
	Hence, the norm bound for both cases reads:  $\tO{\max\{n,~\sqrt{nm} \}}$, which is dominated by $\tO{\max\{n^2,~n\sqrt{m} \}}$.
	    \item[3:] The third case corresponds to the polynomial (since the nonzero expectation term is already handled):
	    \begin{align*}
	        f_3(x):= 3\sum_{i\in[m]}\sum_{\substack{j_1,j_2 \in [n]\\ j_1,j_2\text{ are distinct}}} \left[A_{i,j_1}^2A_{i,j_2}^2 x_{j_1}^2x_{j_2}^2 -x_{j_1}^2x_{j_2}^2\right]\,.
	    \end{align*}
	    Again, we further decompose the polynomial so that each edge variable is a member of the orthonormal basis of $\Nc(0,1)$: $f_3(x):= 2\cdot f_{3\text{-}1} + \sqrt{2}\cdot f_{3\text{-}2}+\sqrt{2}\cdot f_{3\text{-}3}$,
	    where
	    \begin{align*}
	         f_{3\text{-}1} (x)&:= 3\sum_{i\in[m]}\sum_{\substack{j_1,j_2 \in [n]\\ j_1,j_2\text{ are distinct}}} \left[\frac{A_{i,j_1}^2-1}{\sqrt{2}}\frac{A_{i,j_2}^2-1}{\sqrt{2}} x_{j_1}^2x_{j_2}^2\right]\\
	          f_{3\text{-}2} (x)&:= 3\sum_{i\in[m]}\sum_{\substack{j_1,j_2 \in [n]\\ j_1,j_2\text{ are distinct}}} \left[\frac{A_{i,j_1}^2-1}{\sqrt{2}}A_{i,j_2}^2 x_{j_1}^2x_{j_2}^2\right]\\
	           f_{3\text{-}3} (x)&:= 3\sum_{i\in[m]}\sum_{\substack{j_1,j_2 \in [n]\\ j_1,j_2\text{ are distinct}}} \left[A_{i,j_1}^2\frac{A_{i,j_2}^2-1}{\sqrt{2}}x_{j_1}^2x_{j_2}^2\right]\,.
	    \end{align*}
	    One can easily check that the shapes $\alpha_{3\text{-}1}$, $\alpha_{3\text{-}2}$, and $\alpha_{3\text{-}3}$ in Figures~\ref{case31:h}, \ref{case32:h}, and \ref{case33:h} represent $ f_{3\text{-}1}$, $ f_{3\text{-}2} $, and $f_{3\text{-}3}$, respectively. 
	    \begin{figure}[H]
	        \centering
	        	\begin{subfigure}[H]{0.24\textwidth}
			\centering
			\scalebox{0.5}{
				\begin{tikzpicture}[one/.style={circle,fill=blue!25,draw,font=\sffamily\LARGE\bfseries},two/.style={diamond,fill=red!5,draw,font=\sffamily\LARGE\bfseries}]
				\node[two] (c) at (0,-0.5   )  {$w_1$};
				\node[one] (e) at (-1,1.5)  {$\nn_1$}; 
				\node[one] (f) at (1,1.5)  {$\nn_2$}; 
				\node[draw=none,fill=none] (ff) at (-2,2.5) {};
			\node[draw=none,fill=none] (g) at (3,-0.5) {\fontsize{20}{22.4}\selectfont$W$};
		\node[draw=none,fill=none] (h) at (3,1.5) {\fontsize{20}{22.4}\selectfont$U\cap V$};
				\path[every node/.style={font=\sffamily\LARGE},line width=1pt]
				(e) edge node[right]{$2$}  (c) 
				(f) edge node[right]{$2$}  (c); 
				\draw[draw=black] ($(e.north west)+(-0.2,0.2)$) rectangle ($(f.south east)+(0.2,-0.2)$);
				\draw[draw=black] ($(e.north west)+(-0.3,0.3)$) rectangle ($(f.south east)+(0.3,-0.3)$);
				\end{tikzpicture}
			}
			\caption{Figure \ref{case31:h}: $\alpha_{3\text{-}1}$.}
			\label{case31:h}
		\end{subfigure}  
		\begin{subfigure}[H]{0.24\textwidth}
			\centering
			\scalebox{0.5}{
				\begin{tikzpicture}[one/.style={circle,fill=blue!25,draw,font=\sffamily\LARGE\bfseries},two/.style={diamond,fill=red!5,draw,font=\sffamily\LARGE\bfseries}]
				\node[two] (c) at (0,-0.5)  {$w_1$};
				\node[one] (e) at (-1,1.5)  {$\nn_1$}; 
				\node[one] (f) at (1,1.5)  {$\nn_2$}; 
			\node[draw=none,fill=none] (g) at (3,-0.5) {\fontsize{20}{22.4}\selectfont$W$};
		\node[draw=none,fill=none] (h) at (3,1.5) {\fontsize{20}{22.4}\selectfont$U\cap V$};
				\path[every node/.style={font=\sffamily\LARGE},line width=1pt]
				(e) edge node[right]{$2$}  (c) ; 
				\draw[draw=black] ($(e.north west)+(-0.2,0.2)$) rectangle ($(f.south east)+(0.2,-0.2)$);
				\draw[draw=black] ($(e.north west)+(-0.3,0.3)$) rectangle ($(f.south east)+(0.3,-0.3)$);
				\end{tikzpicture}
			}
			\caption{Figure \ref{case32:h}: $\alpha_{3\text{-}2}$.}
			\label{case32:h}
		\end{subfigure}
		\begin{subfigure}[H]{0.24\textwidth}
			\centering
			\scalebox{0.5}{
				\begin{tikzpicture}[one/.style={circle,fill=blue!25,draw,font=\sffamily\LARGE\bfseries},two/.style={diamond,fill=red!5,draw,font=\sffamily\LARGE\bfseries}]
				\node[two] (c) at (0,-0.5)  {$w_1$};
				\node[one] (e) at (-1,1.5)  {$\nn_1$}; 
				\node[one] (f) at (1,1.5)  {$\nn_2$}; 
				\node[draw=none,fill=none] (g) at (3,-0.5) {\fontsize{20}{22.4}\selectfont$W$};
		\node[draw=none,fill=none] (h) at (3,1.5) {\fontsize{20}{22.4}\selectfont$U\cap V$};
				\path[every node/.style={font=\sffamily\LARGE},line width=1pt]
				(f) edge node[right]{$2$}  (c); 
				\draw[draw=black] ($(e.north west)+(-0.2,0.2)$) rectangle ($(f.south east)+(0.2,-0.2)$);
				\draw[draw=black] ($(e.north west)+(-0.3,0.3)$) rectangle ($(f.south east)+(0.3,-0.3)$);
				\end{tikzpicture}
			}
			\caption{Figure \ref{case33:h}: $\alpha_{3\text{-}3}$.}
			\label{case33:h}
		\end{subfigure}
	        \caption{Shapes for Case 3. The edges have label $1$ unless specified otherwise.} 
	    \end{figure}
	    For  $\alpha_{3\text{-}1}$, $\alpha_{3\text{-}2}$, and  $\alpha_{3\text{-}3}$, $\{\nn_1,\nn_2\}$ is a minimum weight separator, and hence, the norm bound reads $\tO{\sqrt{m}}$, which   is clearly dominated by $\tO{\max\{n^2,~n\sqrt{m} \}}$.
	  
	    \item[4:] The fourth case corresponds to the polynomial
	     \begin{align*}
	        f_4(x):= 4\sum_{i\in[m]}\sum_{\substack{j_1,j_2 \in [n]\\ j_1,j_2\text{ are distinct}}} \left[A_{i,j_1}^3A_{i,j_2} x_{j_1}^3x_{j_2}\right]\,.
	    \end{align*}
	    To keep the edge variables within the orthonormal basis, we decompose $A_{i,j_1}^3= \sqrt{6}\cdot \frac{A_{i,j_1}^3-3A_{i,j_1}}{\sqrt{6}} + 3A_{i,j_1}$.
	    Note that both terms in the decomposition can be represented by the shape consisting of three vertices $u_1,w_1, \nn_1$ with edges $\{ \{u_1,w_1\}, \{\nn_1,w_1\}\}$.
	    The only distinction between the two cases is in whether edge $\{\nn_1,w_1\}$ gets label $3$ or $1$. 
	    See  Figure~\ref{case4:h}.
	    \begin{figure}[H]
			\centering
			\scalebox{0.5}{
				\begin{tikzpicture}[one/.style={circle,fill=blue!25,draw,font=\sffamily\LARGE\bfseries},two/.style={diamond,fill=red!5,draw,font=\sffamily\LARGE\bfseries}]
				\node[two] (c) at (1,-0.5)  {$w_1$};
				\node[one] (e) at (1,1.5)  {$\nn_1$}; 
				\node[one] (f) at (-1,-0.5)  {$u_1$}; 
				\node[draw=none,fill=none] (g) at (3,-0.5){\fontsize{20}{22.4}\selectfont$W$};
			\node[draw=none,fill=none] (h) at (3,1.5) {\fontsize{20}{22.4}\selectfont$U\cap V$};
			\node[draw=none,fill=none] (i) at (-1,0.5) {\fontsize{20}{22.4}\selectfont$U\setminus V$};
						\path[every node/.style={font=\sffamily\LARGE},line width=1pt]
				(f) edge  (c)
				(e) edge node[right]{$3$ or $1$}  (c); 
				\draw[draw=black] ($(e.north west)+(-0.2,0.2)$) rectangle ($(e.south east)+(0.2,-0.2)$);
				\draw[draw=black] ($(e.north west)+(-0.3,0.3)$) rectangle ($(e.south east)+(0.3,-0.3)$);
				\draw[draw=black] ($(f.north west)+(-0.2,0.2)$) rectangle ($(f.south east)+(0.2,-0.2)$);
				\end{tikzpicture}
			}
			\caption{ Shapes for Case 4. Here  the only distinction between the two subcases is whether edge $\{\nn_1,w_1\} $ gets label $3$ or $1$. }
			\label{case4:h}
		\end{figure}
	    For both cases, since $\{\nn_1\}$ is the minimum weight separator, the norm bound reads $\tO{\sqrt{nm}}$, which   is  dominated by $\tO{\max\{n^2,~n\sqrt{m} \}}$.

	    \item[5:] Lastly, the fifth case corresponds to the polynomial
	    	     \begin{align*}
	        f_5(x):= \sum_{i\in[m]}\sum_{j_1 \in [n] } \left[\left(A_{i,j_1}^4-3\right)  x_{j_1}^4\right]\,.
	    \end{align*}
	    Again we decompose the random variable $A_{i,j_1}^4$ as 
	    $$A_{i,j_1}^4 -3 = \sqrt{4!}\cdot\frac{A_{i,j_1}^4-6A_{i,j_1}+3}{\sqrt{4!}}+6\sqrt{2}\cdot \frac{A_{i,j_1}^2-1}{\sqrt{2}} $$
	    so that each edge variable becomes a member of the orthonormal basis.
	    Then, following similar reasoning to the previous cases, it turns out each term in the decomposition can be represented by the shapes in Figure~\ref{case51:h}.
	    \begin{figure}[H]
		\centering
		
			\scalebox{0.5}{
				\begin{tikzpicture}[one/.style={circle,fill=blue!25,draw,font=\sffamily\LARGE\bfseries},two/.style={diamond,fill=red!5,draw,font=\sffamily\LARGE\bfseries}]
				\node[two] (c) at (1,-0.5)  {$w_1$};
				\node[one] (e) at (1,1.5)  {$\nn_1$};  
				\node[draw=none,fill=none] (g) at (3,-0.5){\fontsize{20}{22.4}\selectfont$W$};
			\node[draw=none,fill=none] (h) at (3,1.5) {\fontsize{20}{22.4}\selectfont$U\cap V$}; 
						\path[every node/.style={font=\sffamily\LARGE},line width=1pt]
				(c) edge node[right]{$4$ or $2$}  (e); 
				\draw[draw=black] ($(e.north west)+(-0.2,0.2)$) rectangle ($(e.south east)+(0.2,-0.2)$);
				\draw[draw=black] ($(e.north west)+(-0.3,0.3)$) rectangle ($(e.south east)+(0.3,-0.3)$);
		\end{tikzpicture}
				
			}
		 
		\caption{Shapes for Case 4. Here the only distinction between the two subcases is  whether the edge $\{\nn_1,w_1\} $ gets label $4$ or $2$. }
			\label{case51:h}
	 
	\end{figure}
	    For those shapes, $\{\nn_1\}$ is the minimum weight separator, and hence, the norm bound is $\tO{\sqrt{m}}$, which is dominated by $\tO{\max\{n^2,~n\sqrt{m} \}}$.

	\end{enumerate}
	To summarize, since all the remaining cases are dominated by the bound $\tO{\max\{n^2,~n\sqrt{m} \}}$ for the dominant term, Theorem~\ref{baraketal7}  follows. \qed

{\bf 2) Second approach: case reduction via counting arguments.}
	As readers may have noticed, checking cases can become quite laborious, especially when we consider larger shapes.
	Here we offer a simpler approach which can significantly reduce the number of cases to check.
	We note that the counting arguments that we will develop here will be used throughout the rest of this section to simplify cases.

		Our strategy is to show by counting arguments that the remaining cases are bounded by the bound for the dominant case. First, note that for the remaining cases, the middle vertex is not isolated since they have non-zero expected values and thus have at least one edge incident to the middle vertex.
		Moreover, since the case where all four indices are distinct has been taken care of,  there exists a pair of indices that are equal; without loss of generality, let those indices be $j_1$ and $j_2$.
	Now, let us place the vertex corresponding to $j_1=j_2$ in the intersection $\U{}\cap \V{}$ in the shape representation. Then since the vertex is in the intersection, any vertex separator has to contain the vertex, and hence, the norm bound is at most $\tO{n\sqrt{m}}$ since we have assumed the middle vertex is not isolated and there are at most two other vertices of weight $\logn{n}$ and at most one vertex of weight $\logn{m}$.
	This is clearly dominated by the norm bound 	$\tO{\max\{n^2,~n\sqrt{m} \}}$ for the case where all four indices are distinct. \qed

	\subsection{Overcomplete tensor decomposition algorithm analysis} \label{overtensor}
	In this section, we will recover the main component of the overcomplete tensor decomposition analysis in~\cite{ge2015decomposing}.
	The proof therein relies on some clever arguments: it makes use of a clever decoupling technique due to~\cite{de1995decoupling} together with some incoherence properties of random vectors. 
	In contrast, we again emphasize that our proof based on graph matrices is  \emph{mechanical}.

To describe the setting, it was demonstrated in~\cite{ge2015decomposing} that a key component of analyzing an overcomplete tensor decomposition algorithm is showing an upper bound on the injective tensor norm of random tensors.
	Formally, for integers $n,m$ such that $n\leq m \leq c\cdot  n^{1.5}$ for some constant $c>0$, let $a_1, \dots, a_m\in \re^n$ be random vectors whose entries are i.i.d. Rademacher random variables.
	Consider an $n\times n\times n$ tensor $\T=\sum_{i=1}^m a_i^{\otimes 3}$.
	Then, our  goal is to prove an upper bound  on the injective tensor norm defined as  $\norm{\T}_{\text{inj}}:=\sup_{\norm{x}=1} \inp{\T}{x^{\otimes 3}}$.
In this section, we will recover their upper bound of $n^3+\tO{n^{1.5}m}$.
	
	\begin{theorem}[{\cite[Theorem 4.2]{ge2015decomposing}}] \label{gema4.5}
		Suppose that $n\leq m \leq c\cdot  n^{1.5}$ for some constant $c>0$. With high probability, we have $\norm{\T}_{\emph{inj}} \leq n^{1.5}+ \tO{m}$.
	\end{theorem}
 To transform the tensor into a matrix, we use Cauchy-Schwarz inequality following  \cite{ge2015decomposing}:
	$$ (\inp{\T}{x^{\otimes 3}})^2 = \inp{\sum_{i=1}^m \inp{a_i}{x}^2 a_i}{x}^2 \leq \norm{\sum_{i=1}^m \inp{a_i}{x}^2 a_i}^2 \norm{x}^2 =\norm{\sum_{i=1}^m \inp{a_i}{x}^2 a_i}^2 \,. $$ 
	Hence, it suffices to show $\norm{\sum_{i=1}^m \inp{a_i}{x}^2 a_i}^2\leq n^3 +\tO{mn^{1.5}}$. 
	Expanding the term $\norm{\sum_{i=1}^m \inp{a_i}{x}^2 a_i}^2$ together with the fact that  $\inp{a_i}{a_i} =n$, we obtain 
	\begin{align*}
	\left(\inp{T}{x^{\otimes 3}}\right)^2 \leq  n \cdot  \underbrace{\sum_{i=1}^m \inp{a_i}{x}^4}_{=:g(x)} +  \underbrace{\sum_{i_1\neq i_2} \inp{a_{i_1}}{a_{i_2}} \inp{a_{i_1}}{x}^2 \inp{a_{i_2}}{x}^2}_{=:h(x)}\,.
	\end{align*} 
	Our agenda is to prove that $\max_{\norm{x}=1}g(x)=n^2+\tO{m\sqrt{n}}$ and $\max_{\norm{x}=1}h(x)=\tO{mn^{1.5}}$.	
	\subsubsection{Showing \texorpdfstring{$\max_{\norm{x}=1}h(x)=\tO{mn^{1.5}}$}{+1} using graph matrices} \label{hh}
	First, let us write $h(x_1,\dots,x_n)$ as a polynomial (again, we will simply use $a_{i,j}$ to denote the $j$-th coordinate of the vector $a_i$):
	\begin{align*}
	\sum_{\substack{j_1,j_2,j_3,j_4\in[n]}} \left[ 	x_{j_1}x_{j_2}x_{j_3}x_{j_4}\sum_{\substack{i_1\neq i_2\in [m],~ j_5\in[n]    } } a_{i_1,j_5} a_{i_2,j_5} a_{i_1,j_1} a_{i_1,j_2} a_{i_2,j_3} a_{i_2,j_4}
	\right]\,.
	\end{align*} 
Again,  let us  first analyze the dominant terms:
	\begin{itemize}
	    \item ({\bf Distinct indices}:) Note that  the polynomial for the distinct indices corresponds to  Example~\ref{ex2} wherein we have shown the upper bound $\tO{mn^{1.5}}$.
	    Here we recall the illustration of the optimal shape (Figure~\ref{case1a:t}) for the reader's convenience.
	    \begin{figure}[H]
			\centering
			\scalebox{0.5}{
				\begin{tikzpicture}[one/.style={circle,fill=blue!25,draw,font=\sffamily\LARGE\bfseries},two/.style={diamond,fill=red!5,draw,font=\sffamily\LARGE\bfseries}]
				\node[one] (u1) at (-3,2)  {$u_1$};
				\node[one] (u2) at (-3,-2)  {$u_2$};
				\node[two] (w1) at (0,2)  {$w_1$};
				\node[one] (w2) at (0,0)  {$w_2$};
				\node[two] (w3) at (0,-2)  {$w_3$};
				\node[one] (v1) at (3,2)  {$v_1$};
				\node[one] (v2) at (3,-2)  {$v_2$};
				\path[every node/.style={font=\sffamily},line width=1pt]
				(u1) edge  (w1)
				(u2) edge  (w3)
				(w1) edge  (w2)
				(w2) edge  (w3)
				(w1) edge (v1)
				(w3) edge (v2);
						\draw[draw=black] ($(v1.135)+(-0.3,0.3)$) rectangle  ($(v2.-45)+(0.3,-0.3)$);
		\draw[draw=black] ($(u1.135)+(-0.3,0.3)$) rectangle  ($(u2.-45)+(0.3,-0.3)$);
			\node[draw=none,fill=none] (8) at (-3,3.5) {\fontsize{20}{22.4}\selectfont$\U{}$};
			\node[draw=none,fill=none] (10) at (3,3.5) {\fontsize{20}{22.4}\selectfont$\V{}$};
			\node[draw=none,fill=none] (10) at (0,3.5) {\fontsize{20}{22.4}\selectfont$\W{}$};
				\end{tikzpicture}
				
			}
			\caption{An optimal shape representing the case of distinct indices.}
			\label{recall2}
		\end{figure}
	    
	    \item ({\bf Terms with non-zero expected values}:) 
	    In order to have terms with non-zero expected values, note that the shape representation should not contain any edges.
	    In particular, the vertices corresponding to $i$-indices, namely $i_1$ and $i_2$, need to be isolated.
	    
	    	However, the vertices corresponding  $i_1$ and $i_2$  cannot be isolated.
	To see this, first notice that there are three random variables $a_{i_1,j_1}$, $a_{i_1,j_2}$, and $a_{i_1,j_5}$ which has index $i_1$.
	Hence, no matter how $j_1$, $j_2$, and $j_5$ are equal to each other, there will be still an edge incident to the vertex representing $i_1$.
	In particular, even when $j_1=j_2=j_5$, there is a random variable $a_{i_1,j_1}^3=a_{i_1,j_1}$ ($\because$ $a_{i_1,j_1}$ is a Rademacher random variable), which implies that the vertex corresponding to $i_1$ is still incident to an edge.
	The same logic applies to $i_2$. This implies that there are no terms with non-zero expected values.
	
	Actually, one can argue a stronger statement that any shape representation cannot have isolated middle vertices. We develop it here for later use. 
	To prove this, it is enough to demonstrate for the case where $j_5$ is not equal to other $j$-indices that the vertex representing $j_5$ cannot be isolated.
	This statement is indeed true because there are two zero-mean random variables $a_{i_1,j_5}$ and  $a_{i_2,j_5}$, which are distinct since $i_1\neq i_2$; hence, there are two edges incident to the vertex representation of $j_5$.
	\end{itemize}
Now, it suffices to show that the remaining cases are dominated by the dominant case above.
Again, rather than listing all possible equality patterns between indices and analyze them one by one, we follow the previous sections and devise some counting arguments.

As we have demonstrated above, any shape representation does not have isolated middle vertices. Consequently, the upper bound on the maximum value due to our main theorem reads
	$$ \tO{mn^{\frac{1}{2}\#\text{ vertices of weight $\logn{n}$ outside of a weight minimum vertex separator}}  } $$
	since there are two vertices of weight $\logn{m}$ (corresponding to $i_1$ and $i_2$). Below, we will demonstrate that for all cases, there are at most three vertices of weight $\logn{n}$ outside of a minimum weight separator, yielding the desired upper bound of $\tO{mn^{1.5}}$.
	\begin{enumerate}
	    \item First, assume that there is more than one equality between the $j$-indices, i.e.,  $j_1,j_2,j_3,j_4,j_5$. 	In this case, the number of vertices of weight $\logn{n}$ becomes at most $3$.
	Hence, we may assume that there is only one equality between $j_1,j_2,j_3,j_4,j_5$.
	If the equality does not involve $j_5$, i.e., it is between $j_1,j_2,j_3,j_4$, one can place the vertex corresponding to the equal indices in the intersection $\U{}\cap \V{}$.
	Then, that vertex has to be in any  vertex separator, which implies that there are at most three $\logn{n}$-weighted vertices outside of any vertex separator.
	Thus, if the only equality does not involve $j_5$, we get the desired bound.
	
	\item Hence, it suffices to consider the case where the equality includes index $j_5$. Without loss of generality, assume that the equality is $j_1=j_5$.
	Let us place the two vertices corresponding to $j_3$ and $j_4$ in $\U{}$ and $\V{}$, respectively.
	Then, with such a shape representation,  there must be a path between the two vertices as there are  random variables $a_{i_2,j_3}$ and $ a_{i_2,j_4}$.
	From this, one can deduce that between the vertices corresponding to $j_3$ and $j_4$, one of them needs to be in a vertex separator.
	Hence, any minimum weight separator has to contain at least one $\logn{n}$-weight vertex and since another  $\logn{n}$-weight  vertex is removed by the equality, we conclude that there are at most three vertices of weight $\logn{n}$ outside of any vertex separator.
	Therefore, we get the desired bound for this case as well.\qed
 \end{enumerate}
	
	\subsubsection{Showing \texorpdfstring{$\max_{\norm{x}=1}g(x)=n^2+\tO{m\sqrt{n}}$}{+1} using graph matrices}

	Recall $g(x):=\sum_{i=1}^m \inp{a_i}{x}^4$. Once again, following \cite{ge2015decomposing}, we apply Cauchy-Schwarz on $g^2$:	$$ \left( \sum_{i=1}^m \inp{a_i}{x}^4 \right)^2 = \inp{\sum_{i=1}^m \inp{a_i}{x}^3 a_i}{x}^2 \leq \norm{\sum_{i=1}^m \inp{a_i}{x}^3 a_i}^2\,. $$
	Hence, it suffices to show $ \norm{\sum_{i=1}^m \inp{a_i}{x}^3 a_i}^2 \leq n^4+\tO{mn^{2.5}}$.
	Expanding the last term together with the the fact that $\inp{a_i}{a_i}=n$, we obtain 
	\begin{align*}
	    g^2(x) \leq n\cdot \underbrace{\sum_{i=1}^m \inp{a_i}{x}^6 }_{=:g_1(x)} + \underbrace{ \sum_{i\neq j} \inp{a_i}{a_j} \inp{a_i}{x}^3 \inp{a_j}{x}^3}_{=:g_2(x)}\,,
	\end{align*}  
	We will show that $\max_{\norm{x}=1}g_1(x) =n^3+\tO{mn^{1.5}}$ and $\max_{\norm{x}=1}g_2(x) =\tO{mn^{2.5}}$.
	
	{\bf 1) Proof for $g_2$:}
	Actually, our graph matrix analysis yields a better bound than what we need to show:
	we will show that $\tO{\sqrt{m}n^3}$ under a relaxed condition $m
	\leq \OO{n^2}$ rather than $m 
	\leq \OO{n^{1.5}}$;
	here the bound $\tO{\sqrt{m}n^3}$  is certainly better than $\tO{mn^{2.5}}$ 
	since $n\leq m$.

	Let us first write $g_2$ as a sum of monomials:
	\begin{align*}
	\sum_{\substack{j_1,j_2,j_3,j_4,j_5,j_6\in[n]}} \left[ 	x_{j_1}x_{j_2}x_{j_3}x_{j_4}x_{j_5}x_{j_6} \sum_{\substack{i_1\neq i_2\in [m],~ j_7\in[n]    } } a_{i_1,j_7} a_{i_2,j_7} a_{i_1,j_1} a_{i_1,j_2} a_{i_1,j_3} a_{i_2,j_4}
	a_{i_2,j_5}
	a_{i_2,j_6}
	\right]\,.
	\end{align*}
We begin with the dominant terms:
	\begin{itemize}
	    \item ({\bf Distinct indices}:) To begin with, we consider the case where $j_1,j_2,\dots,j_7$ are all distinct.
	Similarly to Example~\ref{ex2}, one can check that the  shape  in Figure~\ref{g21} is an optimal shape representing this case.
		\begin{figure}[H]
			\centering
			\scalebox{0.5}{
				\begin{tikzpicture}[one/.style={circle,fill=blue!25,draw,font=\sffamily\LARGE\bfseries},two/.style={diamond,fill=red!5,draw,font=\sffamily\LARGE\bfseries}, three/.style={circle split,fill=blue!25,draw,font=\sffamily\LARGE\bfseries}]
				\node[one] (u1) at (-2,2)  {$u_1$};
				\node[one] (u2) at (-2,0) {$u_2$};
				\node[one] (u3) at (-2,-2) {$u_3$};
				\node[two] (w1) at (0,2)  {$w_1$};
				\node[one] (w2) at (0,0)  {$w_2$};
				\node[two] (w3) at (0,-2)  {$w_3$};
				\node[one] (v1) at (2,2)  {$v_1$};
				\node[one] (v2) at (2,0)  {$v_2$};
				\node[one] (v3) at (2,-2)  {$v_3$};
				\path[every node/.style={font=\sffamily},line width=1pt]
				(u1) edge  (w1)
				(u2) edge (w1)
				(u3) edge (w3)
				(w1)  edge  (w2)
				(w2) edge  (w3)
				(w1) edge (v1)
				(w3) edge (v2)
				(w3) edge (v3);
					\draw[draw=black] ($(v1.135)+(-0.3,0.3)$) rectangle  ($(v3.-45)+(0.3,-0.3)$);
		\draw[draw=black] ($(u1.135)+(-0.3,0.3)$) rectangle  ($(u3.-45)+(0.3,-0.3)$);
			\node[draw=none,fill=none] (8) at (-2,3.5) {\fontsize{20}{22.4}\selectfont$\U{}$};
			\node[draw=none,fill=none] (10) at (2,3.5) {\fontsize{20}{22.4}\selectfont$\V{}$};
			\node[draw=none,fill=none] (10) at (0,3.5) {\fontsize{20}{22.4}\selectfont$\W{}$};
				\end{tikzpicture}}
			\caption{An optimal shape for the case of  distinct indices.}
			\label{g21}
		\end{figure}
 \noindent Since $\{w_1,u_3\}$ is a minimum weight separator, we obtain the bound of  $\tO{\sqrt{m}n^3}$.
 \item ({\bf Terms with non-zero expected values}:) 
In order to obtain terms with non-zero expected values, the shape representation should not  contain any edges.
In particular, the vertices corresponding to $i$-indices, namely $i_1$ and $i_2$ have to be isolated.
It turns out from this, one can characterize all possible shapes.

To begin with, to make the vertex corresponding to $i_1$ isolated, the four random variables that depend on index $i_1$, i.e. $a_{i_1,j_1}$, $a_{i_1,j_2}$, $a_{i_1, j_3}$, and $a_{i_1,j_7}$, need to be either (i)  paired up into two equal variables or (ii) all equal.
The same two cases  apply to  the four random variables that depend on index $i_2$, i.e. $a_{i_2,j_4}$, $a_{i_2,j_5}$, $a_{i_2, j_6}$, and $a_{i_2,j_7}$.
Hence, there are   four different cases in total, and one can show that each case can be represented by one of the shapes in Figure~\ref{fig:iso}:
	
	\begin{figure}[H]
		\centering
		\begin{subfigure}[H]{0.24\textwidth}
			\centering
			\scalebox{0.5}{
		\begin{tikzpicture}[one/.style={circle,fill=blue!25,draw,font=\sffamily\LARGE\bfseries},two/.style={diamond,fill=red!5,draw,font=\sffamily\LARGE\bfseries}, three/.style={circle split,fill=blue!25,draw,font=\sffamily\LARGE\bfseries}]
			\node[one] (u1) at (-2,2)  {$\nn_1$};
				\node[one] (u2) at (-2,0) {$\nn_2$};
				\node[one] (u3) at (-2,-2) {$\nn_3$};
				\node[two] (w1) at (1,2)  {$w_1$};
				\node[two] (w3) at (1,-2)  {$w_3$};
		 
		\draw[draw=black] ($(u1.135)+(-0.3,0.3)$) rectangle  ($(u3.-45)+(0.3,-0.3)$);
		\draw[draw=black] ($(u1.135)+(-0.4,0.4)$) rectangle  ($(u3.-45)+(0.4,-0.4)$);
			\node[draw=none,fill=none] (8) at (-2,3.5) {\fontsize{20}{22.4}\selectfont$\U{}\cap \V{}$};
			 
			\node[draw=none,fill=none] (10) at (1,3.5) {\fontsize{20}{22.4}\selectfont$\W{}$};
				\end{tikzpicture}}
	\end{subfigure}
	\begin{subfigure}[H]{0.24\textwidth}
			\centering
			\scalebox{0.5}{
		\begin{tikzpicture}[one/.style={circle,fill=blue!25,draw,font=\sffamily\LARGE\bfseries},two/.style={diamond,fill=red!5,draw,font=\sffamily\LARGE\bfseries}, three/.style={circle split,fill=blue!25,draw,font=\sffamily\LARGE\bfseries}]
			\node[one] (u2) at (-2,0) {$\nn_1$};
				\node[one] (u3) at (-2,-2) {$\nn_2$};
				\node[two] (w1) at (1,2)  {$w_1$};
				\node[two] (w3) at (1,-2)  {$w_3$};
		 
		\draw[draw=black] ($(u2.135)+(-0.3,0.3)$) rectangle  ($(u3.-45)+(0.3,-0.3)$);
		\draw[draw=black] ($(u2.135)+(-0.4,0.4)$) rectangle  ($(u3.-45)+(0.4,-0.4)$);
			\node[draw=none,fill=none] (8) at (-2,3.5) {\fontsize{20}{22.4}\selectfont$\U{}\cap \V{}$};
			 
			\node[draw=none,fill=none] (10) at (1,3.5) {\fontsize{20}{22.4}\selectfont$\W{}$};
				\end{tikzpicture}}
	\end{subfigure}
	\begin{subfigure}[H]{0.24\textwidth}
			\centering
			\scalebox{0.5}{
		\begin{tikzpicture}[one/.style={circle,fill=blue!25,draw,font=\sffamily\LARGE\bfseries},two/.style={diamond,fill=red!5,draw,font=\sffamily\LARGE\bfseries}, three/.style={circle split,fill=blue!25,draw,font=\sffamily\LARGE\bfseries}]
			\node[one] (u1) at (-2,2)  {$\nn_1$};
				\node[one] (u2) at (-2,0) {$\nn_2$};
		 		\node[two] (w1) at (1,2)  {$w_1$};
				\node[two] (w3) at (1,-2)  {$w_3$};
		 
		\draw[draw=black] ($(u1.135)+(-0.3,0.3)$) rectangle  ($(u2.-45)+(0.3,-0.3)$);
		\draw[draw=black] ($(u1.135)+(-0.4,0.4)$) rectangle  ($(u2.-45)+(0.4,-0.4)$);
			\node[draw=none,fill=none] (8) at (-2,3.5) {\fontsize{20}{22.4}\selectfont$\U{}\cap \V{}$};
			 
			\node[draw=none,fill=none] (10) at (1,3.5) {\fontsize{20}{22.4}\selectfont$\W{}$};
				\end{tikzpicture}}
	\end{subfigure}
	\begin{subfigure}[H]{0.24\textwidth}
			\centering
			\scalebox{0.5}{
		\begin{tikzpicture}[one/.style={circle,fill=blue!25,draw,font=\sffamily\LARGE\bfseries},two/.style={diamond,fill=red!5,draw,font=\sffamily\LARGE\bfseries}, three/.style={circle split,fill=blue!25,draw,font=\sffamily\LARGE\bfseries}]
					\node[one] (u2) at (-2,0) {$\nn_1$};
			 	\node[two] (w1) at (1,2)  {$w_1$};
				\node[two] (w3) at (1,-2)  {$w_3$};
		 
		\draw[draw=black] ($(u2.135)+(-0.3,0.3)$) rectangle  ($(u2.-45)+(0.3,-0.3)$);
		\draw[draw=black] ($(u2.135)+(-0.4,0.4)$) rectangle  ($(u2.-45)+(0.4,-0.4)$);
			\node[draw=none,fill=none] (8) at (-2,3.5) {\fontsize{20}{22.4}\selectfont$\U{}\cap \V{}$};
			 
			\node[draw=none,fill=none] (10) at (1,3.5) {\fontsize{20}{22.4}\selectfont$\W{}$};
				\end{tikzpicture}}
	\end{subfigure}	\caption{Four shapes  that corresponds to the polynomials with non-zero expected values.}
		\label{fig:iso}
	\end{figure}
	   In all four different cases, one can easily check that we have the norm bounds $\tO{m^2}$, which is certainly $\tO{\sqrt{m}n^3}$ since $m=\OO{n^2}$.
	\end{itemize}
Having analyzed the dominant terms, following the previous sections, we show that the remaining terms are dominated by the dominant terms.
We again divide into cases depending on the existence of isolated middle vertices:
\begin{enumerate}
\item First, consider the case where the shape representation has an isolated middle vertex.
	The vertex corresponding to $j_7$ cannot be isolated because there are two distinct random variables $a_{i_1,j_7}$ and  $a_{i_2,j_7}$.
	Hence, only vertices corresponding to $i_1$ and $i_2$ could be isolated.
	Since the case where both of them are isolated is covered in the dominant terms, we assume that only one of them is isolated.
	
	Without loss of generality, let us say the vertex corresponds to $i_1$ is isolated.
	Then, due to the same reasoning as the dominant term analysis, it should be the case that $j$-indices  appearing with $i_1$ in random variables, i.e.,  $j_1,j_2,j_3,j_7$, are either (a) two pairs of two equal indices or (b) all equal.
	\begin{enumerate}
	    \item Let us  handle the first case. Without loss of generality, assume that $j_1=j_2$ and $j_3=j_7$.
	One can place the vertex representing $j_1$ in the intersection   $\U{}\cap \V{}$.
	Now, if there is an additional equality among the indices $j_3(=j_7), j_4,j_5,j_6$, one can place the vertex corresponding to the equal indices in $\U{}\cap \V{}$ as well.
	Then, there are at most $2$ vertices of weight $\logn{n}$ outside a minimum weight separator.
	Consequently, the upper bound reads $\tO{m^{1.5}n}$, which is certainly of $\tO{\sqrt{m}n^{3}}$ since we have assumed $m= \OO{n^{2}}$. 
	
If there is no additional equality among  $j_3, j_4,j_5,j_6$, one can verify that the following shape represents this case.
\begin{figure}[H]
			\centering
			\scalebox{0.5}{
				\begin{tikzpicture}[one/.style={circle,fill=blue!25,draw,font=\sffamily\LARGE\bfseries},two/.style={diamond,fill=red!5,draw,font=\sffamily\LARGE\bfseries}, three/.style={circle split,fill=blue!25,draw,font=\sffamily\LARGE\bfseries}]
				\node[one] (u1) at (-2,2)  {$u_1$};
				\node[one] (u2) at (-2,0) {$u_2$};
				\node[one] (u3) at (-2,-2) {$u_3$};
				\node[two] (w1) at (0,0)  {$w_1$};
			 
				\node[two] (w3) at (0,-2)  {$w_3$};
		 	\node[one] (v2) at (2,0)  {$v_2$};
				\node[one] (v3) at (2,-2)  {$v_3$};
				\path[every node/.style={font=\sffamily},line width=1pt]
				(u2) edge (w3)
				(u3) edge (w3)
				(w3) edge (v2)
				(w3) edge (v3);
					\draw[draw=black] ($(v2.135)+(-0.3,0.3)$) rectangle  ($(v3.-45)+(0.3,-0.3)$);
		\draw[draw=black] ($(u2.135)+(-0.3,0.3)$) rectangle  ($(u3.-45)+(0.3,-0.3)$);
		\draw[draw=black] ($(u1.135)+(-0.3,0.3)$) rectangle  ($(u1.-45)+(0.3,-0.3)$);	 
		\draw[draw=black] ($(u1.135)+(-0.4,0.4)$) rectangle  ($(u1.-45)+(0.4,-0.4)$);	 
		\node[draw=none,fill=none] (10) at (-4,2) {\fontsize{20}{22.4}\selectfont$\U{}\cap \V{}$};
		\node[draw=none,fill=none] (10) at (-4,-1) {\fontsize{20}{22.4}\selectfont$\U{}\setminus \V{}$};
			\node[draw=none,fill=none] (10) at (4,-1) {\fontsize{20}{22.4}\selectfont$\V{}\setminus \U{}$};
			\node[draw=none,fill=none] (10) at (0,1.5) {\fontsize{20}{22.4}\selectfont$\W{}$};
				\end{tikzpicture}}
			\caption{An optimal shape for the case of no additional equality.}
			\label{nonadd}
		\end{figure}
Since $m=\OO{n^2}$, $\{w_3\}$ is a   minimum weight separator.
Hence, the upper bound reads $\tO{m n^2}$, which is certainly $\tO{\sqrt{m} n^{3}}$ since $m=\OO{n^2}$.
\item When  $j_1,j_2,j_3,j_7$ are all equal, let us again place the equal vertex in the intersection $\U{}\cap \V{}$.
By dividing into cases depending on whether there is an additional equality between $j$-indices, one can similarly conclude that there are at most two $\logn{n}$-weighted vertices outside of a vertex separator.
Hence, the bound reads $\tO{m^{1.5} n}$, which is certainly $\tO{\sqrt{m}n^3}$ since $m=\OO{n^2}$.
\end{enumerate}
	
\item Let us move on to the case where there is no isolated vertex outside of $\U{}\cup \V{}$.
	If there are at least two  equalities between indices $j_1,j_2,\dots, j_7$, then at least one of them is between indices $j_1,j_2,\dots, j_6$.
	Let us  place the  vertex corresponding to the  equal indices in $\U{}\cap \V{}$.
	Then any vertex separator will contain the vertex.
	Since the other equality additionally gets rid of an additional $j$-index, we conclude that there are at most $4$ vertices of weight $\logn{n}$ outside of a vertex separator.
	This implies the bound $\tO{mn^2}$, which is certainly $\tO{\sqrt{m}n^3}$.

\item 	Hence, it suffices to consider the case where there is only one equality.
	Let us first assume that the equality is between indices $j_1,j_2,\dots, j_6$.
	Again, one can place the vertex corresponding to the equal indices in $\U{}\cap \V{}$.
	Let us choose two indices among $j_1,j_2,\dots,j_6$ that are not involved in the equality and place  the vertices corresponding to those two indices in different subsets between $\U{}$ and $\V{}$.
	Then, such two vertices admit a path between them, implying  a vertex separator needs to include one of them.
	Hence, there are at most four $\logn{n}$-weighted vertices, again showing the upper bound of $\tO{mn^{2}}$.
	
	Lastly, consider the case where  the equality involves the index $j_7$.
	In that case, the shape in Figure~\ref{subcase} is an optimal shape.	\begin{figure}[H]
	\centering
			\scalebox{0.5}{
				\begin{tikzpicture}[one/.style={circle,fill=blue!25,draw,font=\sffamily\LARGE\bfseries},two/.style={diamond,fill=red!5,draw,font=\sffamily\LARGE\bfseries}, three/.style={circle split,fill=blue!25,draw,font=\sffamily\LARGE\bfseries}]
				\node[one] (u1) at (-2,2)  {$u_1$};
				\node[one] (u2) at (-2,0)  {$u_2$};
				\node[one] (u3) at (-2,-2)  {$u_3$};
				\node[two] (w1) at (0,2)  {$w_1$};
				\node[two] (w3) at (0,-2)  {$w_3$};
				\node[one] (v1) at (2,2)  {$v_1$};
				\node[one] (v2) at (2,0)  {$v_2$};
				\node[one] (v3) at (2,-2)  {$v_3$};
				\path[every node/.style={font=\sffamily},line width=1pt]
				(u1) edge  (w1)
				(u2) edge (w1)
				(w1)  edge  (v2)
				(w1) edge (v1)
				(w3) edge (u3)
				(w3) edge (v3);
				
					\draw[draw=black] ($(v1.135)+(-0.3,0.3)$) rectangle  ($(v3.-45)+(0.3,-0.3)$);
		\draw[draw=black] ($(u1.135)+(-0.3,0.3)$) rectangle  ($(u3.-45)+(0.3,-0.3)$);
			\node[draw=none,fill=none] (8) at (-2,3.5) {\fontsize{20}{22.4}\selectfont$\U{}$};
			\node[draw=none,fill=none] (10) at (2,3.5) {\fontsize{20}{22.4}\selectfont$\V{}$};
			\node[draw=none,fill=none] (10) at (0,3.5) {\fontsize{20}{22.4}\selectfont$\W{}$};
				
				\end{tikzpicture}}
			\caption{An optimal shape for the corner case}
			\label{subcase}
		\end{figure}
	\noindent Since $\{w_1, u_3\}$ is a minimum weight separator, the upper bound is $\tO{\sqrt{m}n^{2.5}}$. 
\end{enumerate}
	Combining all the cases, the remaining cases are dominated by the upper bounds on the dominant terms, and hence the proof for $g_2$ is completed. \qed 	
	
	{\bf 2) Proof for $g_1$:}
	Here, we use a different trick to recover the bound $\max_{\norm{x}=1}g_1(x) =n^3+\tO{mn^{1.5}}$.
	We note that our trick is helpful in recovering the  term $n^3$ up to its coefficient.
	It demonstrates a general principle of using the graph matrix machinery that depending on one's target bound, one might need to modify the target matrix.
	
	Our graph matrix analysis actually yields a better bound of $\max_{\norm{x}=1}g_1(x)\leq   n^3+\tO{mn}$.
	Consider an $n^3\times m$ matrix $B$ whose $((i_1,i_2,i_3),i)$-th entry is set as $a_{i,i_1} a_{i,i_2} a_{i,i_3}$.
	Then, it is straightforward to check that $BB^\top$ is a matrix representation of $g_1$.
	Hence, we have $\max_{\norm{x}=1}g_1(x)\leq \norm{BB^\top}$.
	Here the trick is to consider $\norm{B^\top B}$ instead of $\norm{BB^\top}$. Note that for $i,j\in[m]$ the $(i,j)$-th entry of $B^\top B$ is given as \begin{align}
	\sum_{i_1,i_2,i_3} a_{i,i_1}a_{i,i_2}a_{i,i_3}a_{j,i_1}a_{j,i_2}a_{j,i_3}\,.\label{exp4}
	\end{align}
	We first decompose $B^\top B$ into  the diagonal part $D$ and the non-diagonal part $N$.
	One can easily see that the diagonal part is equal to $n^3\cdot I_m$ where $I_m$ is the $m\times m$ identity matrix.

	Now, we will decompose $N$ into graph matrices. 
	To do that we will decompose each entry in terms of the equality patterns between indices $i_1,i_2,i_3$.
	More specifically, we examine the decomposition  $N=N_{\emptyset} + N_{i_1=i_2} + N_{i_1=i_3} + N_{i_2=i_3}+N_{i_1=i_2=i_3}$
	where for each matrix, the $(i,j)$-th entry is defined as the value of summation \eqref{exp4} restricted to the corresponding equality pattern.
	For instance, for $i\neq j$, the $(i,j)$-th entry of $N_{\emptyset}$ is equal to $\sum_{i_1,i_2,i_3\text{ are distinct}}a_{i,i_1}a_{i,i_2}a_{i,i_3}a_{j,i_1}a_{j,i_2}a_{j,i_3}$ and the $(i,j)$-th entry of  $N_{i_1=i_2}$  is equal to $ 	\sum_{\substack{i_1=i_2\neq i_3}} a_{i,i_1}a_{i,i_2}a_{i,i_3}a_{j,i_1}a_{j,i_2}a_{j,i_3}$.
	
	Then, one can see that each matrix is a graph matrix. 
	More specifically, $N_\emptyset$ is represented by shape $\gamma_{1\dd1}$ in Figure~\ref{t1-1a}, $N_{i_1=i_2}, N_{i_1=i_3}, N_{i_2=i_3}$
	are represented by shape $\gamma_{1\dd2}$ in Figure~\ref{t1-1b} and $N_{i_1,i_2,i_3}$ is represented by shape $\gamma_{1\dd3}$ in Figure~\ref{t1-1c}.

	\begin{figure}[H]
		\centering
		\begin{subfigure}[H]{0.24\textwidth}
			\centering
			\scalebox{0.5}{
				\begin{tikzpicture}[one/.style={circle,fill=blue!25,draw,font=\sffamily\LARGE\bfseries},two/.style={diamond,fill=red!5,draw,font=\sffamily\LARGE\bfseries}, three/.style={circle split, fill=blue!25,draw,font=\sffamily\LARGE\bfseries}]
				\node[two] (u1) at (-2,0) {$u_1$};
				\node[one] (w1) at (0,2)  {$w_1$};
				\node[one] (w2) at (0,0)  {$w_2$};
				\node[one] (w3) at (0,-2)  {$w_3$};
				\node[two] (v1) at (2,0)  {$v_1$};
				\path[every node/.style={font=\sffamily},line width=1pt]
				(u1) edge  (w1)
				(u1) edge  (w2)
				(u1) edge  (w3)
				(v1) edge  (w1)
				(v1) edge (w2)
				(w3) edge (v1);

			\draw[draw=black] ($(v1.135)+(-0.5,0.5)$) rectangle  ($(v1.-45)+(0.5,-0.5)$);
		\draw[draw=black] ($(u1.135)+(-0.5,0.5)$) rectangle  ($(u1.-45)+(0.5,-0.5)$);
			\node[draw=none,fill=none] (8) at (-2,3.5) {\fontsize{20}{22.4}\selectfont$\U{}$};
			\node[draw=none,fill=none] (10) at (2,3.5) {\fontsize{20}{22.4}\selectfont$\V{}$};
			\node[draw=none,fill=none] (10) at (0,3.5) {\fontsize{20}{22.4}\selectfont$\W{}$};
				\end{tikzpicture}		}
			\caption{Figure~\ref{t1-1a}: $\gamma_{1\dd1}$.}
			\label{t1-1a}
		\end{subfigure}\quad
		\begin{subfigure}[H]{0.24\textwidth}
			\centering
			\scalebox{0.5}{
				\begin{tikzpicture}[one/.style={circle,fill=blue!25,draw,font=\sffamily\LARGE\bfseries},two/.style={diamond,fill=red!5,draw,font=\sffamily\LARGE\bfseries}, three/.style={circle split,fill=blue!25,draw,font=\sffamily\LARGE\bfseries}]
				\node[two] (u1) at (-2,0) {$u_1$};
				\node[one] (w1) at (0,2)  {$w_1$};
				\node[one] (w2) at (0,-2)  {$w_2$};
				\node[two] (v1) at (2,0)  {$v_1$};
				\path[every node/.style={font=\sffamily},line width=1pt]
				(u1) edge  (w2)
				(v1) edge (w2);
									\draw[draw=black] ($(v1.135)+(-0.5,0.5)$) rectangle  ($(v1.-45)+(0.5,-0.5)$);
		\draw[draw=black] ($(u1.135)+(-0.5,0.5)$) rectangle  ($(u1.-45)+(0.5,-0.5)$);
			\node[draw=none,fill=none] (8) at (-2,3.5) {\fontsize{20}{22.4}\selectfont$\U{}$};
			\node[draw=none,fill=none] (10) at (2,3.5) {\fontsize{20}{22.4}\selectfont$\V{}$};
			\node[draw=none,fill=none] (10) at (0,3.5) {\fontsize{20}{22.4}\selectfont$\W{}$};
				\end{tikzpicture}
				
			}
			\caption{Figure~\ref{t1-1b}: $\gamma_{1\dd2}$.}
			\label{t1-1b}
		\end{subfigure}
		\begin{subfigure}[H]{0.24\textwidth}
			\centering
			\scalebox{0.5}{
				\begin{tikzpicture}[one/.style={circle,fill=blue!25,draw,font=\sffamily\LARGE\bfseries},two/.style={diamond,fill=red!5,draw,font=\sffamily\LARGE\bfseries}, three/.style={circle split,fill=blue!25,draw,font=\sffamily\LARGE\bfseries}]
				\node[two] (u1) at (-2,0) {$u_1$};
				\node[one] (w1) at (0,0)  {$w_1$};
				\node[] (w2) at (0,-2.3) {};
				\node[] (w2) at (0,2.3) {};
				\node[two] (v1) at (2,0)  {$v_1$};
				\path[every node/.style={font=\sffamily},line width=1pt]
				(u1) edge  (w1)
				(v1) edge (w1);

		\draw[draw=black] ($(v1.135)+(-0.5,0.5)$) rectangle  ($(v1.-45)+(0.5,-0.5)$);
		\draw[draw=black] ($(u1.135)+(-0.5,0.5)$) rectangle  ($(u1.-45)+(0.5,-0.5)$);
			\node[draw=none,fill=none] (8) at (-2,3.5) {\fontsize{20}{22.4}\selectfont$\U{}$};
			\node[draw=none,fill=none] (10) at (2,3.5) {\fontsize{20}{22.4}\selectfont$\V{}$};
			\node[draw=none,fill=none] (10) at (0,3.5) {\fontsize{20}{22.4}\selectfont$\W{}$};
				\end{tikzpicture}
				
			}
			\caption{Figure~\ref{t1-1c}: $\gamma_{1\dd3}$.}
			\label{t1-1c}
		\end{subfigure}
		\caption{Shapes for the matrices $ N_{\emptyset}, N_{i_1=i_2}, N_{i_1=i_3}, N_{i_2=i_3},$ and $N_{i_1=i_2=i_3}$.}
		\label{t1-1}
	\end{figure}
	It is straightforward to check that a minimum weight separator is $\{u_1\}$ for $\gamma_{1\dd1}$, $\{w_2\}$ for $\gamma_{1\dd2}$ and $\{w_1\}$ for $\gamma_{1\dd3}$.
	Thus, our main theorem concludes that $\norm{N_{\emptyset}}=\tO{\sqrt{m}n^{1.5}}$, $\norm{N_{i_1=i_2}}+ \norm{N_{i_1=i_3}}+ \norm{N_{i_2=i_3}}=\tO{mn}$ and $\norm{N_{i_1=i_2=i_3}} = \tO{m}$. 
	Taken collectively, we have 
	$\max_{\norm{x}=1}g_1(x)\leq \norm{B^\top B} \leq \norm{N}+\norm{D} \leq n^3+\tO{mn}$. \qed
		
	\subsection{Fast spectral algorithm for sparse vector recovery}
	\label{fast:sparse}
	In this section, we will  recover the sparse vector recovery analysis in~\cite{hopkins2016fast}.
	We first describe the setting.
	The standing assumption is that for two dimension parameters $n,m$, we have $n \leq  \sqrt{m}$, 
	i.e., $\logn{n} \leq \frac{1}{2}\logn{m}$, where as before $N:=\max\{n,m\}$.
	Given the dimensions, let $b_1,b_2,\dots b_m\in \re^n$ be  random vectors drawn i.i.d. from $\Nc(0,I_{m} )$.
	Consider an $n\times n$ matrix 
	$$B:=\sum_{i\in [m]} \left(\norm{b_i}_2^2 -n \right)  b_i b_i^\top\,.$$
	The goal of this setting is to come up with a good estimate on  $\norm{B}$.
	Proving a sharp upper bound on $\norm{B}$ is a main technical component for analyzing the fast spectral algorithm for sparse vector recovery in~\cite{hopkins2016fast}. 
	We refer readers to Section 4 (especially Lemma 4.5) therein for precise details.
	\begin{theorem}[{\cite[ Section 4]{hopkins2016fast}}] \label{sparse:result}
	With high probability, $\norm{B} \leq 2  m+\tO{n\sqrt{m}}$.
	\end{theorem}
Our approach is based on decomposing $B$ into  several graph matrices based on equality pattern between indices. 
First, note that the $(j_1,j_2)$-th entry of $B$ is expressed as follows:
	$$ B_{j_1,j_2} =    \sum_{\substack{i\in [m],~j_3\in [n]}} \left(   b_{i,j_3}^2 -1 \right)  b_{i,j_1}b_{i,j_2}\,.$$ 
	Now, similarly to what we did in Section~\ref{fournorm}, let us first analyze the dominant terms:
	\begin{itemize}
	    \item ({\bf Distinct indices}:) First, we consider the  case where all indices are distinct, i.e., $|\{j_1,j_2,j_3\}|=3$.
	    Note that this case corresponds to the graph matrix considered in Example~\ref{ex1} where we have shown the upper bound $\tO{n\sqrt{m}}$.
	    Here we recall the illustration of the shape (Figure~\ref{warmup1}) for the reader's convenience.
	    
	\begin{figure}[H]
		\centering
		\scalebox{0.55}{
			\begin{tikzpicture}[one/.style={circle,fill=blue!25,draw,font=\sffamily\LARGE\bfseries},two/.style={diamond,fill=red!5,draw,font=\sffamily\LARGE\bfseries}]
			\node[one] (u1) at (-2,2)  {$u_1$};
			\node[two] (w1) at (0,2)  {$w_1$};
			\node[one] (v1) at (2,2)  {$v_1$};
			\node[one] (w2) at (0,0) {$w_2$};
			\path[every node/.style={font=\sffamily\LARGE\bfseries},line width=1pt]
			(u1) edge  (w1) 
			(w1) edge    (v1)
			(w1) edge node[right]{$2$}  (w2);
				\draw[draw=black] ($(v1.135)+(-0.3,0.3)$) rectangle  ($(v1.-45)+(0.3,-0.3)$);
		\draw[draw=black] ($(u1.135)+(-0.3,0.3)$) rectangle  ($(u1.-45)+(0.3,-0.3)$);
			\node[draw=none,fill=none] (8) at (-2,3.5) {\fontsize{20}{22.4}\selectfont$\U{}$};
			\node[draw=none,fill=none] (10) at (2,3.5) {\fontsize{20}{22.4}\selectfont$\V{}$};
			\node[draw=none,fill=none] (10) at (0,3.5) {\fontsize{20}{22.4}\selectfont$\W{}$};
			\end{tikzpicture}
		}
		
		\caption{The shape of the graph matrix representing the case of distinct indices.
		}
		\label{recall1}
	\end{figure}
	    \item ({\bf Terms with non-zero expected values}:)  Next, we consider the terms having non-zero expected values. Notice that this is only possible when all three indices are the same, i.e., $j_1=j_2=j_3$.
	    The corresponding shape is depicted in Figure~\ref{allsame}:
	    	\begin{figure}[H]
		\centering
		\scalebox{0.55}{
			\begin{tikzpicture}[one/.style={circle,fill=blue!25,draw,font=\sffamily\LARGE\bfseries},two/.style={diamond,fill=red!5,draw,font=\sffamily\LARGE\bfseries}]
			\node[one] (u1) at (-2,2)  {$\nn_1$};
			\node[two] (w1) at (1,2)  {$w_1$};  
		\draw[draw=black] ($(u1.135)+(-0.3,0.3)$) rectangle  ($(u1.-45)+(0.3,-0.3)$);
		\draw[draw=black] ($(u1.135)+(-0.4,0.4)$) rectangle  ($(u1.-45)+(0.4,-0.4)$);
			\node[draw=none,fill=none] (8) at (-2,3.5) {\fontsize{20}{22.4}\selectfont$\U{}\cap \V{}$};
			\node[draw=none,fill=none] (10) at (1,3.5) {\fontsize{20}{22.4}\selectfont$\W{}$};
			\end{tikzpicture}
		}
		
		\caption{A shape representing the case of all equal indices.
		}
		\label{allsame}
	\end{figure}
	    Invoking  the fact that for $X\sim \Nc(0,1)$, $\ex[X^4-X^2]=3-1=2$, the corresponding term is equal to   $ 2 \sum_{j_1\in [n]}\sum_{i\in [m]}   x_{j_1}^2$, which is equal to  $2m$ since $x $ is a  unit vector.
	
	\end{itemize}
Now, it suffices to show that the remaining cases are dominated by the dominant terms above. 
	Following the counting arguments in Section~\ref{fournorm}, we again devise an argument to show that the remaining cases are dominated by the dominant cases above.
	Since we have taken care of the term with non-zero expected value, for the remaining cases, the middle vertex is not isolated.
	Moreover, since we have already handled the case of distinct indices, we may assume that there exists at least one pair of indices that are equal.
	
	\begin{enumerate}
	    \item If $j_1=j_2$, then one can place the vertex corresponding to $j_1=j_2$ in the left/right intersection, which results in the norm bound at most $\tO{\sqrt{nm}}$.
\item Now if there are more than one equality among $j$-indices, then it must be the case that $j_1=j_2$.
Hence, we are left with the case where there is exactly one equality that is not  $j_1=j_2$.
	Without loss of generality, let the unique equality be $j_1=j_3$.
	This case corresponds to the matrix whose $(j_1,j_2)$-th entry is equal to 
	\begin{align}  \sum_{ i\in [m]} \left(   b_{i,j_1}^3 -b_{i,j_1} \right)  b_{i,j_2} \,. \label{poly:1}
	\end{align}
	Again, we decompose the random variable $b_{i,j_1}^3 -b_{i,j_1}$ so that each edge variable is a member of the orthonormal basis:  
\begin{align*}
    b_{i,j_1}^3-b_{i,j_1}= \sqrt{6}\cdot \frac{b_{i,j_1}^3-3b_{i,j_1}}{\sqrt{6}} + 2b_{i,j_1}\,.
\end{align*}
Then, the graph matrix for this case can be decomposed into two parts accordingly.
In either case, it is straightforward to see that for any shape representation, there  must be two (labeled) edges $\{j_1,i\}$ and $\{i,j_2\}$, which implies that there is a path $j_1\to i\to j_2$ from the left vertex to the right vertex. See Figure~\ref{corner}.
  
	\begin{figure}[H]
		\centering
		\scalebox{0.55}{
			\begin{tikzpicture}[one/.style={circle,fill=blue!25,draw,font=\sffamily\LARGE\bfseries},two/.style={diamond,fill=red!5,draw,font=\sffamily\LARGE\bfseries}]
			\node[one] (u1) at (-4,2)  {$u_1$};
			\node[two] (w1) at (0,2)  {$w_1$};
			\node[one] (v1) at (2,2)  {$v_1$};
		 	\path[every node/.style={font=\sffamily\LARGE},line width=1pt]
			(u1) edge node[above]{$3~\text{or}~1$} (w1) 
			(w1) edge    (v1);
				\draw[draw=black] ($(v1.135)+(-0.3,0.3)$) rectangle  ($(v1.-45)+(0.3,-0.3)$);
		\draw[draw=black] ($(u1.135)+(-0.3,0.3)$) rectangle  ($(u1.-45)+(0.3,-0.3)$);
			\node[draw=none,fill=none] (8) at (-4,3.5) {\fontsize{20}{22.4}\selectfont$\U{}$};
			\node[draw=none,fill=none] (10) at (2,3.5) {\fontsize{20}{22.4}\selectfont$\V{}$};
			\node[draw=none,fill=none] (10) at (0,3.5) {\fontsize{20}{22.4}\selectfont$\W{}$};
			\end{tikzpicture}
		}
		
		\caption{Two shapes representing the decomposition of the polynomial \eqref{poly:1}.
		}
		\label{corner}
	\end{figure}

Hence, the minimum weight separator has to contain at least one vertex;
consequently, the norm bounds for those shapes are $\tO{\sqrt{nm}}$, which is again dominated by the dominant cases.
	\end{enumerate}
	Combining the above argument, we conclude that the dominant cases dominate the remaining cases, which finishes the proof of Theorem~\ref{sparse:result}.  \qed

	\subsection{Fast overcomplete tensor decomposition algorithm analysis} \label{fastovertensor}
	
	In this section, we will recover the main technical ingredient for a faster overcomplete tensor decomposition~\cite{hopkins2016fast}.
	In particular, we will recover Proposition 5.5 therein using graph matrices.
	Let us begin with the setting:
	For integers $n,m$ such that $n\leq m$, let $g, a_1, a_2 \dots, a_m\in \re^n$ be random vectors whose entries are drawn i.i.d. from $\Nc(0,1)$.
	Let $T := \sum_{i\in [m]} a_i (a_i\otimes a_i)^\top$.
	Consider an $n^2\times n^2$ matrix 
	\begin{align}
	\mm := \sum_{i_1\neq i_2 \in [m]} \inp{g}{T(a_{i_1}\otimes a_{i_2})} a_{i_1}a_{i_1}^\top \otimes a_{i_2}a_{i_2}^\top \,.
	\end{align}
The goal is then to come up with a good estimate on $\norm{M}$. The following is the upper bound proved in~\cite{hopkins2016fast}.
	\begin{theorem}[{\cite[Proposition 5.5]{hopkins2016fast}}] \label{hopkins}
	With high probability, $\norm{\mm} \leq \tO{n^{2.5}m^{1.5}} $.
	\end{theorem}
	
Again, similar to Section~\ref{fast:sparse}, our approach is based on decomposing $B$ into  several graph matrices based on equality pattern between indices. 
First, note that the $(j_1,j_2,j_3,j_4)$-th entry of $\mm$ is expressed as follows: 
	\begin{align}
 \mm_{j_1,j_2,j_3,j_4}=\sum_{\substack{i_1,i_2,i_3\in[m]:~i_1\neq i_2,\\ j_5,j_6,j_7\in[n]}} g_{j_5}a_{i_3,j_5}a_{i_3,j_6}a_{i_3,j_7}a_{i_1,j_6} a_{i_2,j_7} a_{i_1,j_1}
	a_{i_1,j_2} a_{i_2,j_3}
	a_{i_2,j_4}  \,. \label{ff}
	\end{align}
As before, let us begin with the dominant terms:
\begin{itemize}
\item ({\bf Distinct indices}:)	First, consider the case where both $i$-indices, i.e.,  $i_1,i_2,i_3$, and $j$-indices, i.e., $j_1,j_2,\dots, j_7$, are all distinct.
	Since there is a random variable $g_{j_5}$ depending only on one index $j_5$, shape representations need to have a dummy vertex of weight $0$; then, $g_{j_5}$ can be represented by an edge whose one end is equal to the dummy vertex.
	Having noticed this, one can check that the graph matrix with the following shape represents this case.
	\begin{figure}[H]
		\centering
		\scalebox{0.5}{
			\begin{tikzpicture}[one/.style={circle,fill=blue!25,draw,font=\sffamily\LARGE\bfseries},two/.style={diamond,fill=red!5,draw,font=\sffamily\LARGE\bfseries}, three/.style={square,fill=red!5,draw,font=\sffamily\LARGE\bfseries}]
			\node[one] (u1) at (-4,8)  {$u_1$};
			\node[one] (v1) at (4,8)	{$v_1$};
			\node[two] (w1) at (0,8)   {$w_1$};
			\node[one] (w2) at (0,6)   {$w_2$};
			\node[two] (w3) at (0,4)  {$w_3$};
			\node[one] (w4) at (0,2)   {$w_4$};
			\node[two] (w5) at (0,0)   {$w_5$};
			\node[one] (u2) at (-4,0)  {$u_2$};
			\node[one] (v2) at (4,0)	{$v_2$};
			\node[one] (w6) at (2,4)	{$w_6$};
			\node[circle,draw,dashed,font=\sffamily\LARGE\bfseries] (w7) at (2,6)	{$w_7$};
			\node[draw=none,fill=none] (f) at (-4,9.5) {\fontsize{20}{22.4}\selectfont$\U{}$};
			\node[draw=none,fill=none] (g) at (0,9.5) {\fontsize{20}{22.4}\selectfont$\W{}$};
			
			\node[draw=none,fill=none] (h) at (4,9.5) {\fontsize{20}{22.4}\selectfont$\V{}$};
			\path[every node/.style={font=\sffamily},line width=1pt]
			(u1) edge (w1)
			(v1) edge (w1)
			(w1) edge (w2)
			(w2) edge (w3)
			(w3) edge (w4)
			(w4) edge (w5)
			(w5) edge (u2)
			(w5) edge (v2)
			(w3) edge (w6)
			(w6) edge (w7);
			\draw[draw=black] ($(u1.135)+(-0.2,0.2)$) rectangle  ($(u2.-45)+(0.2,-0.2)$);
			\draw[draw=black] ($(v1.135)+(-0.2,0.2)$) rectangle  ($(v2.-45)+(0.2,-0.2)$);
			\end{tikzpicture}}
		\caption{The shape of the graph matrix for the case of distinct indices. Here, the vertex representing $w_7$ is depicted with a dotted circle, signifying that it is a dummy vertex  of weight $0$.
		This  dummy vertex is required to represent $g_{j_5}$ in \eqref{ff} since this random variable has  a single index $j_5$; with this dummy vertex, $g_{j_5}$ now can be represented by  edge $\{w_6,w_7\}$. }
		\label{hopkinsetal3}
	\end{figure}
	\noindent Since  $\{u_1,u_2\}$ is a minimum weight separator, we obtain the upper bound $\tO{n^{2.5}m^{1.5}}$.
\item ({\bf Terms with non-zero expected values}:) We demonstrate that there are no terms with non-zero expected values. 
To see this, we show that for any shape representation, there is an $i$-index whose vertex representation is not isolated.
First, notice that since $i_1\neq i_2$, $i$-indices cannot be all equal, and hence, there is an $i$-index that is different for the others.
Let the distinct $i$-index be $i_1$ without loss of generality. 
Then by the similar argument to the previous sections, the vertex representation of $i_1$ cannot be isolated as there are three random variables $a_{i_1, j_1}$, $a_{i_1,j_2}$, and $a_{i_1,j_6}$ that contain $i_1$.
Consequently, there are no terms with non-zero expected values.
\end{itemize}
	Now, we consider the remaining cases. Here we consider the two cases depending on whether there is an equality between $i$-indices.
	\begin{enumerate}
	    \item First, consider the case where there is an equality between $i$-indices. Note that since $i_1\neq i_2$, there could be only one equality present; let the equality be $i_2=i_3$ without loss of generality.
	    If there are no additional equalities between indices, one can check that the graph matrix corresponding to the shape in Figure~\ref{fast:1} represents this case:
	    	\begin{figure}[H]
		\centering
		\scalebox{0.5}{
			\begin{tikzpicture}[one/.style={circle,fill=blue!25,draw,font=\sffamily\LARGE\bfseries},two/.style={diamond,fill=red!5,draw,font=\sffamily\LARGE\bfseries}, three/.style={square,fill=red!5,draw,font=\sffamily\LARGE\bfseries}]
			\node[one] (u1) at (-4,8)  {$u_1$};
			\node[one] (v1) at (4,8)	{$v_1$};
			\node[two] (w1) at (0,8)   {$w_1$};
			\node[one] (w2) at (0,6)   {$w_2$};
			\node[two] (w3) at (0,4)  {$w_3$};
			\node[one] (w4) at (-1,1)   {$w_4$}; 
			\node[one] (u2) at (-4,2)  {$u_2$};
			\node[one] (v2) at (4,2)	{$v_2$};
			\node[one] (w6) at (2,4)	{$w_6$};
			\node[circle,draw,dashed,font=\sffamily\LARGE\bfseries] (w7) at (2,6)	{$w_7$};
			\node[draw=none,fill=none] (f) at (-4,9.5) {\fontsize{20}{22.4}\selectfont$\U{}$};
			\node[draw=none,fill=none] (g) at (0,9.5) {\fontsize{20}{22.4}\selectfont$\W{}$};
			
			\node[draw=none,fill=none] (h) at (4,9.5) {\fontsize{20}{22.4}\selectfont$\V{}$};
			\path[every node/.style={font=\sffamily\LARGE},line width=1pt]
			(u1) edge (w1)
			(v1) edge (w1)
			(w1) edge (w2)
			(w2) edge (w3)
			(w3) edge node[right]{$0 ~\text{or} ~2$} (w4) 
			(w3) edge (u2)
			(w3) edge (v2)
			(w3) edge (w6)
			(w6) edge (w7);
			\draw[draw=black] ($(u1.135)+(-0.2,0.2)$) rectangle  ($(u2.-45)+(0.2,-0.2)$);
			\draw[draw=black] ($(v1.135)+(-0.2,0.2)$) rectangle  ($(v2.-45)+(0.2,-0.2)$);
			\end{tikzpicture}}
		\caption{The shape of the graph matrix representing  the case where the only equality between indices is $i_2=i_3$. }
		\label{fast:1}
	\end{figure}
	A minimum weight separator for both shapes is $\{u_1,u_2\}$, which implies the bound of $\tO{n^3 m}$ ($\because$ $w_4$ could be isolated), which is of $\tO{n^{2.5}m^{1.5}}$ since $n\leq m$.
	
	Now let us consider the case where there are equalities present between $j$-indices.
	\begin{enumerate}
	    \item Let us first consider the case where the equalities are between $j$-indices corresponding to the middle vertices in Figure~\ref{fast:1}, i.e., $w_2$, $w_4$, and $w_6$.
	    One can notice that these equalities will only reduce the number of $\logn{n}$-weighted vertices, resulting in a smaller graph matrix bound.
	   \item Now let us consider the case where there is an equality involving a left or right vertex in Figure~\ref{fast:1}, i.e., $u_1$, $u_2$, $v_1$, and $v_2$.
	   
	   \begin{enumerate}
	       \item Notice first that such equalities can possibly make  the vertex $w_3$ in Figure~\ref{fast:1} isolated.
	   On the other hand, to make $w_3$ isolated, the $j$-indices corresponding to its neighbor, i.e., $u_2$, $v_2$, $w_2$, and $w_6$ need to be paired up into two equal indices.
	   This entails at least two additional equalities.
	   If there is an additional equality,  there would be at most four vertices of weight $\logn{n}$ with one of the middle vertices (i.e. $w_4$) isolated and two vertices of weight $\logn{m}$ with one of them isolated,  resulting in the upper bound of $\tO{n^{2.5} m^{1.5}}$.
	   Otherwise if there is no additional equality, there is  a path $u_1\to w_1\to v_1$ from $U$ to $V$, which implies that at least one $\logn{n}$-weighted vertex needs to be in a minimum weight separator, resulting in the bound $\tO{n^{2.5} m^{1.5}}$. 
	    
	   \item  Thus, we may assume that there is no isolated $\logn{m}$-weighted vertex.
	   \begin{enumerate}
	       \item  If there are at least two equalities among $j$-indices, then we get the upper bound of $\tO{n^{3} m} \leq \tO{n^{2.5}m^{1.5}}$.
	       
	       \item    Hence, it suffices to consider the case of single equality. 	   If the equality 	   is between indices corresponding to $u_1,u_2,v_1,v_2$ in Figure~\ref{fast:1},  then  there could be at most five $\logn{n}$-weighted vertices outside of a vertex separator with one of them possibly isolated, which implies the upper bound of $\tO{n^3 m } \leq \tO{n^{2.5}m^{1.5}}$.
	   
	  \item  Hence, we may assume that the equality is between an index corresponding to one of $u_1,u_2,v_1,v_2$ and an index corresponding to one of  $w_2,w_4,w_6$.
	   In all cases, one can easily check from Figure~\ref{fast:1} that either the path $u_1\to w_1\to v_1$ or  $u_2\to w_3\to v_2$ is unaltered, and hence at least one of the $\logn{n}$-weighted vertex has to be included in a vertex separator.
	   This again implies that  there are at most five $\logn{n}$-weighted vertices outside of a vertex separator with one of them possibly isolated, implying the bound $\tO{n^3 m } \leq \tO{n^{2.5}m^{1.5}}$.
	   \end{enumerate}
	  
	   \end{enumerate}
	   
	\end{enumerate}
	
	\item Now we consider the case where there is no equality between $i$-indices.
	This case corresponds to having no equalities between $\logn{m}$-weight vertices in Figure~\ref{hopkinsetal3}, and hence, one can  notice  that there cannot be isolated middle vertices.
	
	  If there are at least two equalities present among $j$-indices, then   we get the upper bound of $\tO{n^{2.5} m^{1.5}}$ as there are at most five  $\logn{n}$-weighted vertices and at most three $\logn{m}$-weighted vertices.
	   
	   Hence, it suffices to consider the case of a single equality.
	   Without loss of generality, assume that the index corresponding to $u_1$ in Figure~\ref{hopkinsetal3} is involved in the equality.
	   Then the path $u_2\to w_5\to v_2$ in Figure~\ref{hopkinsetal3} from $u_2$ to $v_2$ is unaltered, and hence either $u_2$ or $v_2$ has to be included in a vertex separator.
	   Hence, there could be at most five $\logn{n}$-weighted vertices outside of a vertex separator, which again implies the upper bound of $\tO{n^{2.5} m^{1.5} }$.
	\end{enumerate}
	The above arguments demonstrate that the remaining cases are dominated by the dominant terms, which concludes the proof of Theorem~\ref{hopkins} using graph matrices.
	\qed

	\section{Conclusion and further studies}
	In this paper, we analyzed graph matrices, a general class of random matrices which naturally arise when analyzing the Sum-of-Squares hierarchy. Using the trace power method, we proved probabilistic norm bounds on graph matrices, which are tight up to a poly-logarithmic factor. Notably, we showed that for all shapes $\alpha$, the minimum size (weight) of a separator between the left and right sides of $\alpha$ is the key quantity that characterizes the norm bound on $M_{\alpha}$. We also showed that graph matrices can be used to simplify the proofs of several technical statements in the literature \cite{barak2012hypercontractivity, hopkins2016fast,ge2015decomposing} about the $2 \to 4$ norm of a random operator and tensor decomposition.
	
	This paper leaves several interesting open problems.
	To name a few, one question is to improve our current norm bounds. It is very likely that with a more careful analysis, the poly-logarithmic factor can be reduced or removed. Another interesting question would be to determine the limit of the mean density of the spectra of graph matrices, finding analogues of Wigner's semicircle law \cite{wigner1958distribution} and/or Girko's circular law \cite{girko} for graph matrices.
	
	\appendix

		\section{Lower bounds}\label{lowerboundsection}
		In this appendix, we show that the norm bounds we have obtained are tight up to a factor of $\plog(n)$.
		\begin{theorem}\label{lowerboundtheorem}
			If $\alpha$ is a shape such that for all vertices $w$ in $W_{\alpha} \setminus \iso$, there is a path from $w$ to either $U_{\alpha}$ or $V_{\alpha}$ in $\alpha$, then with high probability $\norm{M_{\alpha}}$ is $\Omega\left(n^{\frac{w(V(\alpha)) - w(\smin) + w(\iso)}{2}}\right)$.
		\end{theorem}
		\begin{remark}
			If $W_{\alpha}$ has vertices which are not isolated and are not connected to $U_{\alpha}$ or $V_{\alpha}$, there is a non-negligible chance that $M_{\alpha}$ has considerably smaller norm than expected. To see this, let $\alpha_{0}$ be the part of $\alpha$ which is not isolated but is not connected to $U_{\alpha}$ or $V_{\alpha}$ and let $\alpha_1$ be the remaining part of $\alpha$. We have that $M_{\alpha} \approx M_{\alpha_0}M_{\alpha_1}$ where $M_{\alpha_0}$ is a scalar depending on the input graph $G$ which has a non-negligible chance of being close to $0$.
		\end{remark}
		\begin{proof}
		To prove this theorem, we analyze the Frobenius norm of $M_{\alpha}$.
		    \begin{definition}
		    The Frobenius norm of an $m \times n$ matrix $M$ is $\norm{M}_{F} = \sqrt{\sum_{i = 1}^{m}{\sum_{j = 1}^{n}{M^{2}_{ij}}}}$.
		    \end{definition}
		    \begin{lemma}\label{lem:frobenius}
		    For any matrix $M$, if the singular value decomposition of $M$ is $M = \sum_{i=1}^{r}{\lambda_i{u_i}{v^\top_i}}$ then $\norm{M}^{2}_{F} = \sum_{i=1}^{r}{\lambda^2_i}$.
		    \end{lemma}
		    \begin{proof}
		    Observe that 
		    \[
		    \norm{M}^{2}_{F} = \sum_{i=1}^{r}{\sum_{j=1}^{r}{{\lambda_i}{\lambda_j}(u_i \cdot u_j)(v_i \cdot v_j)}} = \sum_{i=1}^{r}{\lambda^2_i}\,.
		    \]
		    \end{proof}
		    \noindent Using Lemma~\ref{lem:frobenius}, one can develop a lower bound on the spectral norm of a matrix in terms of the Frobenius norm of the matrix.
		    \begin{corollary} \label{cor:speclower}
		    For any rank $r$ matrix $M$ where $r \geq 1$, $\norm{M} \geq \sqrt{\frac{\norm{M}^{2}_{F}}{r}}$.
		    \end{corollary}
		    With this in mind, our strategy for proving lower bounds on $\norm{M_{\alpha}}$ is as follows:
		    \begin{mdframed}
		    \begin{enumerate}
		        \item Take a matrix $M'_{\alpha}$ such that $M'_{\alpha}$ has rank at most $n^{w(\smin)}$ and with high probability, $\norm{M'_{\alpha} - M_{\alpha}}$ is negligible.
		        \item Observe that $\E\left[\norm{M'_{\alpha}}^{2}_{F}\right]$ is $\Theta(n^{w(V(\alpha)) + w(\iso)})$.
		        \item Show that with high probability, 
		        \[
		        \left|\norm{M'_{\alpha}}^{2}_{F} - \E\left[\norm{M'_{\alpha}}^{2}_{F}\right]\right| \ll n^{w(V(\alpha)) + w(\iso)}\,.
		        \]
		        \item Combining above three statement, the desired probabilistic lower bound  on $\norm{\M_\alpha}$ follows from Corollary~\ref{cor:speclower}.
		    \end{enumerate}
		    \end{mdframed}

		    For simplicity, we first implement the second and third steps with $M_{\alpha}$ rather than $M'_{\alpha}$. We then describe $M'_{\alpha}$ and explain why the same techniques can be used for $M'_{\alpha}$.
		    
		    To analyze $\norm{M_{\alpha}}^{2}_{F}$, we observe that $\norm{M_{\alpha}}^{2}_{F} = \tr\left(M_{\alpha}M_{\alpha}^\top\right)$. Expanding out the term using the definition of $M_{\alpha}$, we have the following proposition:
		    \begin{proposition}\label{prop:Frobeniusdecomposition}
		   For any shape $\alpha$, $\tr\left(M_{\alpha}M_{\alpha}^\top\right) = \sum_{i}{{c_i}M_{\beta_i}}$
		   where each shape $\beta_i$ is obtained through the following steps and each $c_i$ is a constant.
		        \begin{enumerate}
		            \item Start with two copies $\alpha_1$ and $\alpha_2$ of $\alpha$.
		            \item Take some set $E_{\mathrm{constr}}$ of constraint edges between $V(\alpha_1)$ and $V(\alpha_2)$ where $E_{\mathrm{constr}}$ automatically has constraint edges between $U_{\alpha_1}$ and $U_{\alpha_2}$ and between $V_{\alpha_1}$ and $V_{\alpha_2}$.
		            \item Take $V(\beta_i)$ to be $V(\alpha_1) \cup V(\alpha_2)$ where we contract all edges in $E_{\mathrm{constr}}$. Take $U_{\beta_i} = V_{\beta_i} = \emptyset$.
		            \item Consider $E(\alpha_1) \cup E(\alpha_2)$ as a multi-set. If there are two edges between the same pair of vertices (or set of vertices if we have hyperedges) in $V(\alpha')$ and these edges have labels $l_1$ and $l_2$ where $l_1 \leq l_2$ then replace these two edges with a single edge with some label in $[l_2 - l_1,l_2 + l_1]$ (where an edge with label $0$ corresponds to deleting the edge). Take $E(\beta_i)$ to be the resulting set of labeled edges.
		            
		            We do this because $h_{l_1}h_{l_2} = \sum_{j = l_2 - l_1}^{l_2 + l_1}{c_{l_1,l_2,j}h_j}$ for some coefficients $c_{l_1,l_2,j} = \E_{\Omega}[h_{l_1}h_{l_2}h_{j}]$.
		        \end{enumerate}
		    \end{proposition}
		    We now show the following lemma about the shapes $\beta_i$ which we obtain in this way.
		    \begin{lemma}\label{lem:lowerboundkey}
		    Let $\alpha$ be any shape and $\beta_i$ be a shape obtained by decomposing $\tr\left(M_{\alpha}M_{\alpha}^\top\right)$. 
		    Let $\isoo(\beta_i)$ be the set of isolated vertices of $\beta_i$.
		    Then we have 
		    \begin{align*}
		        w(V(\beta_i)) + w(\isoo(\beta_i)) \leq 2(w(V(\alpha)) + w(\iso))\,.
		    \end{align*} Moreover, if $\alpha$ is a shape such that for all vertices $w$ in $W_{\alpha} \setminus \iso$, there is a path from $w$ to either $U_{\alpha}$ or $V_{\alpha}$ in $\alpha$, then for any $\beta_i$ such that $E(\beta_i) \neq \emptyset$, we have
		    \begin{align*}
		        w(V(\beta_i)) + w(\isoo(\beta_i)) < 2(w(V(\alpha)) + w(\iso))\,.
		    \end{align*}
		    \end{lemma}
		    \begin{proof}
Following the argument in Section~\ref{reduc:2}, observe that each vertex $w \in \iso$ just multiplies $M_{\alpha}$ by roughly $n$ (or more generally $n^{w(v)}$). Thus, it is sufficient to consider the case when $\iso$ is empty.

Take $X$ to be the set of vertices in $V(\beta_i)$ which resulted from the contraction of a constraint edge and take $Y = V(\beta_i) \setminus X$. We have that $2w(V(\alpha)) = w(V(\alpha_1)) + w(V(\alpha_2)) = 2w(X) + w(Y)$. Now observe that any $v \in Y$ must be incident to at least one edge in $E(\beta_i)$. Thus, $\isoo(\beta_i) \subseteq X$ so 
\begin{align*}
w(V(\beta_i)) + w(\isoo(\beta_i)) &= 2w(X \cap \isoo(\beta_i)) + w(X \setminus \isoo(\beta_i)) + w(Y)\\
&\leq 2w(X) + w(Y) = 2w(V(\alpha))\,.    
\end{align*}

To show the moreover statement, it is sufficient to show that whenever $E(\beta_i) \neq \emptyset$, $X \setminus \isoo(\beta_i) \neq \emptyset$. To see this, let $w$ be a vertex in $V(\beta_i) \setminus \isoo(\beta_i)$. 
If $w \in X$, then we are done so we may assume that $w \in Y$.
Now consider $\alpha_1$ and $\alpha_2$ before we contract all of the constraint edges to form $\beta_i$. 
$w$ either came from $V(\alpha_1)$ or $V(\alpha_2)$, so without loss of generality we can assume that $w$ came from $V(\alpha_1)$.
By assumption, there is a path $P$ in $\alpha_1$ from $w$ to either $U_{\alpha_1}$ or $V(\alpha_1)$. 
Let $v$ be the first vertex on this path which is incident to a constraint edge and let $u$ be the preceding vertex.
Note that $u$ and $v$ must exist because all of the vertices in $U_{\alpha_1}$ and $V_{\alpha_1}$ are incident to constraint edges while $w$ is not incident to a constraint edge. 
When we contract all of the constraint edges, $v$ (or rather the vertex resulting from contracting $v$ and the other endpoint of its constraint edge) will be in $X$. 
However, the edge between $u$ and $v$ will not vanish, so $v$ will not be isolated. Thus, $X \setminus \isoo(\beta_i) \neq \emptyset$, as needed.
		    \end{proof}
		    \begin{corollary}
		    $\E\left[\norm{M_{\alpha}}^{2}_{F}\right]$ is $\Theta\left(n^{w(V(\alpha)) + w(\iso)}\right)$ and with high probability, $\left|\norm{M_{\alpha}}^{2}_{F} - \E\left[\norm{M_{\alpha}}^{2}_{F}\right]\right|$ is $\widetilde{O}\left(n^{w(V(\alpha)) + w(\iso) - \frac{w_{min}}{2}}\right)$ where $w_{min}$ is the minimum weight of a vertex in $V(\alpha)$.
		    \end{corollary}
		    \begin{proof}
		    For the first part, observe that 
		    \begin{align*}
		        \ex\left[\norm{M_{\alpha}}^{2}_{F}\right] = \sum_{i: E(\beta_i) = \emptyset}{{c_i}M_{\beta_i}}\,.
		    \end{align*}
		    For these $\beta_i$, $M_{\beta_i}$ is a scalar which is $\Theta\left(n^{w(V(\beta_i))}\right) = \Theta\left(n^{\frac{w(V(\beta_i)) + w(\isoo(\beta_i))}{2}}\right)$. By Lemma \ref{lem:lowerboundkey}, $w(V(\beta_i)) + w(\isoo(\beta_i)) \leq 2w(V(\alpha))$ and this is tight as if we have constraint edges between each vertex in $V(\alpha_1)$ and its counterpart in $V(\alpha_2)$ then we will have that $w(V(\beta_i)) = w(\isoo(\beta_i)) = w(V(\alpha))$ and $E(\beta_i) = \emptyset$.
		    
		    For the second part, observe that $\norm{M_{\alpha}}^{2}_{F} = \E\left[\norm{M_{\alpha}}^{2}_{F}\right] = \sum_{i: E(\beta_i) = \emptyset}{{c_i}M_{\beta_i}}$
		    Combining Theorem \ref{boundgeneral} and Lemma \ref{lem:lowerboundkey}, with high probability each of these terms is $\widetilde{O}\left(n^{w(V(\alpha)) + w(\iso) - \frac{w_{min}}{2}}\right)$, as needed.
		    \end{proof}

		     We now implement the second and third steps with an approximation $M'_{\alpha}$.
		     For the approximation, we adopt the idea from \cite[ Section 6]{FinalPlantedClique} and consider a low rank approximation $M'_{\alpha}$ of $M_{\alpha}$ by decomposing $\alpha$ into left, middle, and right parts $\sigma$, $\tau$, and ${\sigma'}^\top$ based on the leftmost and rightmost minimum weight vertex separators of $\alpha$. 
\begin{definition}
Given a shape $\alpha$, we define the leftmost minimum weight vertex separator $S_{\alpha}$ to be the set of vertices such that 
\begin{enumerate}
\item $S_{\alpha}$ is a minimum weight vertex separator between $\U{\alpha}$ and $\V{\alpha}$.
\item For any other minimum weight vertex separator $S'$, $S_{\alpha}$ is a vertex separator between $\U{\alpha}$ and $S'$.
\end{enumerate}
Similarly, we define the rightmost minimum weight vertex separator $T_{\alpha}$ to be the set of vertices such that 
\begin{enumerate}
\item $T_{\alpha}$ is a minimum weight vertex separator between $U_{\alpha}$ and $V_{\alpha}$.
\item For any other minimum weight vertex separator $S'$, $T_{\alpha}$ is a vertex separator between $S'$ and $V_{\alpha}$.
\end{enumerate}
\end{definition}
As shown in \cite{FinalPlantedClique}, the leftmost and rightmost minimum vertex separators are well-defined.
\begin{definition}
Given a shape $\alpha$, we define the left, middle, and right parts of $\alpha$ as follows:
\begin{enumerate}
\item We define the left part $\sigma$ of $\alpha$ to be the part of $\alpha$ between $U_{\alpha}$ and $S_{\alpha}$.
\item We define the middle part $\tau$ of $\alpha$ to be the part of $\alpha$ between $S_{\alpha}$ and $T_{\alpha}$.
\item We define the right part ${\sigma'}^\top$ of $\alpha$ to be the part of $\alpha$ between $T_{\alpha}$ and $V_{\alpha}$.
\end{enumerate}
\end{definition}
\begin{definition}
Given a shape $\alpha$, we define $M'_{\alpha} = M_{\sigma}M_{\tau}M_{{\sigma'}^\top}$ where $\sigma$, $\tau$, and ${\sigma'}^\top$ are the left, middle, and right parts of $\alpha$.
\end{definition}

Using the intersection tradeoff lemma  \cite[Lemma 7.12]{FinalPlantedClique}, we have the following lemma
\begin{lemma}
For any shape $\alpha$, letting $\sigma$, $\tau$, and ${\sigma'}^\top$ be the left, middle, and right parts of $\alpha$, 
\[
M_{\sigma}M_{\tau}M_{{\sigma'}^\top} - M_{\alpha} = \sum_{i}{c_{i}M_{\alpha_i}}
\]
for some shapes $\alpha_i$ and coefficients $c_i$ such that
\begin{enumerate}
\item $\forall i, |c_i|$ is $O(1)$.
\item For each shape $\alpha_i$, $\U{\alpha_i} = \U{\alpha}$, $\V{\alpha_i} = \V{\alpha}$, and $w(V(\alpha_i)) + w(\isoo(\alpha_i)) < w(V(\alpha)) + w(\iso)$, where $\isoo(\alpha_i)$ is the set of isolated vertices of $\alpha_i$.
\end{enumerate}
\end{lemma}
\begin{remark}
Intuitively, this lemma says that $M_{\sigma}M_{\tau}M_{{\sigma'}^\top}$ is a good approximation to $M_{\alpha}$.
\end{remark}
The rank of $M_{\sigma}$, $M_{\tau}$, and $M_{{\sigma'}^\top}$ are all at most $n^{w(S_{min})}$, so it is sufficient to show that $\E\left[\norm{M'_{\alpha}}^{2}_{F}\right]$ is $\Theta(n^{w(V(\alpha)) + w(\iso)})$ and $\left|\norm{M'_{\alpha}}^{2}_{F} - \E\left[\norm{M'_{\alpha}}^{2}_{F}\right]\right| \ll n^{w(V(\alpha)) + w(\iso)}$. To do this, we write $M'_{\alpha} = M_{\alpha} + \sum_{i}{c_{i}M_{\alpha_i}}$ and observe that
\begin{align*}
\tr\left(M'_{\alpha}{M'}_{\alpha}^\top \right) = \tr\left(M_{\alpha}{M}_{\alpha}^\top\right) + \sum_{i}{{c_i} \tr\left(M_{\alpha_i}{M}_{\alpha}^\top \right)} + 
\sum_{i'}{{c_{i'}}\tr\left(M_{\alpha}{M}_{\alpha_{i'}}^\top \right)} + \sum_{i,i'}{{c_i}{c_{i'}}tr\left(M_{\alpha_i}{M}_{\alpha_{i'}}^\top \right)}
\end{align*}
We have already analyzed $\tr\left(M_{\alpha}{M}_{\alpha}^\top \right)$. To show that the other terms are negligible, we combine the following lemma with Theorem \ref{boundgeneral}.
\begin{lemma}
Let $\alpha$ and $\alpha'$ be any shapes such that $U_{\alpha} = U_{\alpha'}$ and $V_{\alpha} = V_{\alpha'}$.
Let $\beta_i$ be a shape obtained from decomposing $\tr\left(M_{\alpha}M_{\alpha'}^\top \right)$.
Letting $\isoo(\beta_i)$ be the set of isolated vertices of $\beta_i$, we have 
\begin{align*}
    w(V(\beta_i)) + w(\isoo(\beta_i)) \leq w(V(\alpha)) + w(\isoo(\alpha)) + w(V(\alpha')) + w(\isoo(\alpha'))\,,
\end{align*} 
where $\isoo(\alpha)$ and $\isoo(\alpha')$ are the isolated middle vertices of $\alpha$ and $\alpha'$ respectively.
\end{lemma}
\begin{proof}
This can be proved in the same way as the first part of Lemma \ref{lem:lowerboundkey}. Again, it is sufficient to consider the case when $\alpha$ and $\alpha'$ have no isolated middle vertices. Taking $X$ to be the set of vertices in $V(\beta_i)$ which resulted from the contraction of a constraint edge and taking $Y = V(\beta_i) \setminus X$, we have that $w(V(\alpha)) + w(V(\alpha')) = 2w(X) + w(Y)$. Now observe that any $v \in Y$ must be incident to at least one edge in $E(\beta_i)$. Thus, $\isoo(\beta_i) \subseteq X$ so 
\begin{align*}
w(V(\beta_i)) + w(\isoo(\beta_i)) &= 2w(X \cap \isoo(\beta_i)) + w(X \setminus \isoo(\beta_i)) + w(Y)\\
&\leq 2w(X) + w(Y) = w(V(\alpha))+ w(V(\alpha'))\,,
\end{align*}
which concludes the proof.
\end{proof}
This completes the proof of Theorem~\ref{lowerboundtheorem}.
\end{proof}
\begin{remark}
A careful reader might notice that \cite{FinalPlantedClique} only analyzes the graph matrices defined in Section \ref{simpledef} and not the generalized graph matrices in Section~\ref{gendef}. 
That said, these results (the existence of leftmost and rightmost minimum weight vertex separators and the intersection tradeoff lemma) hold for generalized graph matrices as well and this can be shown by generalizing the arguments in \cite{FinalPlantedClique}.
\end{remark}
	
	\section{Alternative proof of Lemma~\ref{mingeneral2}} \label{altproof}
Here, we provide an alternative proof of  Lemma~\ref{mingeneral2}. 
We first recall the statement of Lemma~\ref{mingeneral2} for readers' convenience: For any shape $\alpha$ such that $U_{\alpha} \cap V_{\alpha} = \emptyset$, for any well-behaved $C \in \mathcal{C}_{(\alpha,2q)}$ such that $\val{C} \neq 0$, the total weight of the constraint edges is at least $q\cdot \w{\W{\alpha}}+ (q-1)\cdot  \w{\smin } $.

To simplify our analysis, we follow the first paragraph of the proof in the main text, and reduce the statement to showing that:
\begin{align}\label{st}
\text{The total weight of the ``additional'' constraint edges}\geq (q-1)\cdot  \w{\smin }. 
\end{align}
Our proof proceeds by the following three steps:
	\begin{enumerate}
	    \item {\bf Preprocessing the constraint graph:} We will first preprocess the constraint graph so that the task of  counting the number of ``additional'' constraint edges becomes simpler.
	    
	    \item {\bf Computing the total weight of the ``additional'' constraint edges:} After the preprocessing step, one can easily count the number of ``additional'' constraint edges. 
	    Consequently, one can easily compute the total weight of the ``additional'' constraint edges.  
	    \item {\bf Lower bound using vertex separators:} Based on the previous two steps, we will then come up with a lower bound on the total weight in terms of $\w{\smin}$.  This lower bound will exactly match \eqref{st}.
	\end{enumerate}
Below, we show the details of the above three steps:
\subsection{Preprocessing the constraint graph}
	
We first define the notion of equivalence for clarity:
\begin{definition}
In a constraint graph (Definition~\ref{def:constr}), we say two vertices $u$ and $v$ are \emph{equivalent} if $\phi(u)=\phi(v)$. We write $u\sim v$ to denote that $u$ and $v$ are equivalent.
\end{definition}
\noindent With this notion of equivalence,  we now describe our preprocessing steps:

		\begin{enumerate}
			\item We choose a minimum weight separator $\smin$ between $U_{\alpha_1}$ and $V_{\alpha_1}$.
			\item For each vertex $v$ which is not equivalent to some vertex in $\smin$, we rearrange the constraint edges for the vertices which are equivalent to $v$ so that these edges form a path from left to right (with $\smin$ as the reference point).
			\item For each vertex $v \in \smin$, we rearrange the constraint edges for the vertices which are equivalent to $v$ so that these edges form a path from left to right starting at $v$. We then add an additional superfluous constraint edge to make this path into a cycle. Note that if $v$ was not incident to any other constraint edges (which can happen only if $v \in U_{\alpha_1}$ or $v \in V_{\alpha_1}$) then this superfluous constraint edge is a loop.
		\end{enumerate}
Due to the superfluous constraint edges we added at the last step, the statement \eqref{st} we wanted to show now becomes:
	\begin{align}\label{st2}
	    \text{The total weight of ``additional'' constraint edges}\geq q \cdot \w{\smin}\,.
	\end{align}
We now count the number of ``additional'' constraint edges in the preprocessed constraint graph.
 
\subsection{Computing the total weight of ``additional'' constraint edges} 
		
	Due to the preprocessing step, the task of counting the number of ``additional'' constraint edges becomes simpler.
More specifically, note that  the constraint edges between equivalent vertices now form a path from left to right. 
Hence, each vertex in the constraint graph is incident to either $0$, $1$ or $2$ constraint edges. Based on this, we define:
\begin{definition}
For each vertex $v\in \alpha$, let $\numc{v}^{(i)}$ ($i=0,1,2$) be the number of copies of $v$ that are incident to $i$ non-loop constraint edges. Let $\numc{v}^{(loop)}$ be the number of copies of $v$ which are incident to a single loop constraint edge. Note that $\numc{v}^{(loop)}$ is always $0$ or $1$ and we only have that $\numc{v}^{(loop)} = 1$ if there is a copy of $v$ in $\smin$ and this copy of $v$ was originally not incident to any constraint edges.
\end{definition} 	 
\noindent From this definition, it is straightforward to verify the following:
\begin{claim} \label{prop:counting}
The total weight of ``additional'' constraint edges is equal to
\begin{align} \label{counting:0}
   \frac{1}{2} \sum_{v\in \U{\alpha}\cup \V{\alpha} } \w{v}\cdot \left[ \numc{v}^{(1)} + 2 \numc{v}^{(loop)} + 2 \numc{v}^{(2)}  \right] +    \frac{1}{2} \sum_{v\in \W{\alpha} } \w{v} \cdot  \numc{v}^{(2)}  
\end{align}
\end{claim}
\begin{proof}
First, every incident constraint edge to a copy of a vertex in $ \U{\alpha}\cup \V{\alpha}$ is an ``additional'' constraint edge.
Hence, the total weight of ``additional'' constraint edges among the copies of $ \U{\alpha}\cup \V{\alpha}$ is
\begin{align} \label{counting:1}
   \frac{1}{2} \sum_{v\in \U{\alpha}\cup \V{\alpha} } \w{v}\cdot\left[   0\cdot \numc{v}^{(0)}+   1\cdot \numc{v}^{(1)} + 2 \cdot \numc{v}^{(loop)} + 2\cdot   \numc{v}^{(2)}  \right] \,,
\end{align}
where we divide the summation by $2$ to avoid double counting.

Next, consider  vertices that are copies of $\W{\alpha}$. 
Since we are counting only ``additional'' constraint edges, we need to reduce  the number of incident constraint edges by one.
Hence, the total weight among the copies of $\W{\alpha}$ is equal to
\begin{align} \label{counting:2}
   \frac{1}{2} \sum_{v\in \W{\alpha}  } \w{v}\cdot\left[ 0\cdot \numc{v}^{(0)}+  0\cdot \numc{v}^{(1)}+ 1\cdot \numc{v}^{(2)}  \right] \,.
\end{align}
Adding the two terms \eqref{counting:1} and \eqref{counting:2}, \eqref{counting:0} follows.
\end{proof}
\subsection{Lower bound using vertex separators}
Now the last step develops a lower bound on \eqref{counting:0} in terms of $\w{\smin}$. 
To that end, we first classify the vertices that are incident to only one constraint edge:
		\begin{definition}[Endpoints]
		Recall that $\smin$ is the minimum weight vertex separator of $\alpha_1$ which we chose as a reference point. For a vertex $v \notin \smin$ that is incident to only one constraint edge, we say
		\begin{itemize}
		    \item $v$ is a left endpoint if there is no vertex $u$ such that $u \sim v$ and $u$ is to the left of $v$; and 
		    \item $v$ is a right endpoint if there is no vertex $w$ such that $v \sim w$ and $w$ is to the right of $v$.
		\end{itemize} 
		For each vertex $v \in \smin$ which is only incident to a superfluous constraint edge which is a loop,  
		\begin{itemize}
		    \item If $v \in U_{\alpha_1}$ then we consider $v$ to be a left endpoint in $\alpha_1$ and a right endpoint in $\alpha^T_{2q}$.
		    \item If $v \in V_{\alpha_1}$ then we consider $v$ to be a right endpoint in $\alpha_1$ and a left endpoint in $\alpha^T_{1}$.
		\end{itemize} 
		\end{definition}
\noindent 
With this left/right endpoint classification, one can easily prove the following crucial property:
		\begin{claim} \label{claim:edge}
	Let $w_1$ and $w_2$ be the vertices that belong to the same copy $\alpha_i$ in the constraint graph. Assume that $\{w_1,w_2\} \not\subset \U{\alpha_i}$ and $\{w_1,w_2 \} \not\subset \V{\alpha_i}$.
	Then $\{w_1,w_2\} \notin E(\alpha_i)$ if any of the following hold:
	\begin{itemize}
	    \item $w_1$ is a left endpoint and $w_2$ is a right endpoint.
	    \item $w_1\in \U{\alpha_i^\top}$, $w_1$ is not incident to a constraint edge, and $w_2$ is a left endpoint.
	    \item $w_1\in \V{\alpha_i^\top}$, $w_1$ is not incident to a constraint edge, and $w_2$ is a right endpoint.
	      \item $w_1 \in \U{\alpha_i^\top}$ and  $w_2\in \V{\alpha_i^\top}$ and both are not incident to a constraint edge. 
	      \end{itemize} 
	The same conclusion holds when we replace $\alpha_i$ with $\alpha_i^\top$ in the statement. 
		\end{claim}
		\begin{proof}
		Let us consider the first scenario: say $w_1$ is a left endpoint and $w_2$ is a right endpoint.
		Then all of the vertices which are equivalent to $w_1$ are to the right of $w_1$ and all vertices which are equivalent to $w_2$ are to the left of $w_2$ so if $\{w_1,w_2\} \in E(\alpha_i)$ (and $\{w_1,w_2\} \not\subset \U{\alpha_i}$ and $\{w_1,w_2 \} \not\subset \V{\alpha_i}$) then there is no way for the copy $\{w_1,w_2\}$ to appear more than once so we have $\val{C}=0$.
			
			Following similar logic, we can show that in each of the other cases, there is no way for the copy $\{w_1,w_2\}$ to appear more than once so we have $\val{C}=0$.
		\end{proof}
\noindent 
Using Claim~\ref{claim:edge}, we can characterize all left/right traversing paths within a single copy of $\alpha$ in the constraint graph.
More specifically, for $i=1,2,\dots, q$, any path $P$ from $\U{\alpha_i}$ to $\V{\alpha_i}$ within $\alpha_i$ (respectively,  from $\U{\alpha_i^\top}$ to $\V{\alpha_i^\top}$ within $\alpha_i^\top$)  satisfies  one of the following:
\begin{itemize}
    \item $P$ contains at least one vertex that is incident with two constraint edges.
    \item All vertices in $P \setminus V_{\alpha_i}$  (respectively, $P \setminus V_{\alpha_i^\top}$) are left endpoints.
    \item All vertices in $P \setminus U_{\alpha_i}$ (respectively, $P \setminus V_{\alpha_i^\top}$) are right endpoints.
\end{itemize}  
Consequently, one can construct the following vertex separators:
\begin{itemize}
    \item Let $S_i$ be the subset of $\alpha_i$ consisting of left endpoints in $\U{\alpha_i}$, vertices in $\alpha_i$ that are incident to two constraint edges, and right endpoints in $\V{\alpha_i}$. Then $S_i$ is a vertex separator of $\alpha_i$.
    \item Similarly, define $S'_i$ to be the subset of $\alpha_i^{\top}$ consisting of left endpoints in $\U{\alpha_i^{\top}}$, vertices in $\alpha_i^{\top}$ that are incident to two constraint edges, and right endpoints in $\V{\alpha_i^{\top}}$. Then, $S'_i$ is a vertex separator of $\alpha_i^{\top}$.
\end{itemize}
Then, the following identity is straightforward from the construction:
\begin{align}\label{appenb:1}
    \eqref{counting:0} =  \frac{1}{2}\cdot \sum_{i=1}^q \left[ \w{S_i} + \w{S'_i}\right]\,.
\end{align} 
Now since $S_i,S'_i$ are vertex separators of $\alpha_i$ and $\alpha_i^\top$, it follows that 
\begin{align} \label{appenb:2}
     \w{S_i},\w{S'_i} \geq \w{\smin}\,.
\end{align}
Combining \eqref{appenb:1} with \eqref{appenb:2}, we obtain \eqref{st2}, which proves \eqref{st} since we added superfluous constraint edges of total weight $\w{\smin}$ in the preprocessing step.

		\section{Examples of shapes and ribbons} \label{app:example}
		
	In this section, we provide examples that illustrate the concepts from Section~\ref{subsec:simple}.
	To begin with, we illustrate a shape without middle vertices and its corresponding ribbon.

 \begin{example}[A shape without middle vertices] \label{ex:1}
		Let us choose the ground set $N=[10]$.
		Our first example of a shape is $\alpha$ with $\U{\alpha} = (u_1,u_2)$, $
		\V{\alpha} = (v_1,v_2,v_3)$, and $E(\alpha) = \big\{ \{u_1,u_2\}$, $\{u_2,v_1\}$, $ \{u_2,v_3\}\big\}$.
		The graphical representation of $\alpha$ is provided in Figure~\ref{fig:1a}.
		Notice that we draw boxes around distinguished vertices $V(\U{\alpha})$ and $V(\V{\alpha})$.
		
		This shape $\alpha$ can be realized as a ribbon through a realization.
		For instance, consider a realization $\sigma: \{u_1,u_2,v_1,v_2,v_3,w_1,w_2\} \to [10]$ defined as $\sigma(u_1)=5$, $\sigma(u_2)=1$, $\sigma(v_1)=7$, $\sigma(v_2)=9$, and $\sigma(v_3)=3$ .
		Then, $R=\sigma(\alpha)$ is the ribbon depicted in Figure~\ref{fig:1b}.
		Given this, $\chi_R$ is equal to the product of $\chi_e(G)$ for all  edges $e$ in Figure~\ref{fig:1b}.
	\end{example}
		\begin{figure}[H]
			\centering
			\begin{subfigure}[H]{0.40\textwidth}
				\centering
				\scalebox{0.65}{
					\begin{tikzpicture}[shorten >=1pt,auto,node distance=3cm,
					thick,main node/.style={circle,fill=blue!20,draw,font=\sffamily\LARGE\bfseries}]
					\node[main node] (1) at (-3,3)  {$u_1$};
					\node[main node] (2) at (-3,0)  {$u_2$};
					\node[main node] (3) at (3,3)  {$v_1$};
					\node[main node] (4) at (3,0)  {$v_2$};
					\node[main node] (5) at (3,-3)  {$v_3$}; 
					\node[draw=none,fill=none] (8) at (-3,4.5) {\fontsize{20}{22.4}\selectfont$\U{\alpha}$};
					
					\node[draw=none,fill=none] (10) at (3,4.5) {\fontsize{20}{22.4}\selectfont$\V{\alpha}$};
					\path[every node/.style={font=\sffamily},line width=1pt]
				    (1) edge node [] {} (2)
				    (2) edge node [] {} (3)
				    (2) edge node [] {} (5);
		\draw[draw=black] ($(1.135)+(-0.3,0.3)$) rectangle  ($(2.-45)+(0.3,-0.3)$);
		\draw[draw=black] ($(3.135)+(-0.3,0.3)$) rectangle  ($(5.-45)+(0.3,-0.3)$);
	\end{tikzpicture}
					
				}
				\caption{Figure \ref{fig:1a}:  $\alpha$ from Example~\ref{ex:1}.}
				\label{fig:1a}
			\end{subfigure}
			\qquad
			\begin{subfigure}[H]{0.40\textwidth}
				\centering
				\scalebox{0.65}{
					\begin{tikzpicture}[shorten >=1pt,auto,node distance=3cm,
					thick,main node/.style={circle,fill=blue!20,draw,font=\sffamily\LARGE\bfseries}]
					\node[main node] (1) at (-3,3)  {$5$};
					\node[main node] (2) at (-3,0)  {$1$};
					\node[main node] (3) at (3,3)  {$7$};
					\node[main node] (4) at (3,0)  {$9$};
					\node[main node] (5) at (3,-3)  {$3$};  
					\node[draw=none,fill=none] (8) at (-3,4.5) {\fontsize{20}{22.4}\selectfont$A_R$};
					 
					\node[draw=none,fill=none] (10) at (3,4.5) {\fontsize{20}{22.4}\selectfont$B_R$};
						\path[every node/.style={font=\sffamily},line width=1pt]
				  (1) edge node [] {} (2)
				    (2) edge node [] {} (3)
				    (2) edge node [] {} (5);
				\draw[draw=black] ($(1.135)+(-0.3,0.3)$) rectangle  ($(2.-45)+(0.3,-0.3)$);
		\draw[draw=black] ($(3.135)+(-0.3,0.3)$) rectangle  ($(5.-45)+(0.3,-0.3)$);		\end{tikzpicture}
				}
				\caption{Figure \ref{fig:1b}:
			 $R$ from Example~\ref{ex:1}}
				\label{fig:1b}
			\end{subfigure}
			\caption{A shape and a ribbon from Example~\ref{ex:1}.}
		\end{figure}
\noindent Next, we illustrate a shape with middle vertices and its corresponding ribbon.
 \begin{example}[A shape with middle vertices and its realization] \label{ex:1-1}
		Let us choose the ground set $N=[10]$.
		Our first example of a shape is $\alpha$ with $\U{\alpha} = (u_1,u_2)$, $
		\V{\alpha} = (v_1,v_2,v_3)$, $\W{\alpha} = \{ w_1 \}$, and $E(\alpha) = \big\{ \{u_1,w_1\}$, $\{u_2,w_1\}$,  $ \{v_1,w_1\}$, $\{v_2,w_1\}$, $\{v_3,w_1\}\big\}$.
		The graphical representation of $\alpha$ is provided in Figure~\ref{fig:1-1a}.
		Notice  that we again draw boxes around distinguished vertices $V(\U{\alpha})$ and $V(\V{\alpha})$.
		
		This shape $\alpha$ can be realized as a ribbon through a realization $\sigma$.
		For instance, consider a realization $\sigma: \{u_1,u_2,v_1,v_2,v_3,w_1\} \to [10]$ defined as $\sigma(u_1)=5$, $\sigma(u_2)=1$, $\sigma(v_1)=7$, $\sigma(v_2)=9$, $\sigma(v_3)=3$,  and $\sigma(w_1)=2$.
		Then, $R=\sigma(\alpha)$ is the ribbon depicted in Figure~\ref{fig:1-1b}.
		Note that $\chi_R$ is equal to the product of $\chi_e(G)$ for all  edges $e$ in Figure~\ref{fig:1b}.
	\end{example}
		\begin{figure}[H]
			\centering
			\begin{subfigure}[H]{0.40\textwidth}
				\centering
				\scalebox{0.65}{
					\begin{tikzpicture}[shorten >=1pt,auto,node distance=3cm,
					thick,main node/.style={circle,fill=blue!20,draw,font=\sffamily\LARGE\bfseries}]
					\node[main node] (1) at (-3,3)  {$u_1$};
					\node[main node] (2) at (-3,0)  {$u_2$};
					\node[main node] (3) at (3,3)  {$v_1$};
					\node[main node] (4) at (3,0)  {$v_2$};
					\node[main node] (5) at (3,-3)  {$v_3$};
					\node[main node] (7) at (0,0)  {$w_1$};
					\node[draw=none,fill=none] (8) at (-3,4.5) {\fontsize{20}{22.4}\selectfont$\U{\alpha}$};
					\node[draw=none,fill=none] (9) at (0,4.5) {\fontsize{20}{22.4}\selectfont$\W{\alpha}$};
					\node[draw=none,fill=none] (10) at (3,4.5) {\fontsize{20}{22.4}\selectfont$\V{\alpha}$};
					\path[every node/.style={font=\sffamily},line width=1pt]
				    (1) edge node [] {} (7)
				    (2) edge node [] {} (7)
					(3) edge node [] {} (7)
					(4) edge node [] {} (7)
					(5) edge node [] {} (7);
		\draw[draw=black] ($(1.135)+(-0.3,0.3)$) rectangle  ($(2.-45)+(0.3,-0.3)$);
		\draw[draw=black] ($(3.135)+(-0.3,0.3)$) rectangle  ($(5.-45)+(0.3,-0.3)$);
	\end{tikzpicture}
					
				}
				\caption{Figure \ref{fig:1-1a}:  $\alpha$ from Example~\ref{ex:1-1}.}
				\label{fig:1-1a}
			\end{subfigure}
			\qquad
			\begin{subfigure}[H]{0.40\textwidth}
				\centering
				\scalebox{0.65}{
					\begin{tikzpicture}[shorten >=1pt,auto,node distance=3cm,
					thick,main node/.style={circle,fill=blue!20,draw,font=\sffamily\LARGE\bfseries}]
					\node[main node] (1) at (-3,3)  {$5$};
					\node[main node] (2) at (-3,0)  {$1$};
					\node[main node] (3) at (3,3)  {$7$};
					\node[main node] (4) at (3,0)  {$9$};
					\node[main node] (5) at (3,-3)  {$3$}; 
					\node[main node] (7) at (0,0)  {$2$};
					\node[draw=none,fill=none] (8) at (-3,4.5) {\fontsize{20}{22.4}\selectfont$A_R$};
					\node[draw=none,fill=none] (9) at (0,4.5) {\fontsize{20}{22.4}\selectfont$C_R$};
					\node[draw=none,fill=none] (10) at (3,4.5) {\fontsize{20}{22.4}\selectfont$B_R$};
						\path[every node/.style={font=\sffamily},line width=1pt]
				  (1) edge node [] {} (7)
				    (2) edge node [] {} (7)
					(3) edge node [] {} (7)
					(4) edge node [] {} (7)
					(5) edge node [] {} (7);
				\draw[draw=black] ($(1.135)+(-0.3,0.3)$) rectangle  ($(2.-45)+(0.3,-0.3)$);
		\draw[draw=black] ($(3.135)+(-0.3,0.3)$) rectangle  ($(5.-45)+(0.3,-0.3)$);		\end{tikzpicture}
				}
				\caption{Figure \ref{fig:1-1b}:
			 $R$ from Example~\ref{ex:1-1}}
				\label{fig:1-1b}
			\end{subfigure}
			\caption{A shape and a ribbon from Example~\ref{ex:1-1}.}
		\end{figure}
\noindent Lastly, we illustrate a shape with nonempty left/right intersection.	
 \begin{example} [A shape with left/right intersection]\label{ex:2}
    Note from Definition~\ref{def:shape} that left and right vertices could intersect.
    To illustrate, let us add a single variable $\nn_1$ to both  $V(\U{\alpha})$	and $V(\V{\alpha})$ in Example~\ref{ex:1}, and hence, $\U{\alpha} = (u_1,u_2,\nn_1)$ and $\V{\alpha} = (v_1,v_2,v_3,\nn_1)$, and also an edge $\{u_1,\nn_1\}$ to $E(\alpha)$.
	 Then the resulting shape satisfies $V(\U{\alpha})\cap V(\V{\alpha}) = \nn_1$.
	 The graphical representation of the resulting shape is provided in Figure~\ref{fig:2}.
	 Here and below, we  reserve the notation $\nn$ for such vertices in the left/right intersection.
		\begin{figure}[H]
				\centering
				\scalebox{0.65}{
					\begin{tikzpicture}[shorten >=1pt,auto,node distance=3cm,
					thick,main node/.style={circle,fill=blue!20,draw,font=\sffamily\LARGE\bfseries}]
					\node[main node] (1) at (-3,3)  {$u_1$};
					\node[main node] (2) at (-3,0)  {$u_2$};
					\node[main node] (3) at (3,3)  {$v_1$};
					\node[main node] (4) at (3,0)  {$v_2$};
					\node[main node] (5) at (3,-3)  {$v_3$};
					\node[main node] (7) at (0,0)  {$w_2$};
					\node[draw=none,fill=none] (8) at (-3,4.5) {\fontsize{16}{22.4}\selectfont$\U{\alpha}\setminus \V{\alpha}$};
					\node[draw=none,fill=none] (9) at (0,4.5) {\fontsize{16}{22.4}\selectfont$\U{\alpha}\cap \V{\alpha}$};
					\node[draw=none,fill=none] (10) at (3,4.5) {\fontsize{16}{22.4}\selectfont$\V{\alpha}\setminus \U{\alpha}$};
					\node[draw=none,fill=none] (10) at (0,1.5) {\fontsize{20}{22.4}\selectfont$\W{\alpha}$};
					\node[main node] (11) at (0,3)  {$\nn_1$};
					\path[every node/.style={font=\sffamily},line width=1pt]
			        (1) edge node [] {} (7)
				    (2) edge node [] {} (7)
					(3) edge node [] {} (7)
					(4) edge node [] {} (7)
					(5) edge node [] {} (7)
					(1) edge node [] {} (11);
		\draw[draw=black] ($(1.135)+(-0.3,0.3)$) rectangle  ($(2.-45)+(0.3,-0.3)$);
			\draw[draw=black] ($(11.135)+(-0.3,0.3)$) rectangle  ($(11.-45)+(0.3,-0.3)$);
			\draw[draw=black] ($(11.135)+(-0.4,0.4)$) rectangle  ($(11.-45)+(0.4,-0.4)$);
		\draw[draw=black] ($(3.135)+(-0.3,0.3)$) rectangle  ($(5.-45)+(0.3,-0.3)$);
	\end{tikzpicture}
			 }
			 	\caption{ A shape $\alpha$  with left/right intersection from Example~\ref{ex:2}.}
				\label{fig:2}
		\end{figure}
	\end{example}

				\section{More examples of graph matrices and their norm bounds} \label{app:example2}
	In this section, we provide more examples of graph matrices and their norm bounds.
	\begin{example} \label{ex:4}
	Consider $V(\alpha)=\{u_1,v_1,v_2,w_1\}$, $\U{\alpha}=(u_1)$, $\V{\alpha} =(v_1,v_2)$, $\W{\alpha} = \{w_1\}$ and $E(\alpha)=\big\{\{u_1,w_1\},\{v_1,w_1\},\{v_2,w_1\}\big\}$. 
The graphical representation of $\alpha$ is provided in Figure~\ref{fig:3b}.
    Then, $\M_\alpha$ is an $n\times n(n-1)$ matrix such that for $a,b_1,b_2\in [n]$ ($b_1\neq b_2$), the $(a,(b_1,b_2))$-th entry is defined as
    	\begin{align} \label{ex:4:1}
	\M_\alpha \big(a,(b_1,b_2)\big) &= \begin{cases} 
\displaystyle\sum_{c_1\in [n]\setminus\{a,b_1,b_2\}}\chi_{\{a,c_1\}}\chi_{\{b_1,c_1\}}\chi_{\{b_2,c_1\}}
\,,& a\neq b_1\text{ and }a\neq b_2\,,\\
	0, & \text{otherwise}\,.
	\end{cases}
	\end{align}

Now let us consider the shape $\alpha'$ obtained from $\alpha$ by adding an isolated vertex $w_2$ to $\W{\alpha}$.
The graphical representation of $\alpha'$ is provided in Figure~\ref{fig:3c}.
Then, one can easily verify that for $a,b_1,b_2\in [n]$ ($b_1\neq b_2$), the $(a,(b_1,b_2))$-th entry is defined as:
	\begin{align}
	\M_{\alpha'} \big(a,(b_1,b_2)\big) &= \begin{cases} 
\displaystyle\sum_{c_1\in [n]\setminus\{a,b_1,b_2\}}	\sum_{c_2\in [n]\setminus\{a,b_1,b_2,c_1\}}\chi_{\{a,c_1\}}\chi_{\{b_1,c_1\}}\chi_{\{b_2,c_1\}}
\,,& a\neq b_1\text{ and }a\neq b_2\,,\\
	0, & \text{otherwise}\,.
	\end{cases} \nonumber \\
	&= \begin{cases} 
\displaystyle(n-4)\cdot\sum_{c_1\in [n]\setminus\{a,b_1,b_2\}}\chi_{\{a,c_1\}}\chi_{\{b_1,c_1\}}\chi_{\{b_2,c_1\}}
\,,& a\neq b_1\text{ and }a\neq b_2\,,\\
	0, & \text{otherwise}\,. \label{ex:4:2}
	\end{cases}
	\end{align}
	Notice the similarity between \eqref{ex:4:1} and \eqref{ex:4:2}: indeed, one can conclude  $\M_\alpha = (n-4)\cdot \M_{\alpha '}$.
	This similarity will play a crucial role in reducing the case with isolated middle vertices to the case without them; see Section~\ref{reduc:2}.
	\end{example}

		\begin{figure}[H]
			\centering

			\begin{subfigure}[H]{0.40\textwidth}
				\centering
				\scalebox{0.65}{
					\begin{tikzpicture}[shorten >=1pt,auto,node distance=3cm,
					thick,main node/.style={circle,fill=blue!20,draw,font=\sffamily\LARGE\bfseries}]
				\node[main node] (1) at (-3,0)  {$u_1$};
				\node[main node] (2) at (3,0)  {$v_1$};
				\node[main node] (3) at (3,-2)  {$v_2$};
				\node[main node] (4) at (0,0)  {$w_1$};
					\node[draw=none,fill=none] (5) at (-3,1.5) {\fontsize{20}{22.4}\selectfont$\U{\alpha}$};
					\node[draw=none,fill=none] (6) at (3,1.5) {\fontsize{20}{22.4}\selectfont$\V{\alpha}$};
						\node[draw=none,fill=none] (7) at (0,1.5) {\fontsize{20}{22.4}\selectfont$\W{\alpha}$};
					\path[every node/.style={font=\sffamily},line width=1pt]
				    (1) edge node [] {} (4)
				    (2) edge node [] {} (4)
				     (3) edge node [] {} (4);
		\draw[draw=black] ($(1.135)+(-0.3,0.3)$) rectangle  ($(1.-45)+(0.3,-0.3)$);
		\draw[draw=black] ($(2.135)+(-0.3,0.3)$) rectangle  ($(3.-45)+(0.3,-0.3)$);	\end{tikzpicture}
				}
				\caption{Figure \ref{fig:3b}: Shape $\alpha$ from Example~\ref{ex:4}}
				\label{fig:3b}
			\end{subfigure} \quad\quad
			\begin{subfigure}[H]{0.40\textwidth}
				\centering
				\scalebox{0.65}{
					\begin{tikzpicture}[shorten >=1pt,auto,node distance=3cm,
					thick,main node/.style={circle,fill=blue!20,draw,font=\sffamily\LARGE\bfseries}]
				\node[main node] (1) at (-3,0)  {$u_1$};
				\node[main node] (2) at (3,0)  {$v_1$};
				\node[main node] (3) at (3,-2)  {$v_2$};
				\node[main node] (4) at (0,0)  {$w_1$};
				\node[main node] (5) at (0,-2)  {$w_2$};
					\node[draw=none,fill=none] (5) at (-3,1.5) {\fontsize{20}{22.4}\selectfont$\U{\alpha}$};
					\node[draw=none,fill=none] (6) at (3,1.5) {\fontsize{20}{22.4}\selectfont$\V{\alpha}$};
						\node[draw=none,fill=none] (7) at (0,1.5) {\fontsize{20}{22.4}\selectfont$\W{\alpha}$};
					\path[every node/.style={font=\sffamily},line width=1pt]
				    (1) edge node [] {} (4)
				    (2) edge node [] {} (4)
				     (3) edge node [] {} (4);
		\draw[draw=black] ($(1.135)+(-0.3,0.3)$) rectangle  ($(1.-45)+(0.3,-0.3)$);
		\draw[draw=black] ($(2.135)+(-0.3,0.3)$) rectangle  ($(3.-45)+(0.3,-0.3)$);	\end{tikzpicture}
				}
				\caption{Figure \ref{fig:3c}: Shape $\alpha'$ from Example~\ref{ex:4}}
				\label{fig:3c}
			\end{subfigure}
			\caption{Shapes from Examples~\ref{ex:4}.}
		\end{figure}

	Now we illustrate  Theorem~\ref{thm:mainresult} for the graph matrix presented in Example~\ref{ex:4}.
		
			\begin{example}[Revisiting Example~\ref{ex:4}]\label{ex:bd:2}
	For  shape $\alpha$ from Example~\ref{ex:4}, $\{w_1\}$ is a minimum size vertex separator, which implies
	$s= 1$.
By Theorem~\ref{thm:mainresult}, with high probability, $\norm{\M_\alpha}\leq \widetilde{O}(n^{\frac{3}{2}})$. 
Now let us consider  shape $\alpha'$. 
For $\alpha'$, there is one isolated middle vertex, and consequently by the second statement of Theorem~\ref{thm:mainresult}, $\norm{\M_\alpha}\leq \widetilde{O}(n^{\frac{3}{2}+1})=\widetilde{O}(n^{\frac{5}{2}})$  with high probability.
This is consistent with the observation we made in Example~\ref{ex:4}: $\M_{\alpha'} = (n-4)\cdot \M_{\alpha}$
	\end{example}
	
\noindent We also illustrate Theorem~\ref{thm:mainresult} for the graph matrices corresponding to the shapes from Appendix~\ref{app:example}.
	
	\begin{example}[Revisiting Examples~\ref{ex:1}, \ref{ex:1-1} and \ref{ex:2}]\label{ex:bd:3}
	For the shape $\alpha$ from Example~\ref{ex:1}, $\{u_2\}$ is a minimum size vertex separator, so
	$s= 1$.
	By Theorem~\ref{thm:mainresult}, $\norm{\M_\alpha}\leq \widetilde{O}(n^{2})$ with high probability.
	Next, for the shape $\alpha$ from Example~\ref{ex:1-1}, $\{w_1 \}$ is a minimum size vertex separator, so $\norm{\M_\alpha}\leq \widetilde{O}(n^{5/2})$ with high probability.
	Lastly, for the shape $\alpha$ from Example~\ref{ex:2}, $\{w_1, \nn_1 \}$ is a minimum size vertex separator, and hence, 
	$\norm{\M_\alpha}\leq \widetilde{O}(n^{5/2})$ with high probability.
	Notice that the graph matrices corresponding to shapes from Example \ref{ex:1-1} and \ref{ex:2} have the  norm bound of the same order.
	This similarity is not conincidental, and will be formalized in our proof via the reduction argument; see Section~\ref{reduc:1}.
	\end{example}
	
	\bibliographystyle{alpha}
	\bibliography{references}
	 
	 	\section*{Acknowledgements}
	
	This research was initially supported by the Program for Research in Mathematics, Engineering, and Science (PRIMES) at MIT. 
	We thank the PRIMES faculty, including Dr. Tanya Khovanova, Dr. Pavel Etingof, and Dr. Slava Gerovitch for helpful feedback during the initial stage of this work. 
	We also thank Laci Babai, Fernando Granha Jeronimo, and Chris Jones for helpful comments on this manuscript.
	
\end{document}